\documentclass[11pt,a4paper]{paper}
\usepackage[utf8]{inputenc}
\usepackage[margin=2.5cm]{geometry}

\usepackage{amsmath,amssymb,amsfonts}
\usepackage[tight]{subfigure}
\usepackage{color}
\usepackage[utf8]{inputenc}
\usepackage{enumitem}
\usepackage{relsize}
\usepackage[draft]{fixme}
\usepackage{tikz}
\usepackage{tensor}

\newcounter{thm}
\newcounter{ex}
\newtheorem{Theorem}[thm]{Theorem}

\newtheorem{Definition}[thm]{Definition}

\newcommand{\RR}{\mathbb{R}}

\newcommand{\Cc}{{\mathcal C}}
\newcommand{\ip}[1]{\left<#1\right>}
\newcommand{\Log}{\mathrm{Log}}
\newcommand{\Id}{\mathrm{Id}}
\let\oldphi\phi \let\phi\varphi \let\varphi\oldphi
\tikzset{node distance=1.75cm, auto}

\newcommand{\der}[1]{\partial_{#1}}
\usepackage{tensor}

\title{
  Anisotropically Weighted and Nonholonomically Constrained Evolutions on
  Manifolds
}

\author{Stefan Sommer}

\begin{document}
\maketitle

\begin{abstract}
We present evolution equations for a family of paths that results from anisotropically weighting curve energies in non-linear statistics of manifold valued data. This situation arises when performing inference on data that have non-trivial covariance and are anisotropic distributed. The family can be interpreted as most probable paths for a driving semi-martingale that through stochastic development is mapped to the manifold. We discuss how the paths are projections of geodesics for a sub-Riemannian metric on the frame bundle of the manifold, and how the curvature of the underlying connection appears in the sub-Riemannian Hamilton-Jacobi equations. Evolution equations for both metric and cometric formulations of the sub-Riemannian metric are derived. We furthermore show how rank-deficient metrics can be mixed with an underlying Riemannian metric, and we relate the paths to geodesics and polynomials in Riemannian geometry. Examples from the family of paths are visualized on embedded surfaces, and we explore computational representations on finite dimensional landmark manifolds with geometry induced from right-invariant metrics on diffeomorphism groups.
\end{abstract}

\section{Introduction}
When manifold valued data have non-trivial covariance, i.e. when
\emph{anisotropy} asserts higher variance in some directions than others,
non-zero curvature necessitates special care when
generalizing Euclidean space normal distributions to manifold valued
distributions: In the Euclidean situation, normal distributions can be seen as
transition distributions of
diffusion processes but, on the manifold, holonomy makes transport 
of covariance path-dependent in the presence of curvature preventing a global 
notion of spatially constant covariance matrix.
To handle this, in the diffusion PCA framework \cite{sommer_diffusion_2014} and with
the class of anisotropic normal distributions on manifolds defined in
\cite{sommer_anisotropic_2015,sommer_modelling_2016},
data on non-linear manifolds are modelled as being distributed according 
to transition distributions of anisotropic diffusion processes that are
mapped from Euclidean space to the manifold by stochastic
development (see \cite{hsu_stochastic_2002}). The construction is connected to a
non bracket-generating sub-Riemannian metric on the bundle of linear frames of the
manifold, the frame bundle, and the requirement that covariance stays
covariantly constant gives a nonholonomically constrained system.

Velocity vectors and geodesic distances are conventionally used for
estimation and statistics in Riemannian manifolds, for example for estimating the
Frech\'et mean \cite{frechet_les_1948}, for
Principal Geodesic Analysis \cite{fletcher_principal_2004-1}, and for tangent space
statistics \cite{vaillant_statistics_2004}. In contrast to
this, anisotropy as modelled with anisotropic normal distributions 
makes a distance for a sub-Riemannian metric the natural vehicle for 
estimation and statistics. This metric naturally accounts for 
anisotropy in a similar way as the 
precision matrix weights the inner product in the negative log-likelihood of a 
Euclidean normal distribution.
The connection between the weighted distance and statistics
of manifold valued data was
presented in \cite{sommer_anisotropic_2015}, and the underlying
sub-Riemannian and fiber-bundle geometry together with properties of the generated
densities was further explored in \cite{sommer_modelling_2016}. The fundamental
idea is to perform statistics on manifolds by maximum likelihood (ML) instead of
parametric constructions that use e.g. approximating geodesic subspaces: by defining
natural families of probability distributions, in this case using diffusion
processes, ML parameter estimates gives a coherent way to statistically model
non-linear data. The
anisotropically weighted distance and the resulting family of extremal paths
arises in this situation when the diffusion processes have non-isotropic 
covariance, i.e. when the
distribution is not generated from a standard Brownian motion.

In this paper, we focus on the family of \emph{most probable paths} for the 
semi-martingales that drives the stochastic development and in turn the manifold
valued anisotropic stochastic processes. These paths extremize the anisotropically 
weighted action functional.
We present derivations of evolution equations 
for the paths from different viewpoints, and we discuss the role of frames as 
representing either metrics or cometrics. In the derivation, we explicitly see
the influence of the connection and its curvature. We then turn to the 
relation between the 
sub-Riemannian metric and the Sasaki-Mok metric on the frame bundle,
and we develop a construction that allows the sub-Riemannian metric
to be defined as a sum of
a rank-deficient generator and an underlying Riemannian metric. 
Finally, we relate the paths to geodesics and polynomials in Riemannian geometry, and we
explore computational representations on different manifolds including 
a specific case, the finite dimensional manifolds arising in
the LDDMM \cite{younes_shapes_2010} landmark matching problem. The paper ends
with a discussion concerning statistical implications, open questions, and
concluding remarks.
\begin{figure}[t]
  \begin{center}
    \includegraphics[width=.44\columnwidth,trim=170 120 150 220,clip]{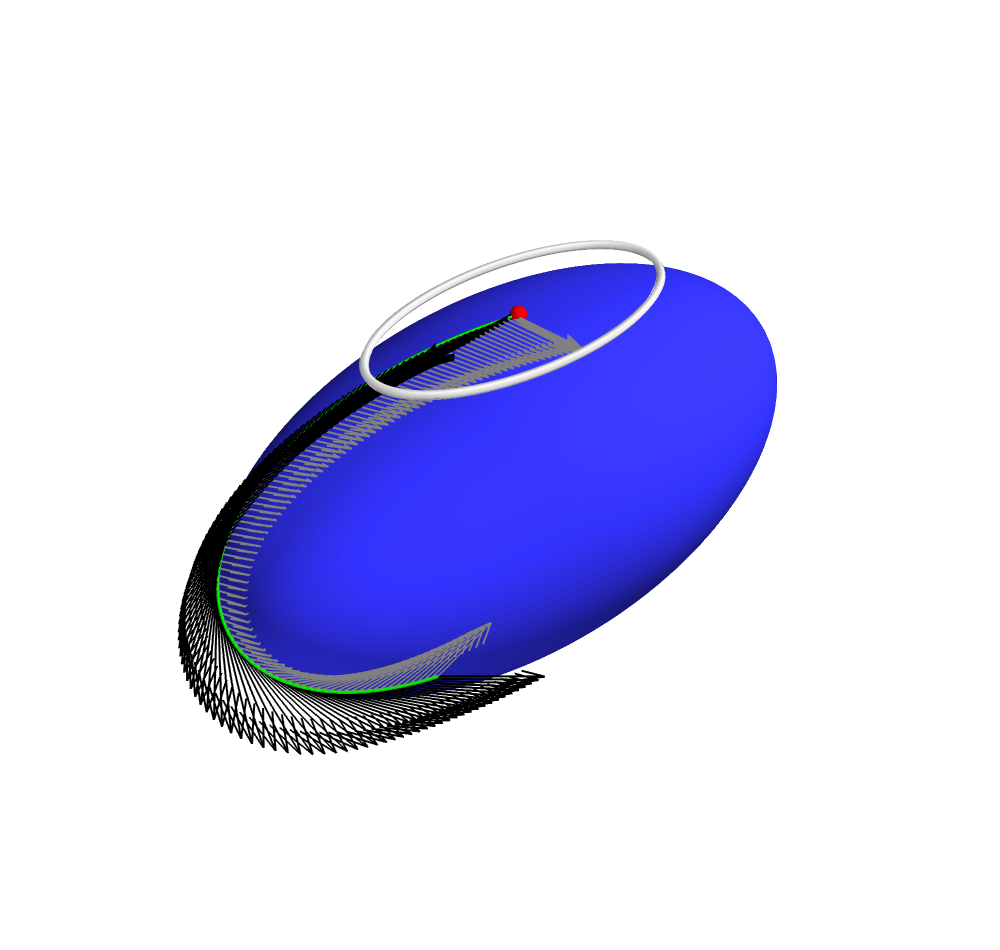}
    \includegraphics[width=.49\columnwidth,trim=100 0 100 50,clip]{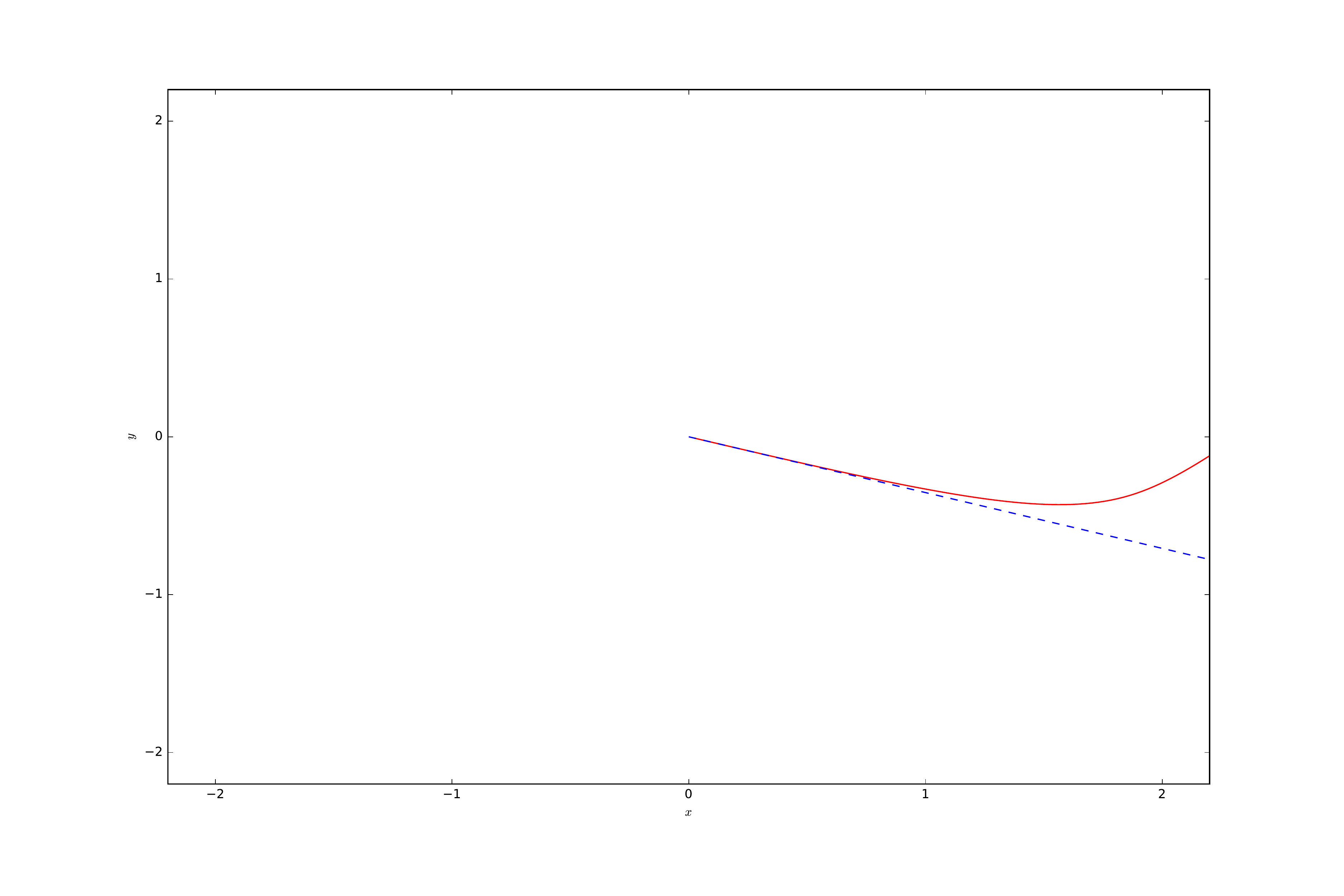}
  \end{center}
  \caption{
    (left) A \emph{most probable path} (MPP) for a driving Euclidean Brownian motion
    on an ellipsoid. The gray ellipsis over the starting point (red dot) indicates
    the covariance of the anisotropic diffusion. A
    frame $u_t$ (black/gray vectors) representing the square root covariance 
    is parallel transported along the curve enabling the
    anisotropic weighting with the precision matrix in the action functional. With isotropic
    covariance, normal MPPs are Riemannian geodesics. In general situations such as the
    displayed anisotropic case, the family of MPPs is much larger. (right) The
    corresponding anti-development in $\RR^2$ (red line) of the MPP. Compare with the
    anti-development of a Riemannian geodesic with same initial velocity (blue
    dotted line). The frames
    $u_t\in\mathrm{GL}(\RR^2,T_{x_t}M)$ provide local frame coordinates for each time
    $t$.
    }
  \label{fig:subrgeo}
\end{figure}

\subsection{Background}
Generalizing common statistical tools for performing inference on
Euclidean space data to manifold valued data has been the subject of extensive
work, see e.g. \cite{pennec_intrinsic_2006}. Perhaps most
fundamental is the notion of Frech\'et or Karcher means
\cite{frechet_les_1948,karcher_riemannian_1977} defined as minimizers of the
square Riemannian distance. Generalizations of the Euclidean 
Principal Component Analysis (PCA) procedure to manifolds are particularly
relevant for data exhibiting anisotropy. Approaches include 
Principal Geodesic Analysis (PGA,
\cite{fletcher_principal_2004-1}), Geodesic PCA (GPCA,
\cite{huckemann_intrinsic_2010}), Principal Nested Spheres (PNS,
\cite{jung_analysis_2012}), Barycentric
Subspace Analysis (BSA, \cite{pennec_barycentric_2015}),
and Horizontal Component Analysis (HCA, \cite{sommer_horizontal_2013}).
Common to these constructions are explicit representations of approximating low-dimensional 
subspaces. The fundamental
challenge is here that the notion of Euclidean linear subspace on which PCA
relies has no direct analogue in non-linear spaces.

A different approach taken by Diffusion PCA 
(DPCA, \cite{sommer_diffusion_2014,sommer_anisotropic_2015}) 
and Probabilistic PGA \cite{zhang_probabilistic_2013} is to base the PCA problem
on a maximum likelihood fit of normal distributions to data. In Euclidean space,
this approach was first introduced with Probabilistic PCA \cite{tipping_probabilistic_1999}.
In DPCA, the process of stochastic development
\cite{hsu_stochastic_2002} is used to define a class of anisotropic distributions 
that generalizes the family of Euclidean space normal distributions to the 
manifold context. DPCA is then a simple 
maximum likelihood fit in this family of distributions mimicking the Euclidean
Probabilistic PCA. The approach transfers the
geometric complexities of defining subspaces common in the approaches listed
above to the problem of defining a geometrically natural notion of normal
distributions. 
\begin{figure}[t]
  \begin{center}
\mbox{
\def\svgwidth{0.68\columnwidth}
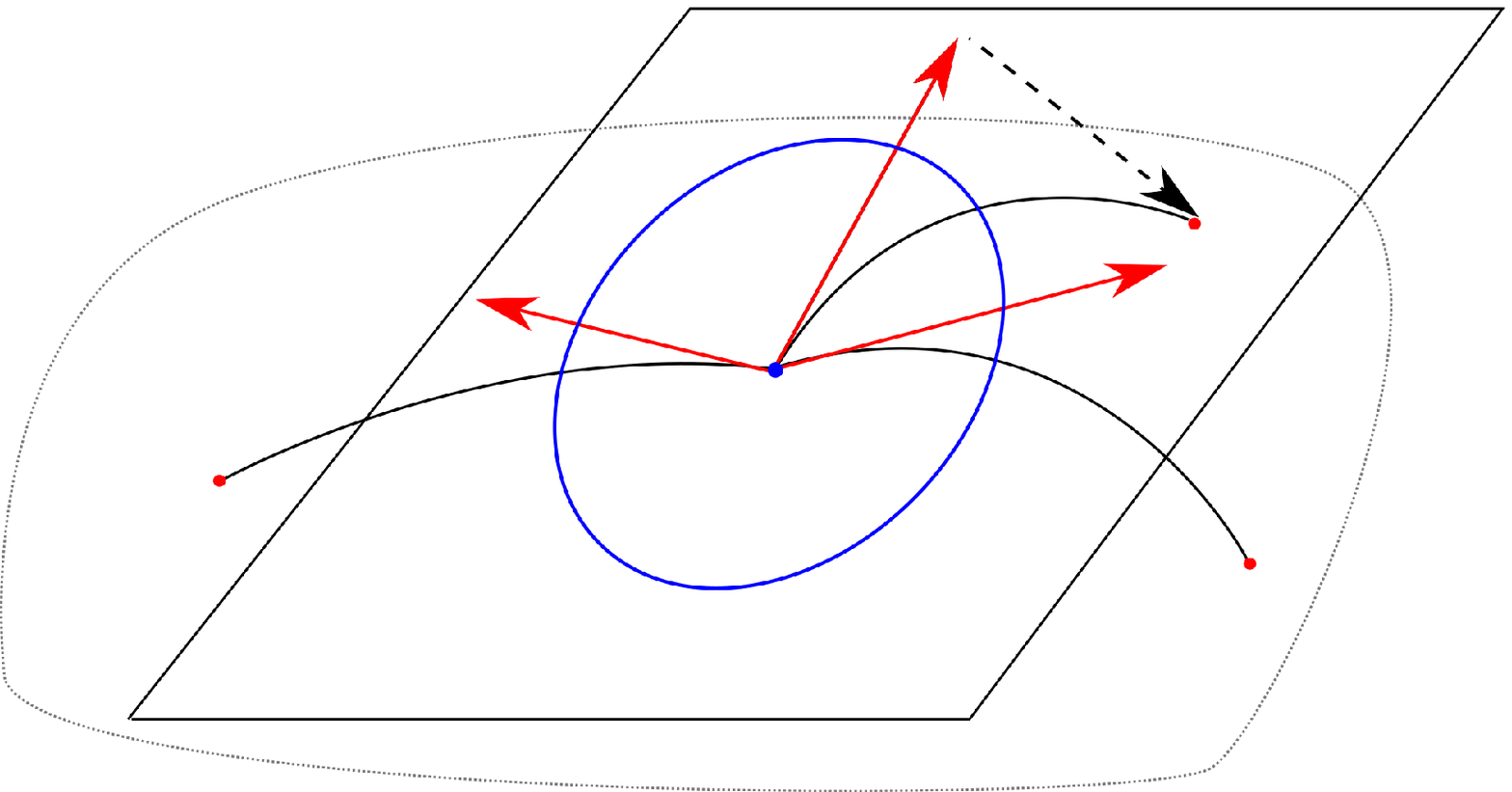
\hspace{-2.9cm}
}
\mbox{
\def\svgwidth{0.64\columnwidth}
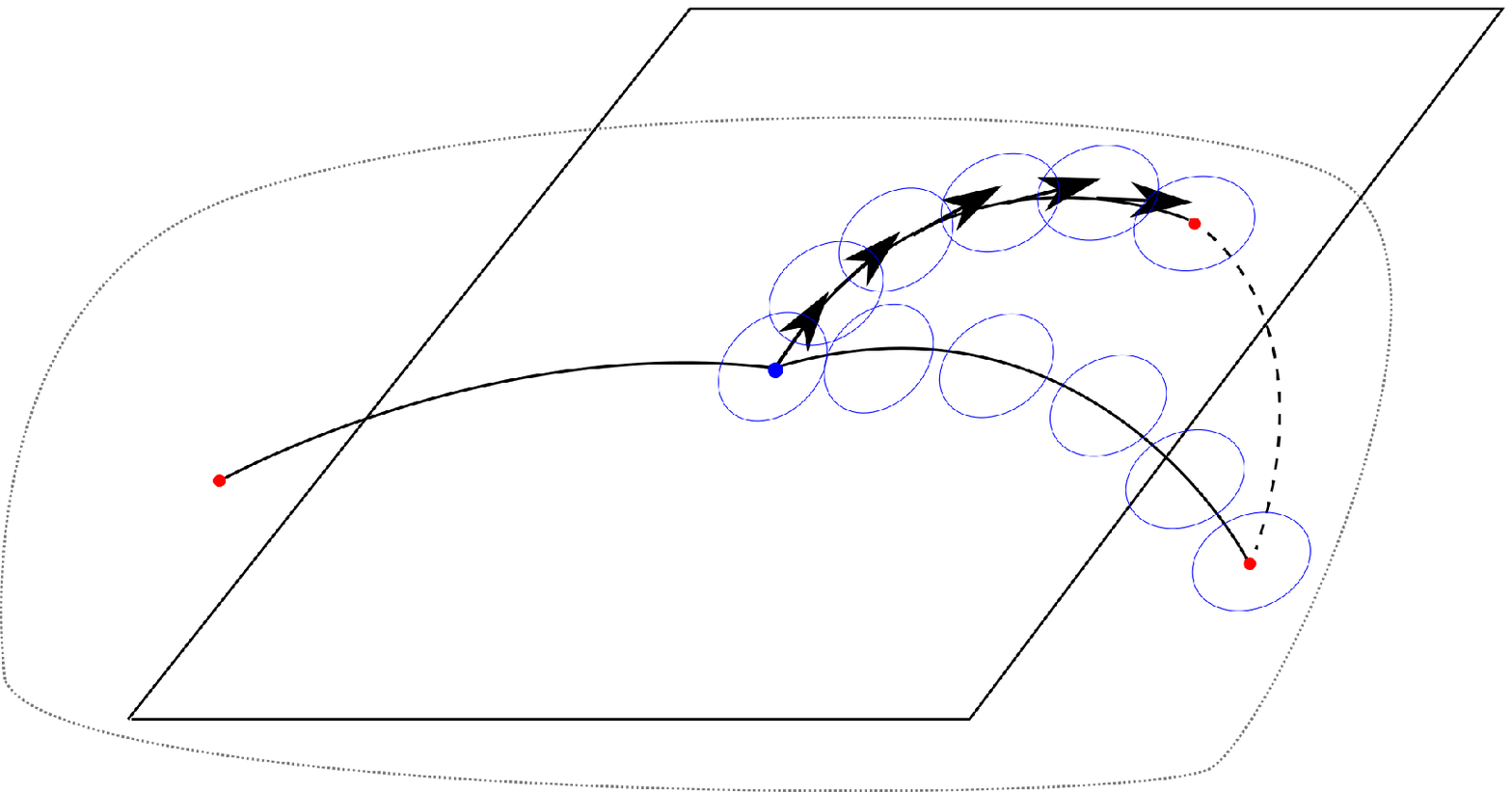
\hspace{-2.9cm}
}
    
  \end{center}
  \caption{
    (left) Normal distributions $u\mathcal N(0,\Id)$ in the tangent space
    $T_{x_0}M$ with covariance $uu^T$ (blue ellipsis) can be mapped to 
    the manifold by applying the exponential map
    $\mathrm{Exp}_{x_0}$ to sampled vectors $v\in T_{x_0}M$ (red vectors). This 
    effectively linearises 
    the geometry around $x_0$. (right)
    The stochastic development map $\phi_u$ maps $\RR^d$ valued paths $w_t$ to $M$
    by transporting the covariance in each step (blue
    ellipses) giving a covariance $u_t$ along the entire sample path. The approach does not linearise around a
    single point. Holonomy of the connection implies that the covariance
    ``rotates'' around closed loops, an effect which can be illustrated by continuing 
    the transport along
    the loop created by the dashed path. The anisotropic metric $g_{FM}$ weights
    step lengths by the transported covariance at each time $t$.
  }
  \label{fig:linearization}
\end{figure}

In Euclidean space, squared distances $\|x-x_0\|^2$ between observations $x$
and the mean $x_0$ are affinely related to the negative log-likelihood of
a normal distribution $\mathcal N(x_0,\Id)$. This makes an ML fit of the mean
such as performed in Probabilistic PCA equivalent to minimizing squared
distances. On a manifold, distances $d_g(x,x_0)^2$ coming from a Riemannian
metric $g$ are equivalent to tangent space distances $\|\Log_{x_0}x\|^2$
when mapping data from $M$ to $T_{x_0}M$ using the inverse exponential map
$\Log_{x_0}$. Assuming $\Log_{x_0}x$ are distributed according to a normal
distribution in the linear space $T_{x_0}M$, this restores the equivalence with a
maximum likelihood fit. Let $\{e_1,\ldots,e_d\}$ be the standard basis for
$\RR^d$. If $u:\RR^d\rightarrow T_{x_0}M$ is a linear invertible map with $ue_1,\ldots,ue_d$
orthonormal with respect to $g$, the normal distribution in $T_{x_0}M$ can be defined 
as $u\mathcal N(0,\Id)$, see Figure~\ref{fig:linearization}.

The map $u$ can be represented as a point in the frame bundle $FM$ of $M$. When
the orthonormal requirement on $u$ is relaxed so that $u\mathcal N(0,\Id)$
is a normal distribution in $T_{x_0}M$ with anisotropic covariance,
the negative log-likelihood in $T_{x_0}M$ is related to
$(u^{-1}\Log_{x_0}x)^T(u^{-1}\Log_{x_0}x)$ in the
same way as the precision matrix $\Sigma^{-1}$ is related to the negative
log-likelihood $(x-x_0)^T\Sigma^{-1}(x-x_0)$ in Euclidean space. The distance is
thus weighted by the anisotropy of $u$, and $u$ can be interpreted as a square
root covariance matrix $\Sigma^{1/2}$.

However, the above approach does not
specify how $u$ changes when moving away from the base point $x_0$.
The use of $\Log_{x_0}x$ effectively linearises the geometry around $x_0$ but a
geometrically natural way to relate $u$ at points nearby to $x_0$ will be to
parallel transport it, equivalently specify that $u$ when transported does 
not change as measured
from the curved geometry. This constraint is
\emph{nonholonomic} and it implies that any path from
$x_0$ to $x$ carries with it a parallel transport of $u$ lifting paths from $M$
to paths in the frame bundle $FM$. It therefore becomes natural to equip
$FM$ with a form of metric that encodes the anisotropy represented by $u$. The result 
is the sub-Riemannian metric on $FM$ defined below that
weights infinitesimal movements on $M$ using the parallel transport of the frame
$u$. Optimal paths for this metric are sub-Riemannian geodesics giving the
family of \emph{most probable paths} for the driving process that this paper
concerns. Figure~\ref{fig:subrgeo} shows one such path for an
anisotropic normal distribution with $M$ an ellipsoid embedded in $\RR^3$.

\section{Frame Bundles, Stochastic Development, and Anisotropic Diffusions}
Let $M$ be a finite dimensional manifold of dimension $d$ with connection $\Cc$, and
let $x_0$ be a fixed point in $M$. When a Riemannian metric is present and $\Cc$
is its Levi-Civita connection, we
denote the metric $g_R$. For a given interval $[0,T]$, we let $W(M)$
denote the Wiener space of continuous paths in $M$ starting at $x_0$. Similarly,
$W(\RR^d)$ is the Wiener space of paths in $\RR^d$. We let $H(\RR^d)$ denote the
subspace of $W(\RR^d)$ of finite energy paths.

Let now $u=(u_1,\ldots,u_d)$ be a frame for $T_xM$, $x\in M$, i.e.
$u_1,\ldots,u_d$ is an ordered set of linearly independent vectors in $T_xM$ with
$\mathrm{span}\{u_1,\ldots,u_d\}=T_xM$. We can regard the frame as an isomorphism 
$u: \mathbb R^d \to T_xM$ with $u(e_i) = u_i$ where $e_1,\ldots,e_d$ denotes the
standard basis in $\mathbb R^d$. Stochastic development, see e.g.
\cite{hsu_stochastic_2002}, provides an invertible map
$\phi_u$ from $W(\RR^d)$ to $W(M)$. Through $\phi_u$, 
Euclidean semi-martingales map to
stochastic processes on $M$. When $M$ is Riemannian and $u$ orthonormal,
the result is the
Eells-Elworthy-Malliavin construction of Brownian motion
\cite{elworthy_geometric_1988}. We here outline
the geometry behind development, stochastic development, the connection and
curvature focusing in particular on frame bundle geometry. 

\subsection{The Frame Bundle}
For each point $x\in M$, let $F_xM$ be the set of frames for $T_xM$,
i.e. the set of ordered bases for $T_xM$. 
The set $\{F_xM\}_{x\in M}$ can be given a natural differential structure as a
fiber bundle on $M$ called the
frame bundle $FM$. It can equivalently be defined as the principal bundle
$\mathrm{GL}(\RR^d,TM)$. We let the map $\pi:FM\rightarrow M$ denote the
canonical projection. The kernel of $\pi_*:TFM\rightarrow TM$ is the subbundle of
$TFM$ that consists of
vectors tangent to the fibers $\pi^{-1}(x)$. It is denoted the
vertical subspace $VFM$.
We will often work in a local trivialization $u=(x,u_1,\ldots,u_d)\in FM$ where
$x=\pi(u)\in M$ denotes the base point and for each $i=1,\ldots,d$,
$u_i\in T_xM$ is the $i$th frame vector. For $v\in T_xM$ and $u\in
FM$ with $\pi(u)=x$, the vector $u^{-1}v\in\RR^d$ expresses $v$ in components in
terms of the frame $u$. We will denote the vector $u^{-1}v$ frame coordinates of $v$.

For a differentiable curve $x_t$ in $M$ with $x=x_0$,
a frame $u$ for $T_{x_0}M$ can be parallel transported along 
$x_t$ by parallel transporting each vector in the frame thus giving a path
$u_t\in FM$.
Such paths are called horizontal and have zero acceleration in the sense
$\Cc(\dot{u}_{i,t})=0$. 
For each $x\in M$, their derivatives form a $d$-dimensional subspace of the
$d+d^2$-dimensional tangent space $T_uFM$. This horizontal subspace
$HFM$ and the vertical subspace $VFM$ 
together split the tangent bundle of $FM$, i.e.
$TFM=HFM\oplus VFM$. 
The split induces a map
$\pi_*:HFM\rightarrow TM$,
see Figure~\ref{fig:FM}. For fixed $u\in FM$,
the restriction
$\pi_*|_{H_uFM}:H_uFM\rightarrow T_{x}M$ is an isomorphism. Its inverse
is called the horizontal lift and denoted $h_u$ in the following.
Using $h_u$, horizontal vector fields $H_e$ on $FM$ are defined for
vectors $e\in\RR^d$ by $H_e(u)=h_u(ue)$. In particular, 
the standard basis $(e_1,\ldots,e_d)$ on $\RR^d$ gives
$d$ globally defined horizontal vector fields $H_i\in HFM$, $i=1,\ldots,d$ by
$H_i=H_{e_i}$. Intuitively, the fields $H_i(u)$ model infinitesimal transformations in $M$ of
$x_0$ in direction $u_i=ue_i$ with corresponding infinitesimal parallel transport of 
the vectors $u_1,\ldots,u_d$ of the frame along the direction $u_i$.
A \emph{horizontal lift} of a differentiable curve
$x_t\in M$ is a curve in $FM$ tangent to $HFM$ that projects to $x_t$.
Horizontal lifts are unique up to the choice of initial frame $u_0$.

\subsection{Development and Stochastic Development}
Let $x_t$ be a differentiable curve on $M$ and $u_t$ a horizontal lift. If
$s_t$ is a curve in $\RR^d$ with components $s_t^i$ such that $\dot{x}_t=H_i(u)s_t^i$, $x_t$
is said to be a development of $s_t$. Correspondingly, $s_t$ is the
anti-development of $x_t$. For each $t$, the vector $s_t$ contains frame coordinates of $\dot{x}_t$ as defined above. 
Similarly, let $W_t$ be an $\RR^d$ valued Brownian motion 
so that sample paths $W_t(\omega)\in W(\RR^d)$. A solution to the stochastic differential 
equation
$
  dU_t
  =
  \sum_{i=1}^d
  H_i(U_t)
  \circ
  dW^i_t
$
in $FM$ is called a stochastic development of $W_t$ in $FM$. The solution projects to
a stochastic development $X_t=\pi(U_t)$ in $M$. We call the
process $W_t$ in $\RR^d$, that through $\phi$ maps to $X_t$, the driving
process of $X_t$. Let $\phi_u:W(\RR^d)\rightarrow W(M)$ be the map that for
fixed $u$ sends a path in $\RR^d$ to its development on $M$. Its inverse $\phi_u^{-1}$
is the anti-development in $\RR^d$ of paths on $M$ given $u$.

Equivalent to the fact that normal
distributions $\mathcal{N}(0,\Sigma)$ in $\RR^d$ can be obtained as the transition
distributions of diffusion processes $\Sigma^{1/2}W_t$ stopped at time $t=1$, 
a general class of
distributions on the manifold $M$ can be defined by stochastic development of
processes $W_t$ resulting in $M$-valued random variables $X=X_1$. This family of
distributions on $M$ introduced in \cite{sommer_anisotropic_2015} is denoted 
\emph{anisotropic normal distributions}. The stochastic development by
construction ensures that $U_t$ is horizontal, and the frames
are thus parallel transported along the stochastic displacements. The effect is that
the frames stay covariantly constant thus resembling the Euclidean situation
where $\Sigma^{1/2}$ is spatially constant and therefore does not change as 
$W_t$ evolves. Thus, as further discussed in Section~\ref{sec:covnonhol},
the covariance is kept constant at each of the infinitesimal stochastic
displacements. The existence of a smooth density for the target process
$X_t$ and small time asymptotics are discussed in \cite{sommer_modelling_2016}.
\begin{figure}[t]
\begin{center}
\begin{tikzpicture}
  \node (TFM) {$TFM$};
  \node (TSFM) [below of=TFM] {$T^*FM$};
  \node (empty) [right of=TFM] {};
  \node (HFM) [right of=empty] {$HFM$};
  \node (empty2) [left of=TFM] {};
  \node (VFM) [left of=empty2] {$VFM$};
  \node (FMlie) [below of=VFM] {$FM\times\mathfrak{gl}(n)$};
  \node (FM) [right of=HFM] {$FM$};
  \node (TM) [below of=HFM] {$TM$};
  \node (TSM) [left of=TM] {$T^*M$};
  \node (M) [right of=TM] {$M$};
  \draw[->] (TFM) to node {$h+v\mapsto h$} (HFM);
  \draw[->] (TFM) to node [swap] {$h+v\mapsto v$} (VFM);
  \draw[->] (HFM) to node {} (FM);
  \draw[->] (HFM) to node [swap] {$\pi_*$} (TM);
  \draw[<->] (VFM) to node [swap] {$\psi$} (FMlie);
  \draw[->] (FM) to node {$\pi_{FM}$} (M);
  \draw[->] (TM) to node {$\pi_{TM}$} (M);
  \draw[->] (TSFM) to node [swap] {$g_{FM}$} (HFM);
  \draw[<->] (TSM) to node [swap] {$g_R$} (TM);
\end{tikzpicture}
\end{center}
\caption{Relations between the manifold, frame bundle, the horizontal
  distribution $HFM$, the vertical bundle $VFM$, a Riemannian metric $g_R$ and 
  the sub-Riemannian metric $g_{FM}$ defined below.
  The connection $\Cc$ provides the splitting $TFM=HFM\oplus VFM$.
  The restrictions $\pi_*|_{H_uM}$ are invertible maps
  $H_uM\rightarrow T_{\pi(u)}M$ with inverse $h_u$, the horizontal lift.
  Correspondingly, the vertical bundle $VFM$ is isomorphic to the trivial bundle
  $FM\times\mathfrak{gl}(n)$.  The metric $g_{FM}:T^*FM\rightarrow TFM$ has image in the
  subspace $HFM$.
}
\label{fig:FM}
\end{figure}

Stochastic development gives a map $\int_{\mathrm{Diff}}:FM\rightarrow\mathrm{Prob}(M)$ 
to the space of probability distributions on $M$. For each
point $u\in FM$, the map sends a
Brownian motion in $\RR^d$ to a distribution $\mu_u$ by stochastic development of
the process $U_t$ in $FM$ starting at $u$ and letting $\mu_u$ be the
distribution of $X=\pi(U_1)$.
The pair $(x,u)$, $x=\pi(u)$ is analogous to the parameters
$(\mu,\Sigma)$ for a Euclidean normal distribution: The point $x\in M$
represents the starting  point of the diffusion, and the frame $u$ represents a
square root $\Sigma^{1/2}$ of the covariance $\Sigma$. In the general situation 
where $\mu_u$ has smooth density, the construction can be used to fit the 
parameters $u$ to data
by maximum likelihood. As an example,
Diffusion PCA fits distributions obtained through $\int_{\mathrm{Diff}}$ by 
maximum likelihood to observed samples in $M$, i.e. it optimizes for the most
likely parameters $u=(x,u_1,\ldots,u_d)$ for the anisotropic diffusion process giving
a fit to the data of the manifold generalization of the Euclidean normal
distribution.

\subsection{Adapted Coordinates}
\label{sec:coords}
For concrete expressions of the geometric constructions related to frame bundles
and for computational purposes, it is useful to apply coordinates that are
adapted to the horizontal bundle $HFM$ and the vertical bundle $VFM$ together
with their duals $H^*FM$ and $V^*FM$. The notation below follows the notation
used in for example
\cite{mok_differential_1978}. Let
$z=(u,\xi)$ be a local trivialization of $T^*FM$, and let 
$(x^i,u^i_\alpha)$ be coordinates on $FM$
with $u^i_\alpha$ satisfying $u_\alpha=u^i_\alpha\der{x^i}$ for each
$\alpha=1,\ldots,d$.

To find a basis that is adapted to the horizontal distribution, define the
$d$ linearly independent vector fields
$
  D_j
  =
  \der{x^j}
  -
  \Gamma^{h_\gamma}_j\der{u^h_\gamma}
$
where
$
  \Gamma^{h_\gamma}_j
  =
  \Gamma\indices{^h_{ji}}u^i_\gamma
$
is the contraction of the Christoffel symbols $\Gamma\indices{^h_{ij}}$ for the connection
$\Cc$ with $u_\alpha^i$. We denote this adapted frame $D$. The vertical
distribution is correspondingly spanned by 
$D_{j_\beta}=\partial_{u^j_\beta}$.
The vectors 
$
  D^h=dx^h,
  $
and
$
  D^{h_\gamma}
  =
  \Gamma^{h_\gamma}_jdx^j+du^h_\gamma
  $
constitutes a dual coframe $D^*$.
The map $\pi_*:HFM\rightarrow TM$ is in coordinates of the adapted frame
$\pi_*(w^jD_j)=w^j\der{x^j}$. Correspondingly, the horizontal lift $h_u$ is
$h_u(w^j\der{x^j})=w^jD_j$.
The map $u:\RR^d\rightarrow T_xM$ is given by
the matrix $[u^i_\alpha]$ so that 
$uv=u^i_\alpha v^\alpha\der{x^i}=u_\alpha v^\alpha$. 

Switching between standard coordinates and the adapted frame and coframes can
be expressed in terms of the component matrices $A$ below of
the frame and coframe induced by the coordinates $(x^i,u_\alpha^i)$ and the adapted
frame $D$ and coframe $D^*$. We have
\begin{equation*}
  \,_{(\der{x^i},\der{u_\alpha^i})}A_D
  =
  \begin{bmatrix}
    I & 0 \\ -\Gamma & I
  \end{bmatrix}
  \mbox{ with inverse }
  \,_DA_{(\der{x^i},\der{u_\alpha^i})}
  =
  \begin{bmatrix}
    I & 0 \\ \Gamma & I
  \end{bmatrix}
\end{equation*}
writing $\Gamma$ for the matrix $[\Gamma^{h_\gamma}_j]$. Similarly, the component
matrices of the dual frame $D^*$ are
\begin{equation*}
  \,_{(\der{x^i},\der{u_\alpha^i})^*}A_{D^*}
  =
  \begin{bmatrix}
    I & \Gamma^T \\ 0 & I
  \end{bmatrix}
  \mbox{ and }
  \,_{D^*}A_{(\der{x^i},\der{u_\alpha^i})^*}
  =
  \begin{bmatrix}
    I & -\Gamma^T \\ 0 & I
  \end{bmatrix}
  \ .
\end{equation*}

\subsection{Connection and Curvature}
The $TM$ valued connection $\Cc:TM\times TM\rightarrow TM$ lifts to 
a principal connection $TFM\times TFM\rightarrow VFM$ on the principal bundle $FM$. 
$\Cc$ can then be identified with the $\mathfrak{gl}(n)$-valued connection
form $\omega$ on $TFM$. The identification
occurs by the isomorphism $\psi$ between $FM\times\mathfrak{gl}(n)$ and $VFM$
given by $\psi(u,v)=\tfrac{d}{dt}u\exp(tv)|_{t=0}$, see e.g.
\cite{taubes_differential_2011,kolar_natural_1993}. 

The map $\psi$ is equivariant with
respect to the $\mathrm{GL}(n)$ action $g\mapsto ug^{-1}$ on $FM$. In order to
explicitly see the connection between the usual covariant derivative
$\nabla:\Gamma(TM)\times\Gamma(TM)\rightarrow\Gamma(TM)$ on $M$ determined by
$\Cc$ and $\Cc$ regarded as a connection on the principal bundle $FM$, 
following \cite{taubes_differential_2011}, we let
$s:M\rightarrow TM$ be a local vector field on $M$, equivalently 
$s\in \Gamma(TM)$ is a local section of $TM$. $s$ determines a map
$s^{FM}:FM\rightarrow\RR^d$ by $s^{FM}(u)=u^{-1}s(\pi(u))$, i.e. it gives the
coordinates of $s(x)$ in the frame $u$ at $x$. The pushforward
$(s^{FM})_*:TFM\rightarrow\RR^d$ has in its $i$th component the exterior
derivative $d(s^{FM})^i$. Let now $w(x)$ be a local section of $FM$.
The composition 
$w\circ(s^{FM})_*\circ h_w:TM\rightarrow TM$ is identical to the covariant
derivative $\nabla_\cdot s:TM\rightarrow TM$. The
construction is independent of the choice of $w$ because of the
$\mathrm{GL}(n)$-equivariance of $s^{FM}$. The connection form $\omega$ can be
expressed as the matrix $(s_1^{FM}\circ h_w,\ldots,s_d^{FM}\circ h_w)$ when
letting $s_i^{FM}(u)=e_i$.

The identification becomes particularly simple if the covariant derivative is
taken along a curve $x_t$ on which $w_t$ is the horizontal lift. In this case,
we can let $s_t=w_{t,i}s^i_t$. Then $s^{FM}(w_t)=(s^1_t,\ldots,s^d_t)^T$ and
\begin{equation}
  w_t^{-1}\nabla_{\dot x_t}s
  =
  (s^FM)_*(h_{w_t}(\dot x_t))=\tfrac{d}{dt}
  (s^1_t,\ldots,s^d_t)^T
  \ ,
  \label{eq:framecovder}
\end{equation}
i.e. the covariant derivative takes the form of the standard derivative applied
to the frame coordinates $s^i_t$.
  
The curvature tensor $R\in \mathcal T^3_1(M)$ gives the
$\mathfrak{gl}(n)$-valued curvature form 
$\Omega:TFM\times TFM\rightarrow \mathfrak{gl}(n)$ on $TFM$ by
\begin{equation*}
  \Omega(v_u,w_u)
  =
  u^{-1}R(\pi_*(v_u),\pi_*(w_u))u
  \ ,\ 
  v_u,w_v\in TFM
  \ .
\end{equation*}
Note that $\Omega(v_u,w_u)=\Omega(h_u(\pi_*(v_u)),h_u(\pi_*(w_u)))$
which we can use to write the curvature $R$ as the $\mathfrak{gl}(n)$-valued map
$R_u:T^2(T_{\pi(u)}M)\rightarrow \mathfrak{gl}(n)$,
$(v,w)\mapsto\Omega(h_u(\pi_*(v_u)),h_u(\pi_*(w_u)))$ for fixed $u$.
In coordinates, the curvature is 
\begin{equation*}
  R\indices{_{ijk}^s}
  =
  \Gamma\indices{^l_{ik}}
  \Gamma\indices{^s_{jl}}
  -
  \Gamma\indices{^l_{jk}}
  \Gamma\indices{^s_{il}}
  +
  \Gamma\indices{^s_{ik;j}}
  -
  \Gamma\indices{^s_{jk;i}}
\end{equation*}
where $\Gamma\indices{^s_{ik;j}}=\der{x^j}\Gamma\indices{^s_{ik}}$.

Let $x_{t,s}$ be a family of paths in $M$ and let $u_{t,s}\in\pi^{-1}(x_{t,s})$ 
be horizontal lifts of $x_{t,s}$ for each fixed $s$.
Write $\dot{x}_{t,s}=\partial_t x_{t,s}$ and $\dot{u}_{t,s}=\partial_t u_{t,s}$.
The $s$-derivative of $u_{t,s}$ can be regarded a pushforward of the horizontal lift 
and is the curve in $TFM$
\begin{equation}
  \begin{split}
  \partial_su_{t,s}
  &=
  \psi\big(u_{t,s},
  \psi_{u_{0,s}}^{-1}(\Cc(\partial_su_{0,s}))
  +\int_0^s\Omega(\dot{u}_{r,s},\partial_s u_{r,s})dr\big)
  +
  h_{u_{t,s}}(\partial_s x_{t,s})
  \\&
  =
  \psi\big(u_{t,s},
  \psi_{0,s}^{-1}(\Cc(\partial_su_{0,s}))
  +\int_0^s R_{u_{r,s}}(\dot{x}_{r,s},\partial_s x_{r,s}) dr\big)
  +
  h_{u_{t,s}}(\partial_s x_{t,s})
  \ .
  \end{split}
  \label{eq:hdercurve}
\end{equation}
This follows from the structure equation 
  $d\omega
  =
  -\omega
  \wedge
  \omega
  +
  \Omega
  $,
see e.g. \cite{andersson_finite_1999}. Note that the curve depends on the
vertical variation $\Cc(\partial_su_{0,s})$ at only one point along the curve.
The remaining terms depend on the horizontal variation or, equivalently,
$\partial_sx_{t,s}$.
The $t$-derivative of $\partial_su_{t,s}$
is the curve in $TTFM$ satisfying
\begin{equation}
  \begin{split}
    \partial_sh_{u_{t,s}}(\dot{x}_{t,s})
  &
  =
  \psi\big(u_{t,s},
  R_{u_{t,s}}(\dot{x}_{t,s},\partial_s x_{t,s})\big)
  +
  \partial_t\psi\big(u_{t,s},
  \psi_{0,s}^{-1}(\Cc(\partial_su_{0,s}))
  \big)
  +
  \partial_t\big(h_{u_{t,s}}(\partial_s x_{t,s})\big)
  \\
  &
  =
  \psi\big(u_{t,s},
  R_{u_{t,s}}(\dot{x}_{t,s},\partial_s x_{t,s})\big)
  +
  \partial_t\psi\big(u_{t,s},
  \psi_{0,s}^{-1}(\Cc(\partial_su_{0,s}))
  \big)
  \\&\qquad
  +
  h_{u_{t,s}}(\partial_t\partial_s x_{t,s})
  +
  (\partial_th_{u_{t,s}})(\partial_s x_{t,s})
  \ .
  \end{split}
  \label{eq:hder}
\end{equation}
Here, the first and third term in the last expression are identified with elements of
$T_{\partial_su_{t,s}}TFM$ by the natural mapping 
$T_{u_{t,s}}FM\rightarrow T_{\partial_su_{t,s}}TFM$.
When $\Cc(\partial_su_{0,s})$ is zero,
the relation reflects the property that the
curvature arise when computing brackets between horizontal vector fields.
Note that the first term of \eqref{eq:hder} has values in $VFM$ while the third
term has values in $HFM$.

\section{The Anisotropically Weighted Metric}
For a Euclidean driftless diffusion process with 
spatially constant stochastic generator $\Sigma$, the 
log-probability of a sample path can formally be written
\begin{equation}
  \ln \tilde{p}_\Sigma(x_t)
    \propto
    -
    \int_0^1
    \|\dot{x}_t\|_{\Sigma}^2
    dt
    +
    c_\Sigma
    \label{eq:pathdens}
\end{equation}
with the norm $\|\cdot\|_{\Sigma}$ given by the inner product
$\ip{v,w}_{\Sigma}=\ip{\Sigma^{-1/2}v,\Sigma^{-1/2}w}=v\Sigma^{-1}w$, i.e. the inner product
weighted by the precision matrix $\Sigma^{-1}$.
Though only formal as the sample paths are almost surely nowhere differentiable,
the interpretation can be given a precise meaning by
taking limits of piecewise linear curves \cite{andersson_finite_1999}.
Turning to the manifold situation with the processes mapped to $M$ by stochastic
development, the probability of observing a differentiable path can either be
given a precise meaning in the
manifold by taking limits of small tubes around the curve, or in $\RR^d$ by
considering infinitesimal tubes around the anti-development of the curves. With the 
former formulation,
a scalar curvature correction term must be added to \eqref{eq:pathdens} giving
the Onsager-Machlup function \cite{fujita_onsager-machlup_1982}. The latter
formulation corresponds to defining a notion of path density for the driving
$\RR^d$-valued process $W_t$.
When $M$ is Riemannian and $\Sigma$ unitary, taking the maximum of 
\eqref{eq:pathdens} gives geodesics
as most probable paths for the driving process.

Let now $u_t$ be a path in $FM$ and choose a local trivialization
$u_t=(x_t,u_{1,t},\ldots,u_{d,t})$ such that the matrix
$[u_{\alpha,t}^i]$ represents the square root covariance matrix $\Sigma^{1/2}$ at $x_t$. 
Since $u_t$ being a
frame defines an invertible map $\RR^d\rightarrow T_{x_t}M$, the norm
$\|\cdot\|_{\Sigma}$ above has a direct analogue in the norm
$\|\cdot\|_{u_t}$ defined by the inner product
\begin{equation}
  \ip{v,w}_{u_t}
  = 
  \ip{
    u_t^{-1}v
    ,
    u_t^{-1}w
  }_{\RR^d}
  \label{eq:mppdens}
\end{equation}
for vectors $v,w\in T_{x_t}M$. The transport of the frame along paths in effect defines a
transport of inner product along sample paths: The paths carry with them the
inner product weighted by the precision matrix which in turn is a transport of the square 
root covariance $u_0$ at $x_0$. 

The inner product can equivalently be defined as a metric
$g_u:T^*_{x}M\rightarrow T_{x}M$. Again using that $u$ can be
considered a map $\RR^d\rightarrow T_xM$, $g_u$ is defined by 
$\xi\mapsto u_((\xi\circ u)^\sharp)$ where $\sharp$ is the standard identification
$(\RR^d)^*\rightarrow\RR^d$. The sequence of mappings defining $g_u$ is illustrated
below:
\begin{equation}
  \begin{matrix}
    T^*_xM &\rightarrow& (\RR^d)^* &\rightarrow& \RR^d &\rightarrow& T_xM \\
    \xi \quad&\mapsto&\quad \xi\circ u \quad&\mapsto&\quad (\xi\circ u)^\sharp
    \quad&\mapsto&\quad u(\xi\circ u)^\sharp 
  \ .
  \end{matrix}
  \label{eq:g}
\end{equation}
This definition uses the $\RR^d$ inner product in
the definition of $\sharp$.
Its inverse gives the cometric
$g_u^{-1}:T_xM\rightarrow T_x^*M$, i.e. $v\mapsto (u^{-1}v)^\flat\circ
u^{-1}$.
\begin{equation}
  \begin{matrix}
    T_xM &\rightarrow& \RR^d &\rightarrow& (\RR^d)^* &\rightarrow& T_x^*M \\
    v &\mapsto& u^{-1}v &\mapsto& (u^{-1})^\flat &\mapsto&
    (u^{-1})^\flat \circ u^{-1}
  \ .
  \end{matrix}
  \label{eq:ginv}
\end{equation}

\subsection{Sub-Riemannian Metric on the Horizontal Distribution}
We now lift the path-dependent metric defined above to a sub-Riemannian metric
on $HFM$. For any $w,\tilde{w}\in H_uFM$, the lift of \eqref{eq:mppdens}
by $\pi_*$ is the inner product 
\begin{equation*}
  \ip{w,\tilde{w}}
  = 
  \ip{
    u^{-1}\pi_*w,
    u^{-1}\pi_*\tilde{w}
  }_{\RR^d}
  \ .
\end{equation*}
The inner product induces a sub-Riemannian metric 
$g_{FM}:TFM^*\rightarrow HFM\subset TFM$ by
\begin{equation}
  \ip{w,g_{FM}(\xi)}
  =
  (\xi|w)
  \ ,\ 
  \forall w\in H_uFM
  \label{eq:subrieg}
\end{equation}
with $(\xi|w)$ denoting the evaluation $\xi(w)$ for the covector $\xi\in T^*FM$. The 
metric $g_{FM}$ gives $FM$ a
non bracket-generating sub-Riemannian structure
\cite{strichartz_sub-riemannian_1986} on $FM$,
see also Figure~\ref{fig:FM}.
It is equivalent to the lift 
\begin{equation}
  \xi\mapsto h_u(g_u(\xi\circ h_u))\ ,\ \xi\in T_uFM
  \label{eq:lift}
\end{equation}
of the metric $g_u$ above. In frame coordinates, the metric 
takes the form
\begin{equation}
  u^{-1}\pi_*g_{FM}(\xi)
  =
  \begin{pmatrix}
    \xi(H_1(u)) \\
    \vdots \\
    \xi(H_d(u))
  \end{pmatrix}
  \ .
  \label{eq:metricframecomp}
\end{equation}
In terms of the adapted coordinates for $TFM$ described in Section~\ref{sec:coords}, with
$w=w^jD_j$ and $ \tilde{w}=\tilde{w}^jD_j$, we have
\begin{align*}
  \ip{w,\tilde{w}}
  &=
  \ip{w^iD_i,\tilde{w}^jD_j}
  =
  \ip{u^{-1}w^i\der{x^i},u^{-1}\tilde{w}^j\der{x^j}}
  \\
  &=
  \ip{w^iu^\alpha_i,\tilde{w}^ju^\alpha_j}_{\RR^d}
  =
  \delta_{\alpha\beta}w^iu^\alpha_i\tilde{w}^ju^\beta_j
  =
  W_{ij}w^i\tilde{w}^j
\end{align*}
where $[u^\alpha_i]$ is the inverse of $[u^i_\alpha]$ and
$W_{ij}=\delta_{\alpha\beta}u^\alpha_iu^\beta_j$. Define
now $W^{kl}=\delta^{\alpha\beta}u^k_\alpha u^l_\beta$ so that 
$W^{ir}W_{rj}=\delta^i_j$ and $W_{ir}W^{rj}=\delta_i^j$.
We can then write the metric $g_{FM}$ directly as 
\begin{equation}
  g_{FM}(\xi_hD^h+\xi_{h_\gamma}D^{h_\gamma})
  =
  W^{ih}\xi_hD_i
  \ ,
  \label{eq:gadapted}
\end{equation}
because 
$
  \ip{w,g_{FM}(\xi)}
  =
  \ip{w,W^{jh}\xi_hD_j}
  =
  W_{ij}w^iW^{jh}\xi_h
  =w^i\xi_i
  =\xi_hD^h(w^jD_j)
  =\xi(w)
$.
One clearly recognizes the dependence on the horizontal $H^*FM$ part of $T^*FM$ only and
the fact that $g_{FM}$ has image in $HFM$.
The sub-Riemannian energy of an a.e. horizontal path $u_t$ is
\begin{equation*}
  l_{FM}(u_t)
  =
  \int g_{FM}(\dot{u}_t,\dot{u}_t)dt
  \ ,
\end{equation*}
i.e., the line element is $ds^2=W_{ij}D^iD^j$ in adapted coordinates. The
corresponding distance is given by
\begin{equation*}
  d_{FM}(u_1,u_2) = \inf\{l_{FM}(\gamma) \mid \gamma(0)=u_1, \gamma(1)=u_2\}
\ . 
\end{equation*}

If we wish to express $g_{FM}$ in canonical coordinates
on $T^*FM$, we can switch between the
adapted frame and the coordinates $(x^i,u_\alpha^i,\xi^i,\xi_\alpha^i)$.
From \eqref{eq:gadapted}, $g_{FM}$ has $D,D^*$ components
\begin{equation*}
  \,_Dg_{FM,D^*}
  =
  \begin{bmatrix}
    W^{-1} & 0 \\
    0 & 0
  \end{bmatrix}
  \ .
\end{equation*}
Therefore, $g_{FM}$ has the following components in the coordinates
$(x^i,u_\alpha^i,\xi_h,\xi_{h_\gamma})$
\begin{equation*}
  \,_{(\der{x^i},\der{u_\alpha^i})}g_{FM,(\der{x},\der{u_\alpha^i})^*}
  =
  \,_{(\der{x^i},\der{u_\alpha^i})}A_D\,_Dg_{FM,D^*}\,_{D^*}A_{(\der{x^i},\der{u_\alpha^i})^*}
  =
  \begin{bmatrix}
    W^{-1} & -W^{-1}\Gamma^T  \\
    -\Gamma W^{-1} & \Gamma W^{-1}\Gamma^T
  \end{bmatrix}
\end{equation*}
or
$
  g_{FM}^{ij}=W^{ij}
  $,\ \ 
  $
  g_{FM}^{i j_\beta}=-W^{ih}\Gamma^{j_\beta}_h
  $,\ \ 
  $
  g_{FM}^{i_\alpha j}=-\Gamma^{i_\alpha}_hW^{hj}
  $, and
  $
  g_{FM}^{i_\alpha j_\beta}=\Gamma^{i_\alpha}_kW^{kh}\Gamma^{j_\beta}_h
  $.

\subsection{Covariance and Nonholonomicity}
\label{sec:covnonhol}
The metric $g_{FM}$ encodes the anisotropic weighting given the
frame $u$ thus up to an affine transformation measuring the energy of horizontal 
paths equivalently to
the negative log-probability of sample paths of Euclidean anisotropic
diffusions as formally given in \eqref{eq:pathdens}. In addition, the
requirement that paths must stay almost everywhere horizontal enforces that
$\Cc(\dot{u}_t)=0$ a.e., i.e. that \emph{no change of the covariance is measured
by the connection}. The intuitive effect is that covariance is covariantly
constant as seen by the connection. Globally, \emph{curvature} of $\Cc$ will imply that
the covariance changes when transported along closed loops, and \emph{torsion}
will imply that the base point ``slips'' when travelling along covariantly closed
loops on $M$. However, the zero acceleration requirement implies that the
covariance is as close to spatially constant as possible with the given
connection. This is enabled by the parallel transport of the frame, and it
ensures that the model closely resembles the Euclidean case with spatially constant 
stochastic generator.

With non-zero curvature of $\Cc$, the horizontal distribution is non-integrable,
i.e. the brackets $[H_i,H_j]$ are non-zero for some $i,j$. This prevents
integrability of the horizontal distribution $HFM$ in the sense of the 
Frobenius theorem. In this case,
the horizontal constraint is \emph{nonholonomic} similarly to nonholonomic
constraints appearing in geometric mechanics, see e.g. \cite{bloch_nonholonomic_2003}.
The requirement of covariantly constant covariance thus results in a nonholonomic
system.

\subsection{Riemannian Metrics on $FM$}
If the horizontality constraint is relaxed, a related Riemannian metric on
$FM$ can be defined by pulling back a metric on $\mathfrak{gl}(n)$ to each
fiber using the isomorphism
$\psi(u,\cdot)^{-1}:V_uFM\rightarrow\mathfrak{gl}(n)$.
Therefore, the metric on $HFM$ can be extended to a Riemannian metric on $FM$. 
Such metrics incorporate the anisotropically weighted metric on $HFM$
however allowing vertical variations and thus that covariances can change
unrestricted.

When $M$ is Riemannian, the metric $g_{FM}$ is in addition related to the Sasaki-Mok
metric on $FM$ \cite{mok_differential_1978} that extends the Sasaki metric on
$TM$. As for the above Riemannian metric on $FM$, the Sasaki-Mok metric allows paths in 
$FM $ to have derivatives in the vertical space $VFM$. On $HFM$, the Riemannian metric
$g_R$ is here lifted to the metric
$g_{SM}=(v_u,w_u)=g_R(\pi_*(v_u),\pi_*(w_u))$, i.e. the metric is not 
anisotropically weighted. The line element is in this case
$ds^2=g_{ij}dx^idx^j+X_{\beta\alpha}g_{ij}D^{\alpha_i}D^{\beta_j}$.

Geodesics for $g_{SM}$ are lifts of Riemannian geodesics for $g_R$ on $M$ in
contrast to the sub-Riemannian normal geodesics for $g_{FM}$ we will
characterize below. The family of curves arising as projections to $M$ of normal 
geodesics for $g_{FM}$ includes Riemannian geodesics for $g_R$ and thus
projections of geodesics for $g_{SM}$ but the family is in general larger than
geodesics for $g_R$.

\section{Constrained Evolutions}
Extremal paths for \eqref{eq:mppdens} can be interpreted as most probable paths for
the driving process $W_t$ when $u_0$ defines an anisotropic diffusion.
This is captured in the following definition \cite{sommer_modelling_2016}:

\begin{Definition}
  \label{def:mpp}
A \emph{most probable path for the driving process} (MPP) from $x=\pi(u_0)\in
M$ to $y\in M$ is a smooth path $x_t: [0,1]\to M $ with $x_0=x$ and
$x_1=y$ such that its anti-development $\phi^{-1}_{u_0}(x_t)$ is a
most probable path for $W_t$, i.e.
\begin{equation*}
x_t\in\textrm{argmin}_{\sigma, \sigma_0=x, \sigma_1=y} \int_0^1
-L_{\RR^d}(\phi_{u_0}^{-1}(\sigma_t),\tfrac{d}{dt}\phi_{u_0}^{-1}(\sigma_t))\, dt
\end{equation*}
with $L_{\RR^d}$ being the
Onsager-Machlup function for the process $W_t$ on $\RR^d$
\cite{fujita_onsager-machlup_1982}.
\end{Definition}
The definition uses the one-to-one relation between $W(\RR^d)$ and $W(M)$ provided by
$\phi_{u_0}$ to characterize the paths using the $\RR^d$ Onsager-Machlup function $L_{\RR^d}$.
When $M$ is Riemannian with metric $g_R$, the Onsager-Machlup function for a
$g$-Brownian motion on $M$ is 
$L(x_t,\dot{x}_t)=-\frac{1}{2}\|\dot x_t\|_{g_R}^2+\tfrac{1}{12}S_{g_R}(x_t)$
with $S_{g_R}$ denoting the scalar curvature. This curvature term vanishes on $\RR^d$
and therefore
$L_{\RR^d}(\gamma_t,\dot{\gamma}_t)=-\frac{1}{2}\|\dot\gamma_t\|^2$
for a curve $\gamma_t\in\RR^d$.

By pulling $x_t\in M$ back to $\RR^d$ using $\phi^{-1}_{u_0}$, the construction removes 
the $\tfrac{1}{12}S_{g_R}(x_t)$ scalar curvature correction term present
in the non-Euclidean Onsager-Machlup function. It thereby provides a relation between
geodesic energy and most probable paths for the driving process. This is
contained in the following characterization of most
probable paths for the driving process as extremal paths of the sub-Riemannian
distance \cite{sommer_modelling_2016} that follows from the Euclidean 
space Onsager-Machlup theorem \cite{fujita_onsager-machlup_1982}.
\begin{Theorem}[\cite{sommer_modelling_2016}]\label{driving}
  Let $Q(u_0)$ denote the principal subbundle of $FM$ of points $z\in FM$
  reachable form $u_0\in FM$ by horizontal paths. Suppose
the H\"{o}rmander condition is satisfied on $Q(u_0)$ and that $Q(u_0)$ has compact fibers. Then most probable paths from $x_0$ to $y\in M$ for the driving process of $X_t$ exist, and they are projections of sub-Riemannian geodesics in $FM$ minimizing the sub-Riemannian distance from $u_0$ to $\pi^{-1}(y)$.
\end{Theorem}
Below, we will derive evolution equations for the set of such extremal paths that 
correspond to normal sub-Riemannian geodesics.

\subsection{Normal Geodesics for $g_{FM}$}
Connected to the metric $g_{FM}$ is the Hamiltonian
\begin{equation}
  H(z)
  =
  \frac{1}{2}
  (z|g_{FM}(z))
  \label{eq:hamiltonian}
\end{equation}
on the symplectic space $T^*FM$. Letting $\hat{\pi}$ denote the 
projection on the bundle $T^*FM\rightarrow FM$, \eqref{eq:subrieg} gives
\begin{equation*}
  H(z)
  =
  \frac{1}{2}\ip{g_{FM}(z)|g_{FM}(z)}
  =
  \frac{1}{2}
  \|z\circ h_{\hat{\pi}(z)}\circ \hat{\pi}(z)\|_{(\RR^d)^*}^2
  =
  \frac{1}{2}
  \sum_{i=1}^d\xi(H_i(u))^2
  \ .
\end{equation*}

Normal geodesics in sub-Riemannian manifolds satisfy the
Hamilton-Jacobi equations \cite{strichartz_sub-riemannian_1986} with
Hamiltonian flow
\begin{equation}
  \dot{z}_t
  =
  X_H
  =
  \Omega^\#dH(z)
  \label{eq:hameq}
\end{equation}
where $\Omega$ here is the canonical symplectic form on $T^*FM$, see e.g.
\cite{marsden_introduction_1999}. We denote \eqref{eq:hameq} the MPP equations,
and we let projections $x_t=\pi_{T^*FM}(z_t)$ of minimizing curves
satisfying \eqref{eq:hameq} be denoted normal MPPs. The system
\eqref{eq:hameq} has $2(d+d^2)$ degrees of freedom in contrast to the usual $2d$
degrees of freedom for the classical geodesic equation. Of these, $d^2$
describes the current frame at time $t$ while the remaining $d^2$ allows the
curve to ``twist'' while still being horizontal. We will see this effect
visualized in Section~\ref{sec:experiments}.

In a local canonical trivialization
$z=(u,\xi)$, \eqref{eq:hameq} gives the Hamilton-Jacobi equations
\begin{equation}
  \begin{split}
  \dot{u}
  &=
  \partial_\xi H(u,\xi) 
  =
  g_{FM}(u,\xi)
  =
  h_u\big(
  u
  (\,\xi(H_1(u)),\ldots,\xi(H_d(u))\,)^T
  \big)
  \\
  \dot{\xi}
  &=
  -\partial_u H(u,\xi)
  =
  -\partial_u 
  \frac{1}{2}
  \|\xi\circ h_u\circ u\|_{(\RR^d)^*}^2
  =
  -\partial_u 
  \frac{1}{2}
  \sum_{i=1}^d\xi(H_i(u))^2
  \ .
\end{split}
\label{eq:hamloctriv}
\end{equation}
Using \eqref{eq:hder}, we have for the second equation
\begin{equation}
  \begin{split}
  \dot{\xi}
  &=
  -
  \sum_{i=1}^d
  \xi(H_i(u))\xi(\partial_uh_u(ue_i))
  \\
  &=
  -
  \sum_{i=1}^d
  \xi(H_i(u))\xi\big(
  \psi(u,R_u(ue_i,\pi_*(\der{u})))
  +\partial_{h_u(ue_i)}\psi\big(u,
  \psi^{-1}(\Cc(\partial_u))
  \big)
  +\partial_{h_u(ue_i)}h_u(\pi_*(\der{u}))
  \big)
  \\
  &=
  -
  \xi\big(
  \psi(u,R_u(\pi_*(\dot{u}),\pi_*(\der{u})))
  +\partial_{\dot{u}}\psi\big(u,
  \psi^{-1}(\Cc(\partial_u))
  \big)
  +\partial_{\dot{u}}h_u(\pi_*(\der{u}))
  \big)
  \ .
  \end{split}
  \label{eq:dotxi}
\end{equation}
Here $\partial_{\dot{u}}$ denotes $u$-derivative in the direction $\dot{u}$,
equivalently
$\partial_{\dot{u}}h_u(v)=\partial_t(h_u)(v)$.
While the first equation of \eqref{eq:hamloctriv} involves only the horizontal 
part of $\xi$, the second equation couples the vertical part of
$\xi$ through the evaluation of $\xi$ on the term
$\psi(u,R_u(\pi_*(\dot{u}),\pi_*(\der{u}))$.
If the connection is curvature-free which in non-flat cases implies
that it carries torsion, this vertical term vanishes.
Conversely, when $M$ is Riemannian, $\Cc$ the $g_R$ Levi-Civita connection, and 
$u_0$ is $g_R$ orthonormal,
$g_{FM}(h_u(v),h_u(w))=g_R(v,w)$ for all $v,w\in T_{\pi(u_t)}M$. In this case,
a normal MPP $\pi(u_t)$ will be a Riemannian $g_R$ geodesic.

\subsection{Evolution in Coordinates}
In coordinates $u=(x^i,u_\alpha^i,\xi_i,\xi_{i_\alpha})$
for $T^*FM$, we can equivalently write
\begin{align*}
  &
  \dot{x}^i
  =
  g^{ij}\xi_j
  +
  g^{ij_\beta}\xi_{j_\beta}
  =
  W^{ij}\xi_j
  -
  W^{ih}\Gamma^{j_\beta}_h\xi_{j_\beta}
  \\&
  \dot{X}^i_\alpha
  =
  g^{i_\alpha j}\xi_j
  +
  g^{i_\alpha j_\beta}\xi_{j_\beta}
  =
  -\Gamma^{i_\alpha}_hW^{hj}\xi_j
  +
  \Gamma^{i_\alpha}_kW^{kh}\Gamma^{j_\beta}_h\xi_{j_\beta}
  \\&
  \dot{\xi}_i
  =
  -\frac{1}{2}
  \left(
  \der{y^i}g_y^{hk}
  \xi_h\xi_k
  +
  \der{y^i}g_y^{h k_\delta}
  \xi_h\xi_{k_\delta}
  +
  \der{y^i}g_y^{h_\gamma k}
  \xi_{h_\gamma}\xi_k
  +
  \der{y^i}g_y^{h_\gamma k_\delta}
  \xi_{h_\gamma}\xi_{k_\delta}
  \right)
  \\&
  \dot{\xi}_{i_\alpha}
  =
  -\frac{1}{2}
  \left(
  \der{y^{i_\alpha}}g_y^{hk}
  \xi_h\xi_k
  +
  \der{y^{i_\alpha}}g_y^{h k_\delta}
  \xi_h\xi_{k_\delta}
  +
  \der{y^{i_\alpha}}g_y^{h_\gamma k}
  \xi_{h_\gamma}\xi_k
  +
  \der{y^{i_\alpha}}g_y^{h_\gamma k_\delta}
  \xi_{h_\gamma}\xi_{k_\delta}
  \right)
  \\&\qquad
\end{align*}
with 
$\Gamma^{h_\gamma}_{k,i}$ for $\der{y^i}\Gamma^{h_\gamma}_k$
and where
\begin{align*}
  &\der{y^{l}}g^{ij}=0
  \ ,\quad
  \der{y^{l}}g^{i j_\beta}=-W^{ih}\Gamma^{j_\beta}_{h,l}
  \ ,\quad
  \der{y^{l}}g^{i_\alpha j}=-\Gamma^{i_\alpha}_{h,l}W^{hj}
  \ ,\quad
  \der{y^{l}}g^{i_\alpha j_\beta}
  =\Gamma^{i_\alpha}_{k,l}W^{kh}\Gamma^{j_\beta}_h+\Gamma^{i_\alpha}_kW^{kh}\Gamma^{j_\beta}_{h,l}
  \ ,\\
  &\der{y^{l_\zeta}}g^{ij}=W\indices{^{ij}_{,l_\zeta}}
  \ ,\quad
  \der{y^{l_\zeta}}g^{i j_\beta}
  =
  -W\indices{^{ih}_{,l_\zeta}}\Gamma^{j_\beta}_h
  -W^{ih}\Gamma^{j_\beta}_{h,l_\zeta}
  \ ,\quad
  \der{y^{l_\zeta}}g^{i_\alpha j}
  =
  -\Gamma^{i_\alpha}_{h,l_\zeta}W^{hj}
  -\Gamma^{i_\alpha}_hW\indices{^{hj}_{,l_\zeta}}
  \ ,\\
  &
  \der{y^{l_\zeta}}g^{i_\alpha j_\beta}
  =
  \Gamma^{i_\alpha}_{k,l_\zeta}W^{kh}\Gamma^{j_\beta}_h
  +\Gamma^{i_\alpha}_kW\indices{^{kh}_{,l_\zeta}}\Gamma^{j_\beta}_h
  +\Gamma^{i_\alpha}_kW^{kh}\Gamma^{j_\beta}_{h,l_\zeta}
  \ ,\\
  &\Gamma^{i_\alpha}_{h,l_\zeta}
  =
  \der{y^{l_\zeta}}
  \left(\Gamma\indices{^i_{hk}}u^k_\alpha\right) 
  =
  \delta^{\zeta\alpha}\Gamma\indices{^i_{hl}}
  \ ,\quad
  W\indices{^{ij}_{,l_\zeta}}
  =
  \delta^{il}u^j_\zeta+\delta^{jl}u^i_\zeta
  \ .
\end{align*}
Combining these expressions, we obtain 
\begin{align*}
  &
  \dot{x}^i
  =
  W^{ij}\xi_j
  -
  W^{ih}\Gamma^{j_\beta}_h\xi_{j_\beta}
  \quad,\quad
  \dot{X}^i_\alpha
  =
  -\Gamma^{i_\alpha}_hW^{hj}\xi_j
  +
  \Gamma^{i_\alpha}_kW^{kh}\Gamma^{j_\beta}_h\xi_{j_\beta}
  \\&
  \dot{\xi}_i
  =
  W^{hl}\Gamma^{k_\delta}_{l,i}
  \xi_h\xi_{k_\delta}
  -
  \frac{1}{2}
  \left(
  \Gamma^{h_\gamma}_{k,i}W^{kh}\Gamma^{k_\delta}_h+\Gamma^{h_\gamma}_kW^{kh}\Gamma^{k_\delta}_{h,i}
  \right)
  \xi_{h_\gamma}\xi_{k_\delta}
  \\&
  \dot{\xi}_{i_\alpha}
  =
  \Gamma^{h_\delta}_{k,i_\alpha}W^{kh}\Gamma^{k_\delta}_h
  \xi_{h_\gamma}\xi_{k_\delta}
  -\left(
  W\indices{^{hl}_{,i_\alpha}}\Gamma^{k_\delta}_l
  +W^{hl}\Gamma^{k_\delta}_{l,i_\alpha}
  \right)
  \xi_h\xi_{k_\delta}
  -\frac{1}{2}
  \left(
  W\indices{^{hk}_{,i_\alpha}}
  \xi_h\xi_k
  +
  \Gamma^{h_\delta}_kW\indices{^{kh}_{,i_\alpha}}\Gamma^{k_\delta}_h
  \xi_{h_\gamma}\xi_{k_\delta}
  \right)
  \ .
\end{align*}

\subsection{Acceleration and Polynomials for $\Cc$}
We can identify the covariant acceleration $\nabla_{\dot x_t}\dot{x}_t$ of
curves satisfying the MPP equations and hence normal
MPPs through their frame coordinates. Let $(u_t,\xi_t)$ satisfy \eqref{eq:hameq}.
Then $u_t$ is a horizontal lift of $x_t=\pi(u_t)$ and hence by
\eqref{eq:framecovder}, \eqref{eq:hder}, \eqref{eq:metricframecomp}, and
\eqref{eq:dotxi},
\begin{equation}
  \begin{split}
  u_t^{-1}\nabla_{\dot x_t}\dot{x}_t
  &=
  \frac{d}{dt}
  \begin{pmatrix}
    \xi(h_{u_t}(u_te_1)) \\
    \vdots \\
    \xi(h_{u_t}(u_te_d)) \\
  \end{pmatrix}
  =
  \begin{pmatrix}
    \dot \xi(h_{u_t}(u_te_1)) \\
    \vdots \\
    \dot \xi(h_{u_t}(u_te_d)) \\
  \end{pmatrix}
  +
  \begin{pmatrix}
    \xi(\partial_th_{u_t}(u_te_1)) \\
    \vdots \\
    \xi(\partial_th_{u_t}(u_te_d)) \\
  \end{pmatrix}
  \\&
  =
  -\begin{pmatrix}
    \xi(\partial_{h_{u_t}(u_te_1)}h_{u_t}(\pi_*(\dot{u}_t)) \\
    \vdots \\
    \xi(\partial_{h_{u_t}(u_te_d)}h_{u_t}(\pi_*(\dot{u}_t)) \\
  \end{pmatrix}
  +
  \begin{pmatrix}
    \xi(\partial_{h_{u_t}(\pi_*(\dot{u}_t))}h_{u_t}(u_te_1)) \\
    \vdots \\
    \xi(\partial_{h_{u_t}(\pi_*(\dot{u}_t))}h_{u_t}(u_te_d)) \\
  \end{pmatrix}
  \\&
  =
  \begin{pmatrix}
    \xi(\psi(u_t,R_{u_t}(u_te_1,\pi_*(\dot{u}_t)))) \\
    \vdots \\
    \xi(\psi(u_t,R_{u_t}(u_te_d,\pi_*(\dot{u}_t)))) \\
  \end{pmatrix}
  \ .
  \end{split}
  \label{eq:mppcovder}
\end{equation}

The fact that the covariant derivative vanishes for classical geodesic leads to
a definition of higher-order polynomials through the covariant derivative
by requiring $(\nabla_{\dot x_t})^k\dot{x}_t=0$ for a $k$th order polynomial, see 
e.g. \cite{leite_covariant_2008,hinkle_intrinsic_2014}. As
discussed above, compared to
classical geodesics, curves satisfying the MPP equations have extra $d^2$ degrees of freedom allowing the
curves to twist and deviate from being geodesic with respect to $\Cc$ while still 
satisfying the
horizontality constraint on $FM$. This makes it natural to ask if normal MPPs
relate to polynomials defined using $\Cc$.
For curves satisfying the
MPP equations, using \eqref{eq:mppcovder} and \eqref{eq:dotxi}, we have
\begin{equation*}
  \begin{split}
  u_t^{-1}(\nabla_{\dot x_t})^2\dot{x}_t
  &=
  \frac{d}{dt}
  \begin{pmatrix}
    \xi(\psi(u_t,R_{u_t}(u_te_1,\pi_*(\dot{u}_t)))) \\
    \vdots \\
    \xi(\psi(u_t,R_{u_t}(u_te_d,\pi_*(\dot{u}_t)))) \\
  \end{pmatrix}
  =
  \begin{pmatrix}
    \xi(\psi(u_t,\frac{d}{dt}R_{u_t}(u_te_1,\pi_*(\dot{u}_t)))) \\
    \vdots \\
    \xi(\psi(u_t,\frac{d}{dt}R_{u_t}(u_te_d,\pi_*(\dot{u}_t)))) \\
  \end{pmatrix}
  \ .
  \end{split}
\end{equation*}
Thus, in general, normal MPPs are not second order polynomials in the sense
$(\nabla_{\dot x_t})^2\dot x_t=0$ unless the curvature
$R_{u_t}(u_te_i,\pi_*(\dot{u}_t))$ is constant in $t$.

For comparison, in the Riemannian case, a variational formulation placing a
cost on covariant acceleration \cite{noakes_cubic_1989,camarinha_geometry_2001}
leads to cubic splines
\begin{equation*}
  (\nabla_{\dot x_t})^2\dot x_t
  =
  -R(\nabla_{\dot x_t}\dot x_t,x_t,)\dot x_t
  \ .
\end{equation*}
In \eqref{eq:mppcovder}, the curvature terms appear in
the covariant acceleration for normal MPPs while cubic splines leads to
the curvature term appearing in the third order derivative.

\section{Cometric Formulation and Low-Rank Generator}
\label{sec:rankdef}
We now investigate a cometric $g_{F^kM}+\lambda g_R$ where $g_R$ is Riemannian, 
$g_{F^kM}$ is a rank $k$ positive semi-definite inner product arising from $k$ linearly independent tangent
vectors, and $\lambda>0$ a weight. We assume that $g_{F^kM}$ is chosen so that $g_{F^kM}+\lambda g_R$ 
is invertible even though $g_{F^kM}$ is rank-deficient. 
The situation corresponds to extracting the first $k$ eigenvectors in
Euclidean space PCA. If the eigenvectors are estimated statistically from
observed data, this allows the estimation to be restricted to only the first $k$
eigenvectors.
In addition, an important practical implication
of the construction is that a numerical implementation need not transport a full $d\times
d$ matrix for the frame but a potentially much lower dimensional $d\times k$
matrix. This point is essential when dealing with high-dimensional data,
examples of which are landmark manifolds as discussed in
Section~\ref{sec:experiments}.

When using the frame bundle to model covariances, the sum formulation 
is natural to express as a cometric compared to a metric because, with the cometric
formulation, $g_{F^kM}+\lambda g_R$ represents a sum of covariance matrices
instead of a sum of precision matrices. Thus $g_{F^kM}+\lambda g_R$ can be intuitively thought
of as adding isotropic noise of variance $\lambda$ to the covariance represented
by $g_{F^kM}$.

To pursue this, let $F^kM$ denote the bundle of rank $k$ linear maps $\RR^k\rightarrow T_xM$.
We define a cometric by
\begin{equation*}
  \ip{\xi,\tilde{\xi}}
  = 
  \delta^{\alpha\beta}
  (\xi|h_u(u_\alpha))
  (\tilde{\xi}|h_u(u_\beta))
  +
  \lambda
  \ip{
    \xi,
    \tilde{\xi}
  }_{g_R}
\end{equation*}
for $\xi,\tilde{\xi}\in T_u^*F^kM$.
The sum over $\alpha,\beta$ is for $\alpha,\beta=1,\ldots,k$. The first term is
equivalent to the lift \eqref{eq:lift} of the cometric
$\ip{\xi,\tilde{\xi}}=\left(\xi|g_u(\hat{\xi})\right)$ given
$u:\RR^k\rightarrow T_xM$. Note that
in the definition \eqref{eq:g} of $g_u$, the map $u$ is not inverted,
thus the definition of the metric
immediately carries over to the rank-deficient case.
%

Let $(x^i,u^i_\alpha)$, $\alpha=1,\ldots,k$ be a coordinate system on $F^kM$.
The vertical distribution is in this case spanned by the $dk$ vector fields
$
  D_{j_\beta}=\partial_{u^j_\beta}
$.
Except for index sums being over $k$ instead of $d$ terms, the situation is 
thus similar to the full-rank case.
Note that $(\xi|\pi_*^{-1}w)=(\xi|w^jD_j)=w^i\xi_i$.
The cometric in coordinates is
\begin{equation*}
  \ip{\xi,\tilde{\xi}}
  = 
  \delta^{\alpha\beta}
  u_\alpha^i\xi_i
  u_\beta^j\tilde{\xi}_j
  +
  \lambda
  g_{R}^{ij}
  \xi_i
  \tilde{\xi}_j
  = 
  \xi_i
  \left(
  \delta^{\alpha\beta}
  u_\alpha^i
  u_\beta^j
  +
  \lambda
  g_{R}^{ij}
  \right)
  \tilde{\xi}_j
  =
  \xi_i
  W^{ij}
  \tilde{\xi}_j
\end{equation*}
with $W^{ij}=\delta^{\alpha\beta}u_\alpha^iu_\beta^j+\lambda g_{R}^{ij}$.
We can then write the corresponding sub-Riemannian metric $g_{F^kM}$ in terms of the adapted frame
$D$
\begin{equation}
  g_{F^kM}(\xi_hD^h+\xi_{h_\gamma}D^{h_\gamma})
  =
  W^{ih}\xi_hD_i
\end{equation}
because
$
  (\xi|g_{F^kM}(\tilde{\xi}))=\ip{\xi,\tilde{\xi}}
  =
  \xi_i
  W^{ij}
  \tilde{\xi}_j
  $.
That is, the situation is analogous to \eqref{eq:gadapted} except
the term $\lambda g_{R}^{ij}$ is added to $W^{ij}$.
\begin{figure}[t!]
  \begin{center}
      \includegraphics[width=.3\columnwidth,trim=0 140 0 0,clip=true]{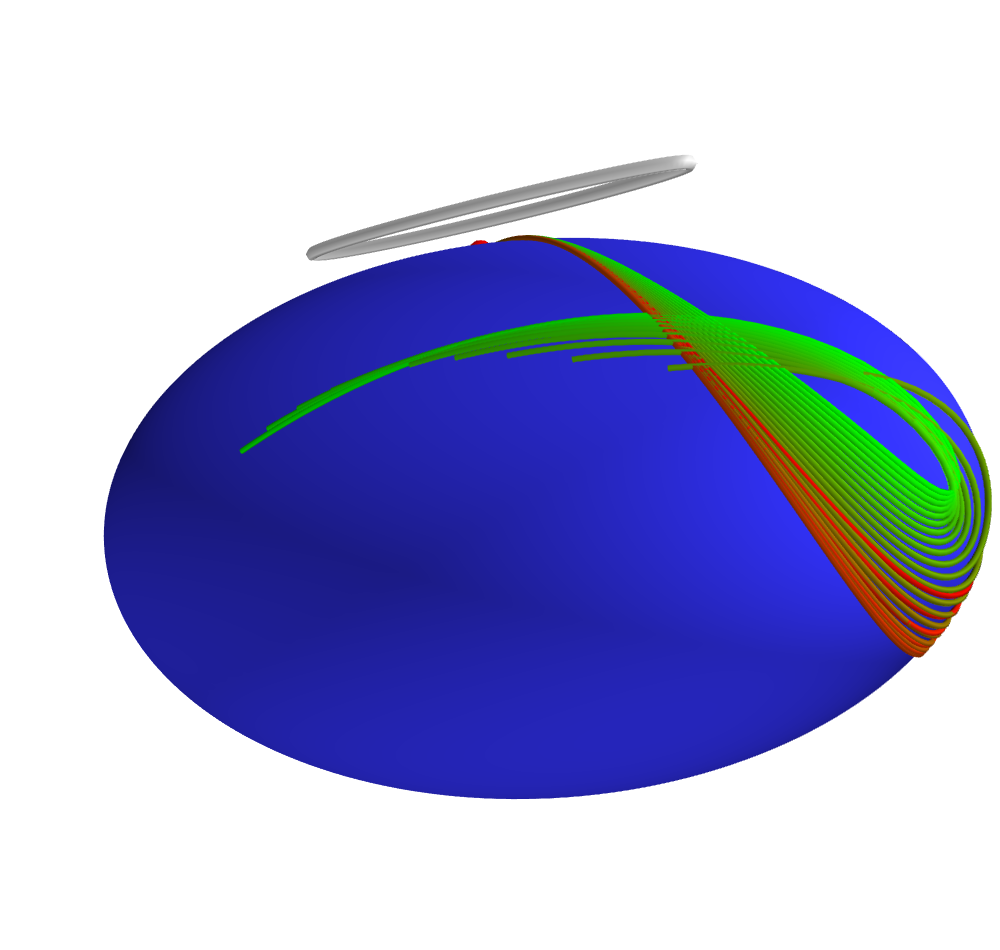}
      \includegraphics[width=.27\columnwidth,trim=0 140 200 0,clip=true]{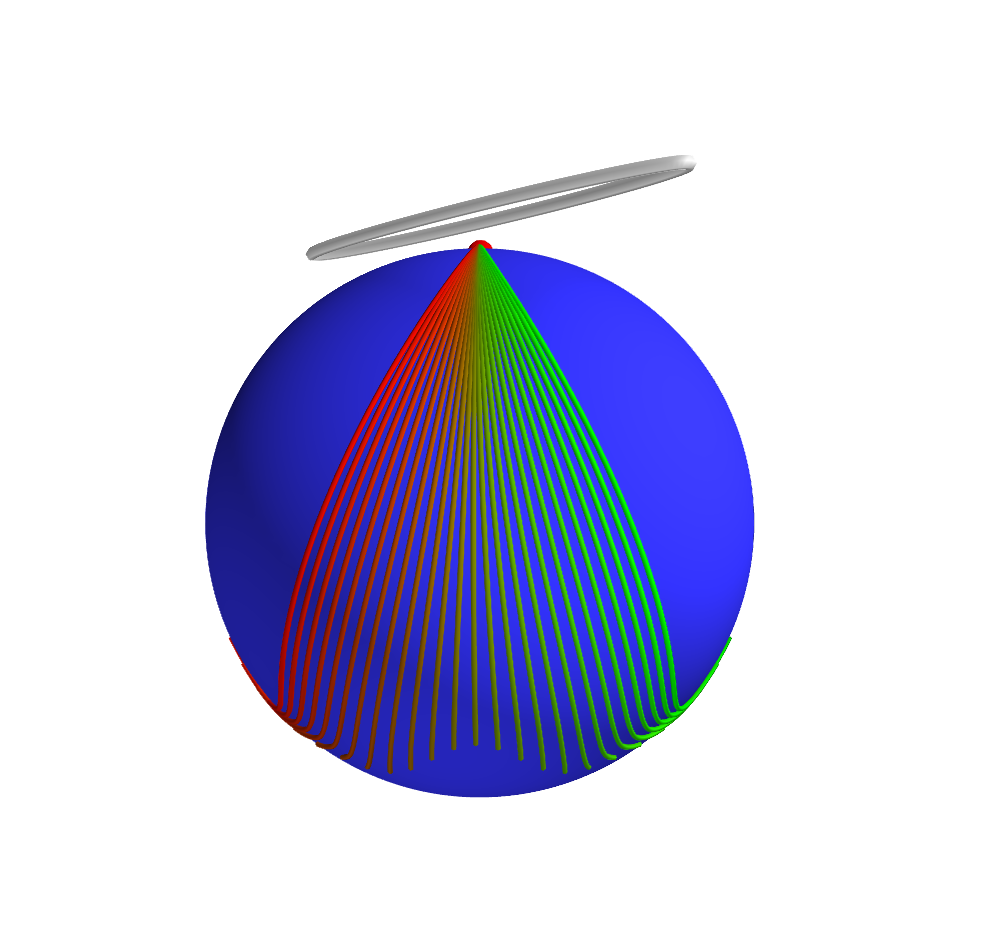}
      \includegraphics[width=.37\columnwidth,trim=0 200 0 200,clip=true]{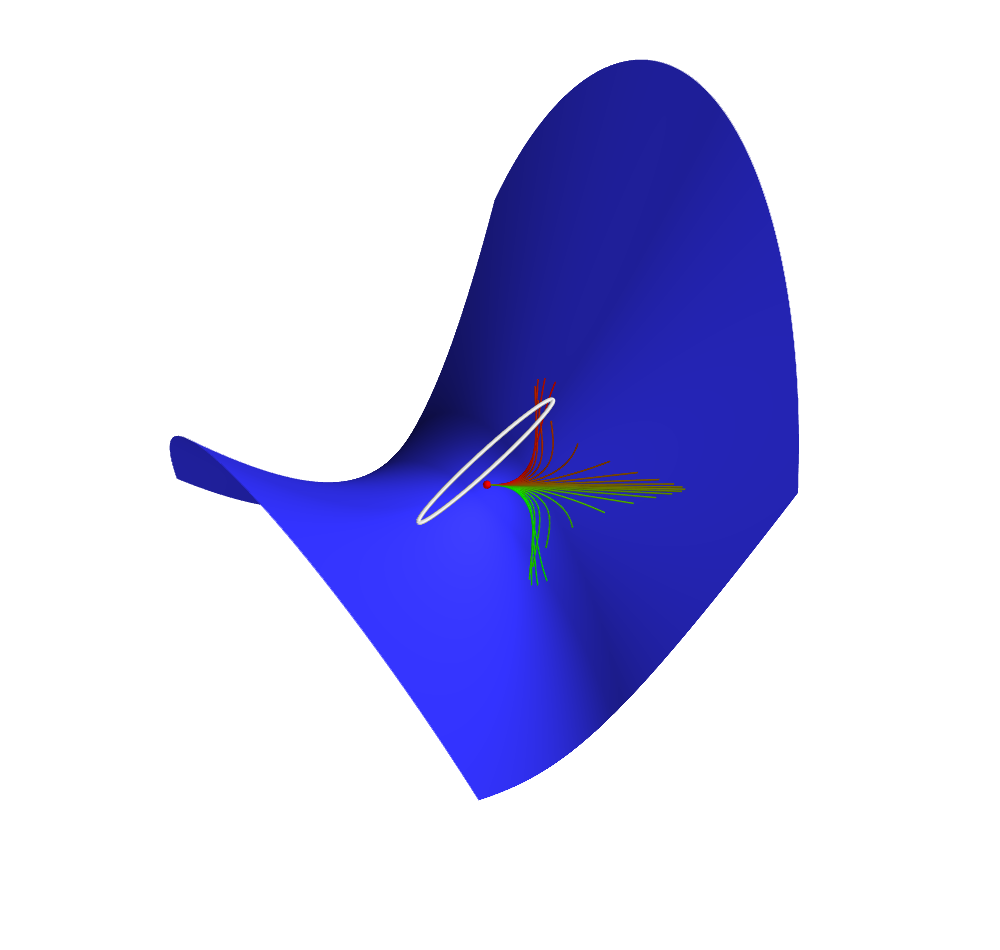}
      \\
      \includegraphics[width=.3\columnwidth,trim=0 0 0 0,clip=true]{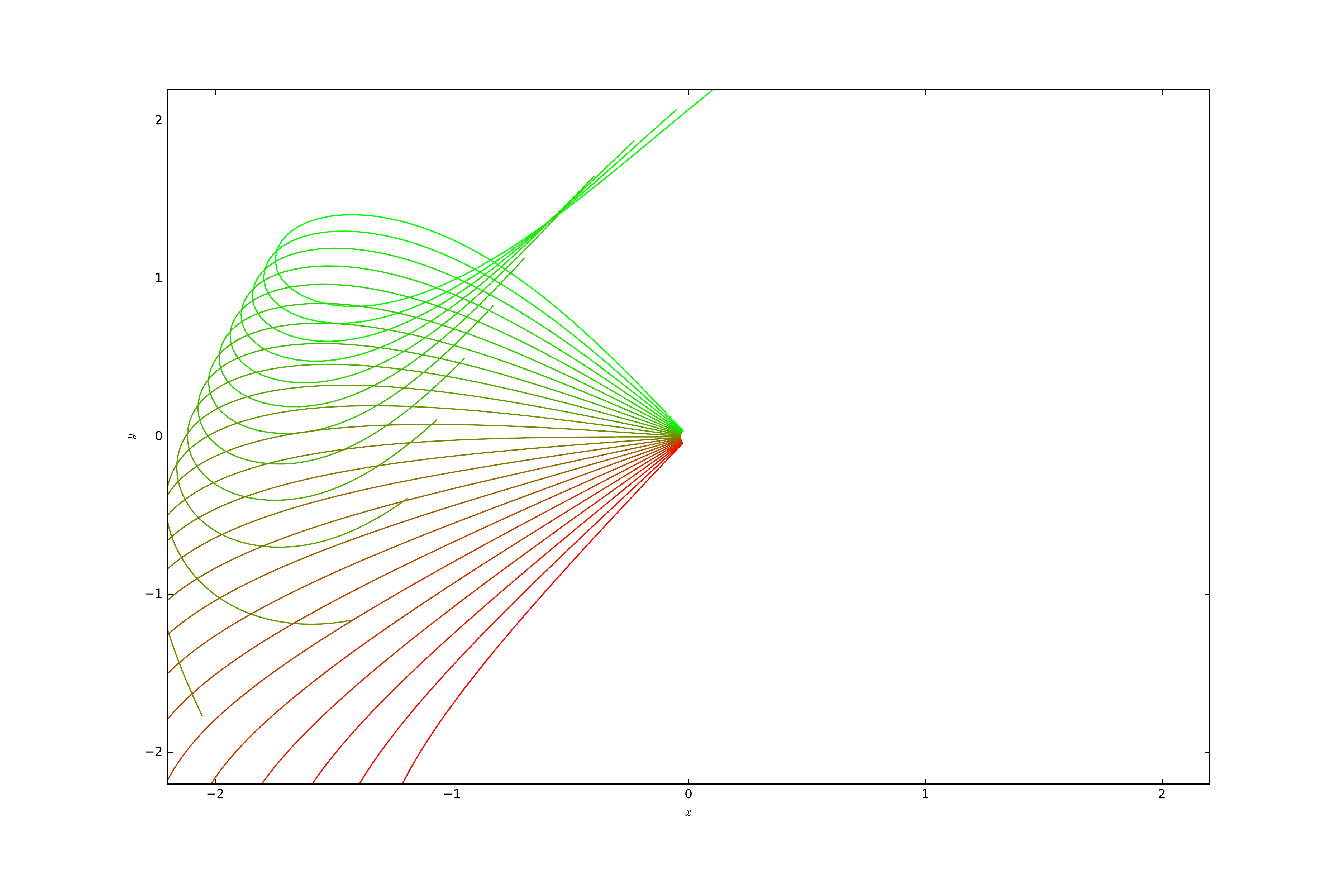}
      \includegraphics[width=.3\columnwidth,trim=0 0 0 0,clip=true]{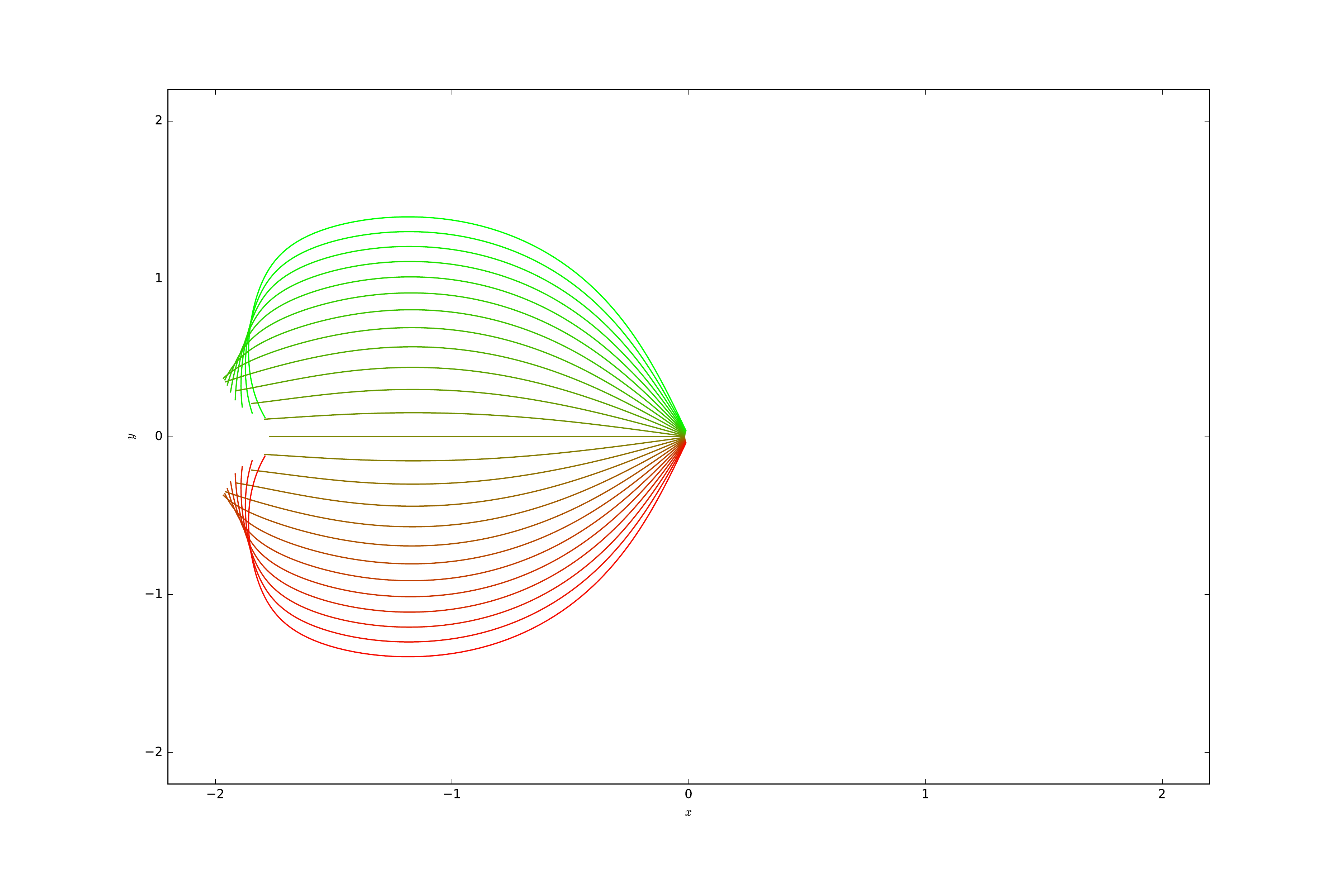}
      \hspace{.2cm}
      \includegraphics[width=.3\columnwidth,trim=0 0 0 0,clip=true]{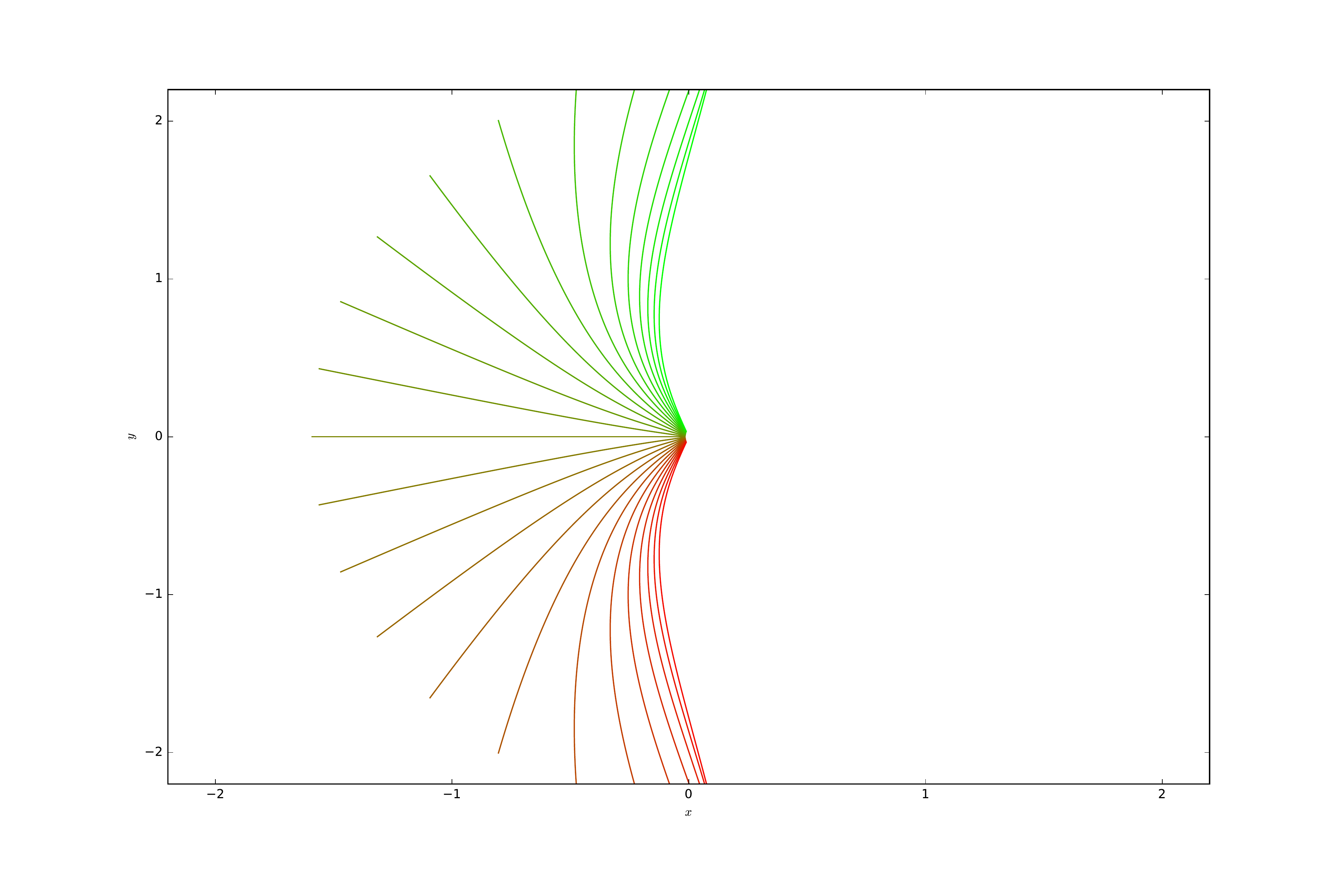}
      \hspace{.4cm}
  \end{center}
  \caption{
    Curves satisfying the MPP equations
    (top row) and corresponding anti-development (bottom row) on
    three surfaces embedded in $\RR^3$: An ellipsoid (left), a sphere (middle),
    a hyperbolic surface (right). The family of curves is generated by
    rotating by $\pi/2$ radians the anisotropic covariance represented in the initial 
    frame $u_0$ and displayed in the gray ellipse.
    }
  \label{fig:varmom}
\end{figure}

The geodesic system is again given by the Hamilton-Jacobi equations. As in the
full-rank case, the system is
specified by the derivatives of $g_{F^kM}$:
\begin{align*}
  &\der{y^{l}}g_{F^kM}^{ij}=W\indices{^{ij}_{,l}}
  \ ,\quad
  \der{y^{l}}g_{F^kM}^{i j_\beta}=
  -W\indices{^{ih}_{,l}}\Gamma^{j_\beta}_h
  -W^{ih}\Gamma^{j_\beta}_{h,l}
  \ ,\quad
  \der{y^{l}}g_{F^kM}^{i_\alpha j}=
  -\Gamma^{i_\alpha}_{h,l}W^{hj}
  -\Gamma^{i_\alpha}_hW\indices{^{hj}_{,l}}
  \ ,\\
  &\der{y^{l}}g_{F^kM}^{i_\alpha j_\beta}
  =
  \Gamma^{i_\alpha}_{k,l}W^{kh}\Gamma^{j_\beta}_h
  +\Gamma^{i_\alpha}_kW\indices{^{kh}_{,l}}\Gamma^{j_\beta}_h
  +\Gamma^{i_\alpha}_kW^{kh}\Gamma^{j_\beta}_{h,l}
  \ ,\\
  &\der{y^{l_\zeta}}g_{F^kM}^{ij}=W\indices{^{ij}_{,l_\zeta}}
  \ ,\quad
  \der{y^{l_\zeta}}g_{F^kM}^{i j_\beta}
  =
  -W\indices{^{ih}_{,l_\zeta}}\Gamma^{j_\beta}_h
  -W^{ih}\Gamma^{j_\beta}_{h,l_\zeta}
  \ ,\quad
  \der{y^{l_\zeta}}g_{F^kM}^{i_\alpha j}
  =
  -\Gamma^{i_\alpha}_hW\indices{^{hj}_{,l_\zeta}}
  -\Gamma^{i_\alpha}_{h,l_\zeta}W^{hj}
  \ ,\\
  &\der{y^{l_\zeta}}g_{F^kM}^{i_\alpha j_\beta}
  =
  \Gamma^{i_\alpha}_{k,l_\zeta}W^{kh}\Gamma^{j_\beta}_h
  +\Gamma^{i_\alpha}_kW\indices{^{kh}_{,l_\zeta}}\Gamma^{j_\beta}_h
  +\Gamma^{i_\alpha}_kW^{kh}\Gamma^{j_\beta}_{h,l_\zeta}
  \ ,\\
  &\Gamma^{i_\alpha}_{h,l_\zeta}
  =
  \der{y^{l_\zeta}}
  \left(\Gamma\indices{^i_{hk}}u^k_\alpha\right) 
  =
  \delta^{\zeta\alpha}\Gamma\indices{^i_{hl}}
  \ ,\quad
  W\indices{^{ij}_{,l}}
  =
  \lambda
  g\indices{_R^{ij}_{,l}}
  \ ,\quad
  W\indices{^{ij}_{,l_\zeta}}
  =
  \delta^{il}u^j_\zeta+\delta^{jl}u^i_\zeta
  \ .
\end{align*}
Note that the introduction of the Riemannian metric $g_R$ implies that $W^{ij}$
are now dependent on the manifold coordinates $x^i$.

\section{Numerical Experiments}
\label{sec:experiments}
We aim at visualizing most probable paths for the driving process and projections 
of curves satisfying the MPP equations \eqref{eq:hameq}
in two cases: On 2D surfaces embedded in $\RR^3$ and on finite dimensional
landmark manifolds that arise from equipping a subset of the diffeomorphism
group with a right-invariant metric and letting the action descend to the
landmarks by a left action. The surface examples are implemented in Python using
the Theano \cite{the_theano_development_team_theano:_2016} framework for symbolic 
operations, automatic differentiation, and numerical evaluation.
The landmark equations are detailed below and implemented in Numpy using Numpy's
standard ODE integrators. The code for running the experiments is available at 
{\texttt http://bitbucket.com/stefansommer/mpps/}.
\begin{figure}[t]
  \begin{center}
      \includegraphics[width=.3\columnwidth,trim=220 250 250 200,clip=true]{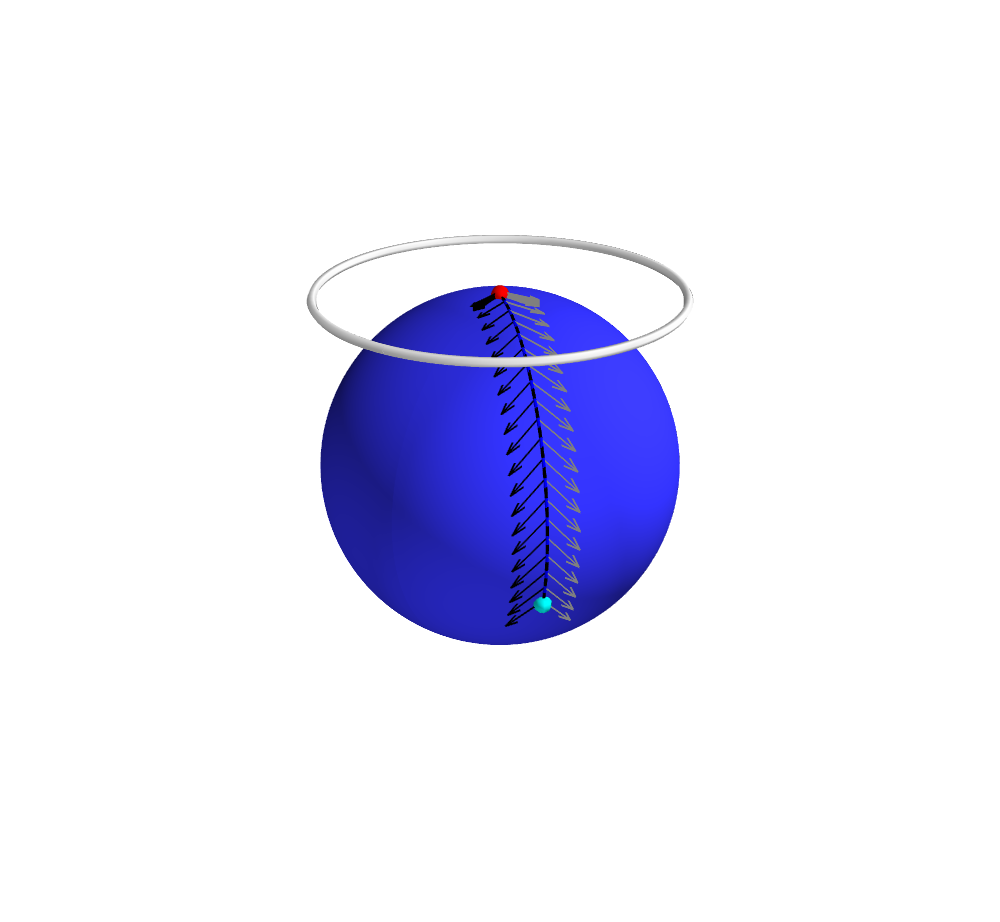}
      \includegraphics[width=.3\columnwidth,trim=220 250 250 200,clip=true]{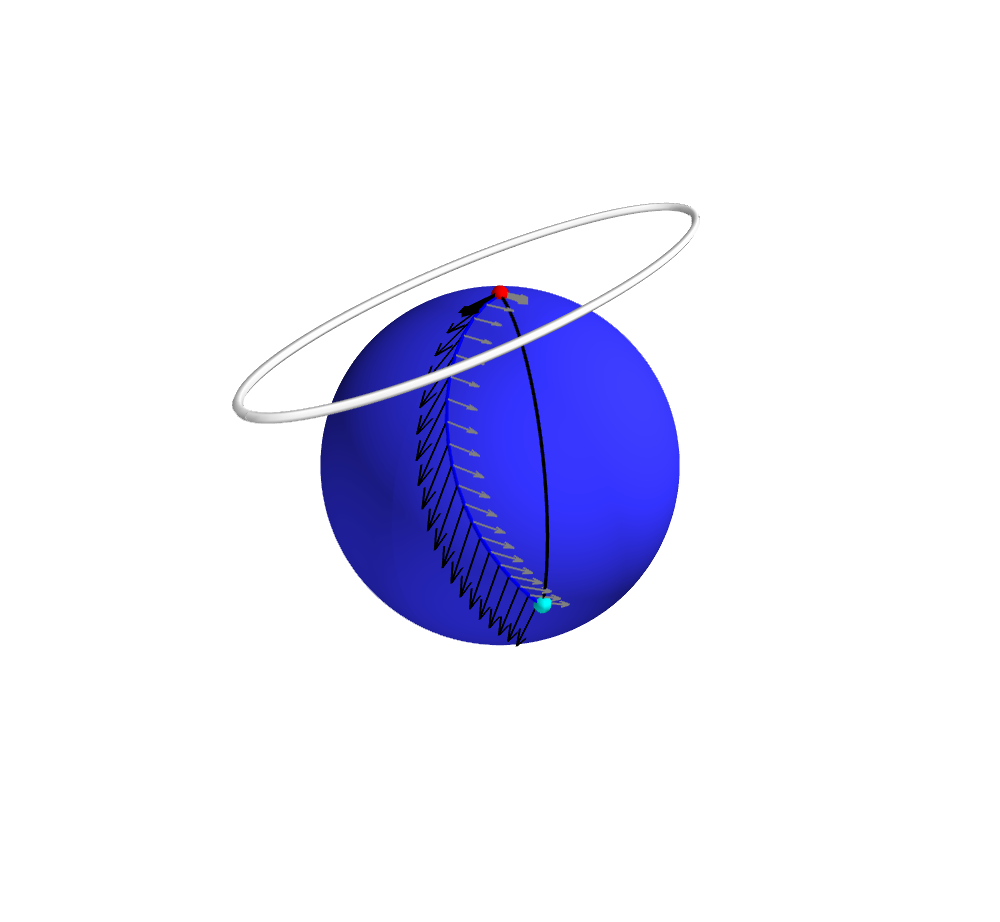}
      \includegraphics[width=.3\columnwidth,trim=220 250 250 200,clip=true]{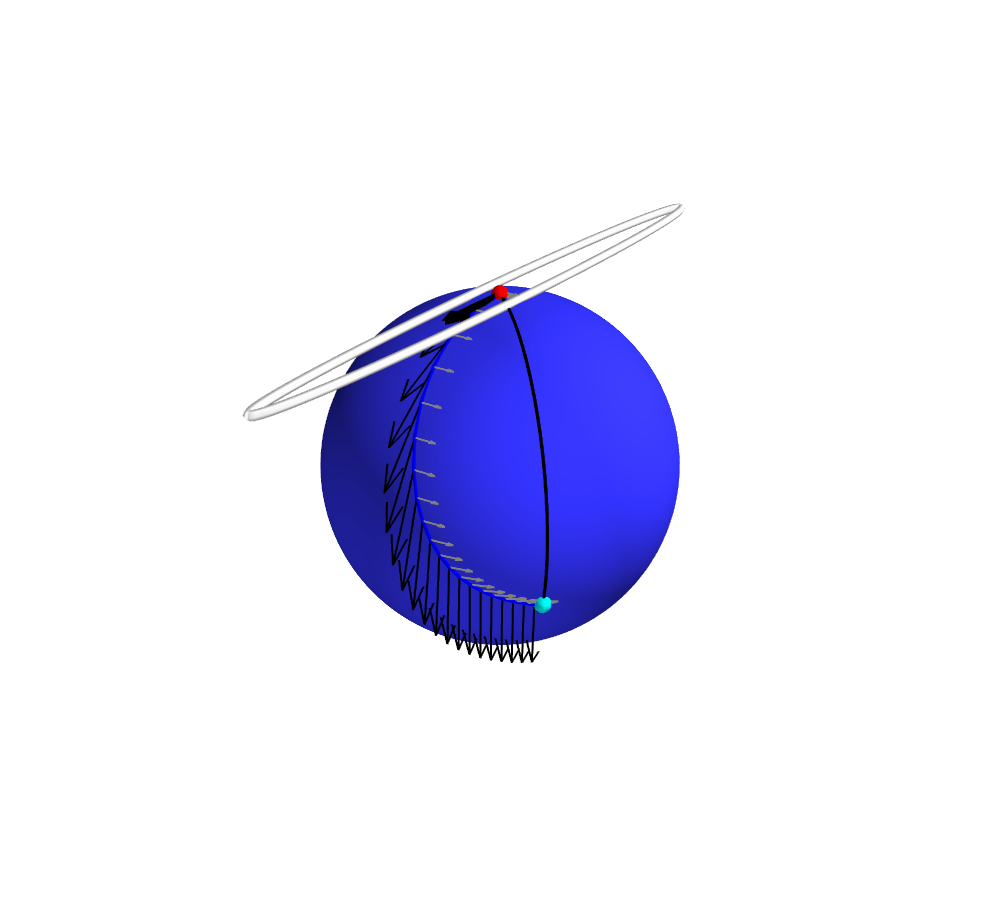}
      \\
      \vspace{-.2cm}
      \includegraphics[width=.3\columnwidth,trim=0 50 0 0,clip=true]{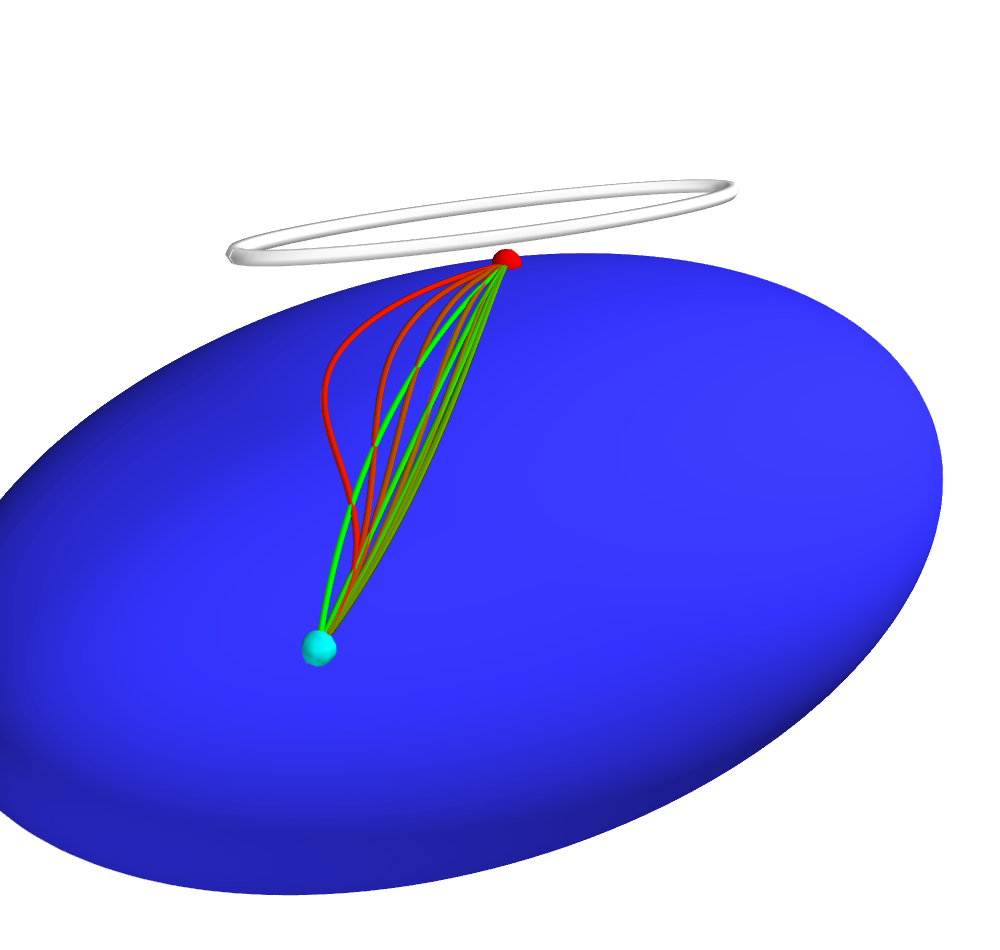}
      \includegraphics[width=.3\columnwidth,trim=0 50 0 0,clip=true]{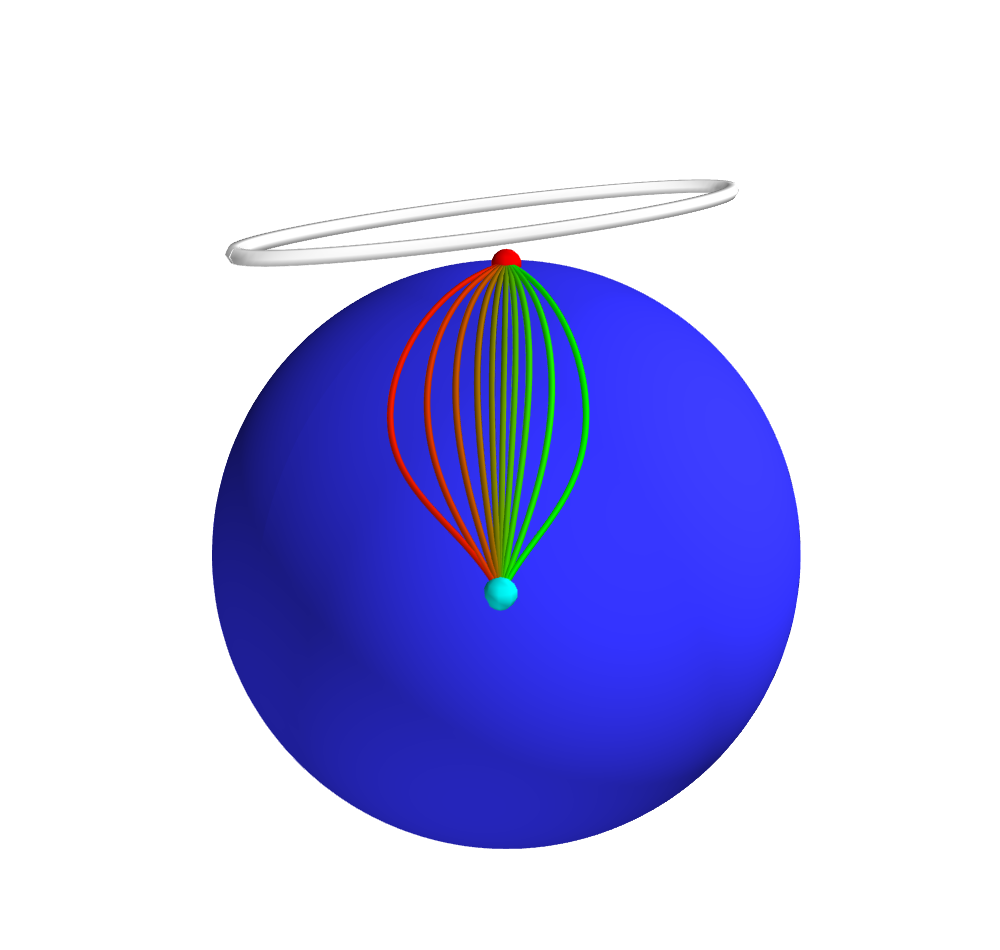}
      \includegraphics[width=.32\columnwidth,trim=200 240 200 300,clip=true]{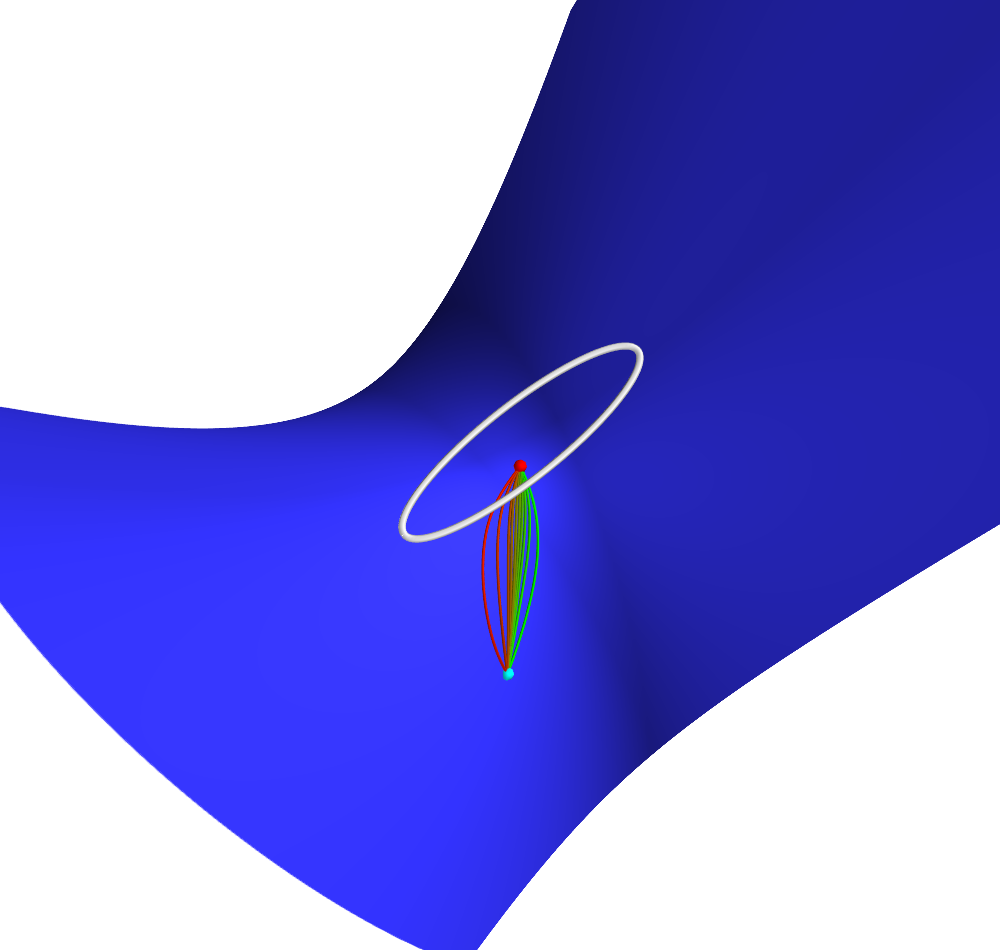}
      \\
      \hspace{-.2cm}
      \includegraphics[width=.3\columnwidth,trim=0 0 0 0,clip=true]{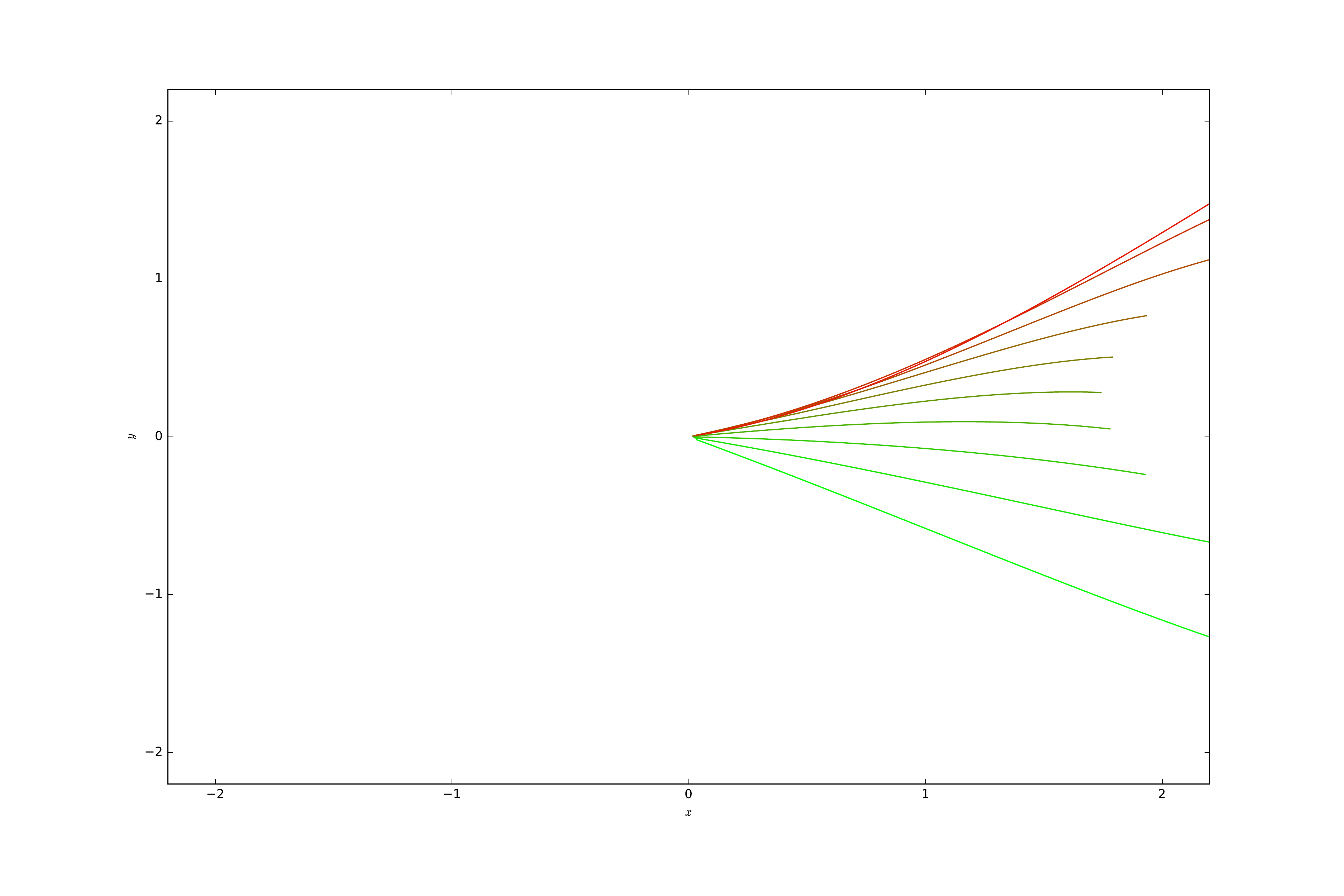}
      \hspace{.2cm}
      \includegraphics[width=.3\columnwidth,trim=0 0 0 0,clip=true]{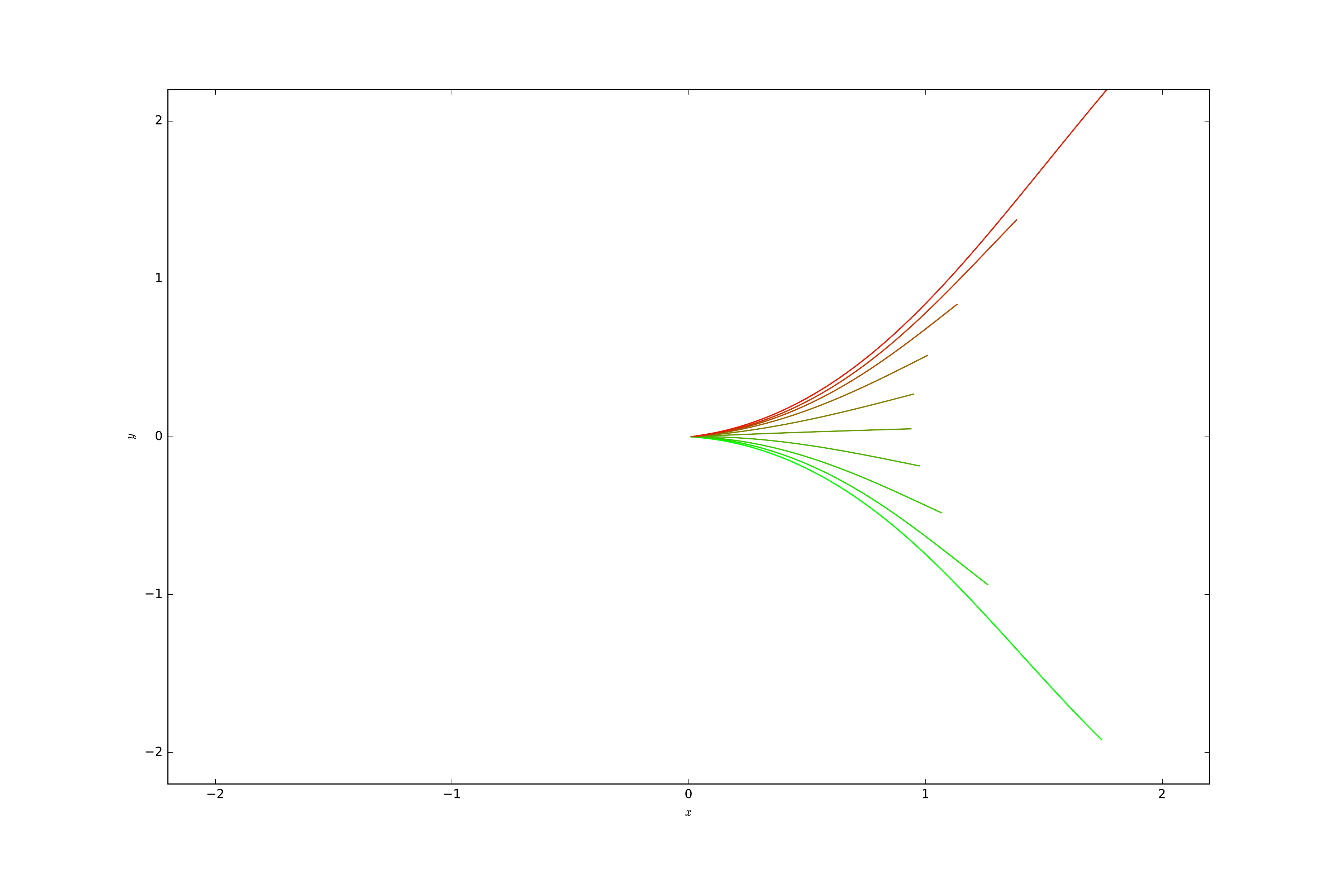}
      \includegraphics[width=.3\columnwidth,trim=0 0 0 0,clip=true]{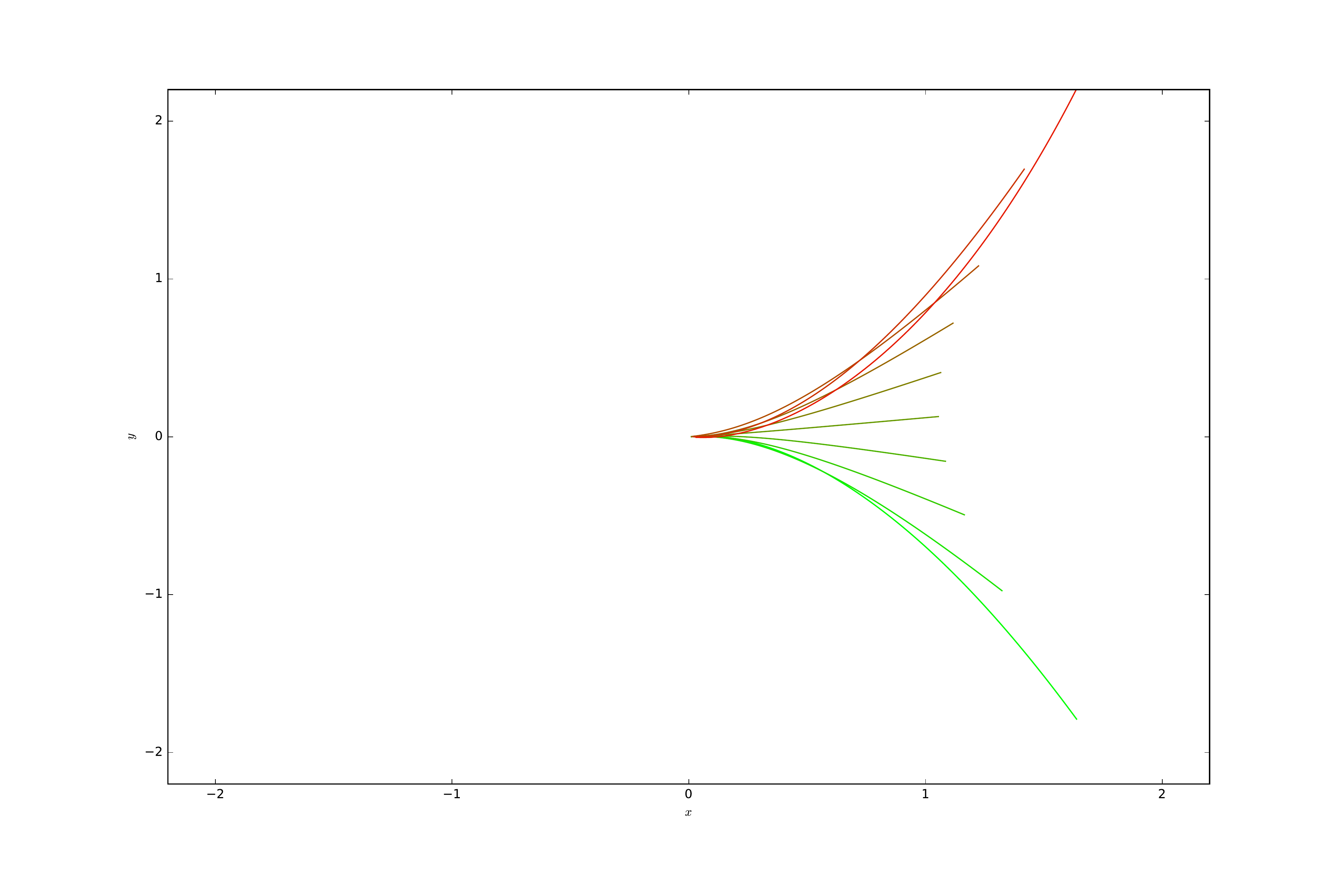}
      \hspace{.6cm}
  \end{center}
  \caption{
    Minimizing normal MPPs between two fixed points (red/cyan). 
    From isotropic covariance (top row, left) to anisotropic (top row, right) on
    $\mathbb{S}^2$. Compare with minimizing Riemannian geodesic (black curve). 
    The MPP travels longer in the directions of high variance.
    Families of curves 
    (middle row) and corresponding 
    anti-development (bottom row) on
    the three surfaces in Figure~\ref{fig:varmom}. The family of curves is generated by rotating the
    covariance matrix as in Figure~\ref{fig:varmom}. Notice how the varying
    anisotropy affects the resulting minimizing curves, and how the
    anti-developed curves end at different points in $\RR^2$.
    }
  \label{fig:varlog}
\end{figure}

\subsection{Embedded Surfaces}
We visualize normal MPPs and projections 
of curves satisfying the MPP equations \eqref{eq:hameq}
on surfaces embedded in $\RR^3$ in three cases:
The sphere $\mathbb{S}^2$, on an ellipsoid, and on a hyperbolic surface. The 
surfaces are chosen
in order to have both positive and negative curvature, and to have varying
degree of symmetry. In all cases, an open subset of the surfaces are 
represented in a single
chart by a map $F:\RR^2\rightarrow\RR^3$. For the sphere and ellipsoid, this
gives a representation of the surface except for the south
pole. The metric and Christoffel symbols are calculated using the symbolic
differentiation features of Theano. The integration are performed by a simple 
Euler integrator.
\begin{figure}[t!]
  \begin{center}
      \includegraphics[width=.3\columnwidth,trim=0 0 0 0,clip=true]{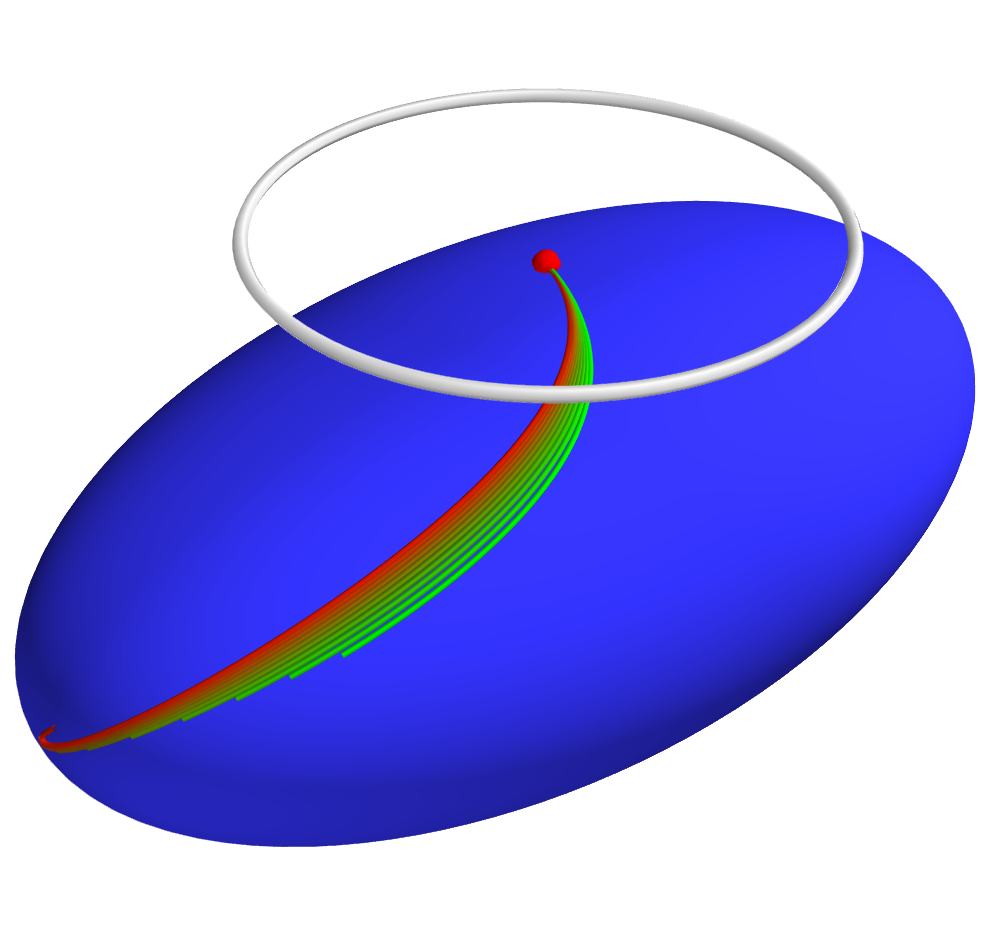}
      \includegraphics[width=.3\columnwidth,trim=0 0 0 0,clip=true]{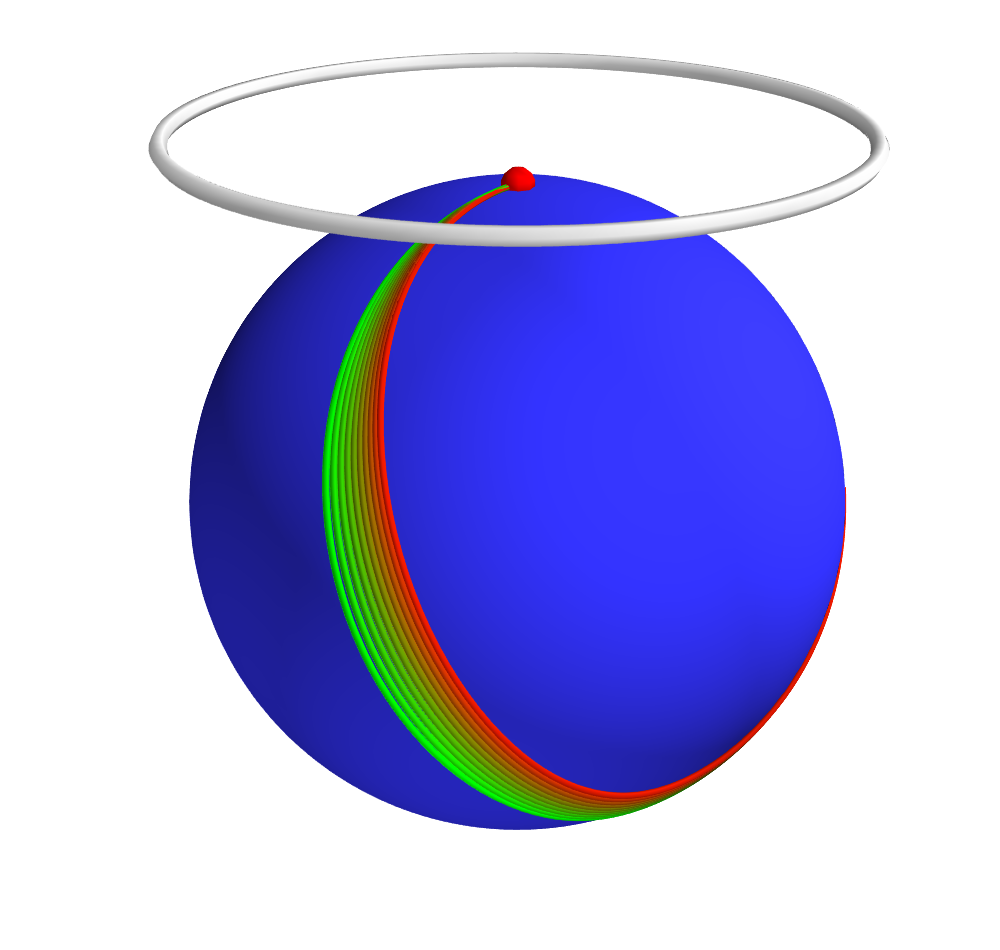}
      \includegraphics[width=.3\columnwidth,trim=100 100 100 100,clip=true]{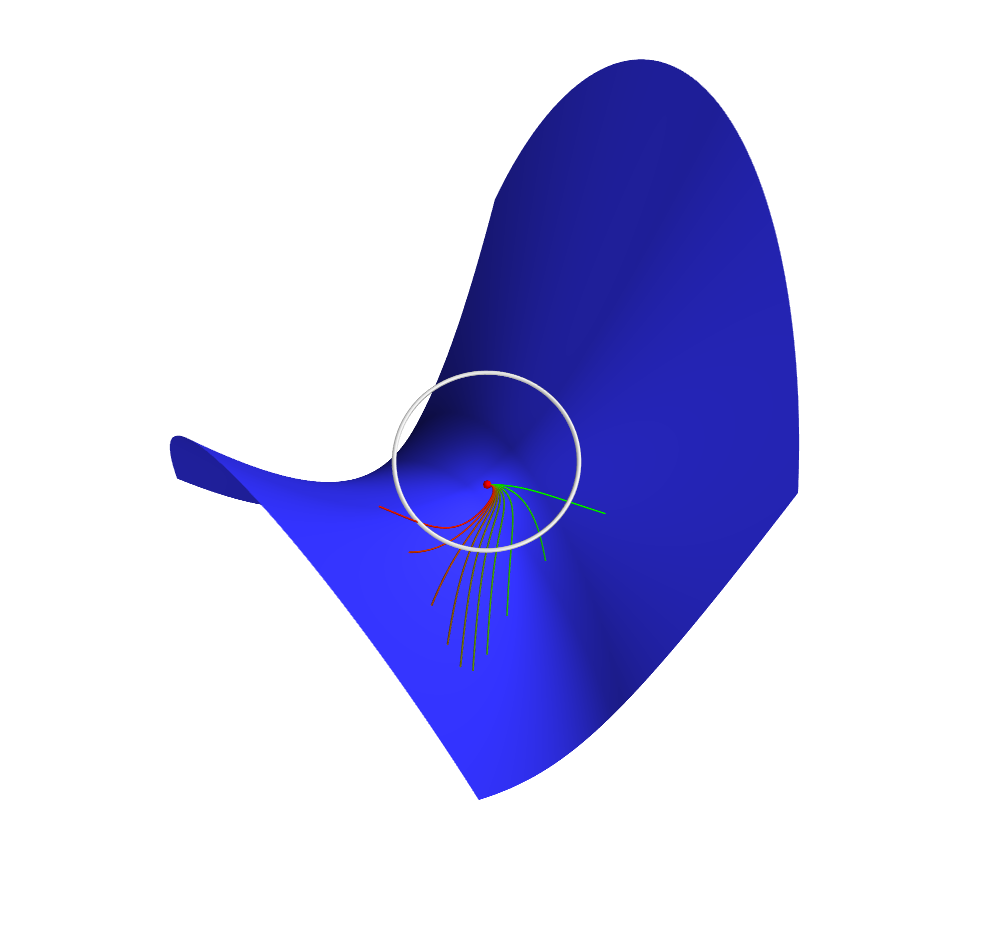}
      \\
  \vspace{-.1cm}
      \includegraphics[width=.3\columnwidth,trim=0 0 0 0,clip=true]{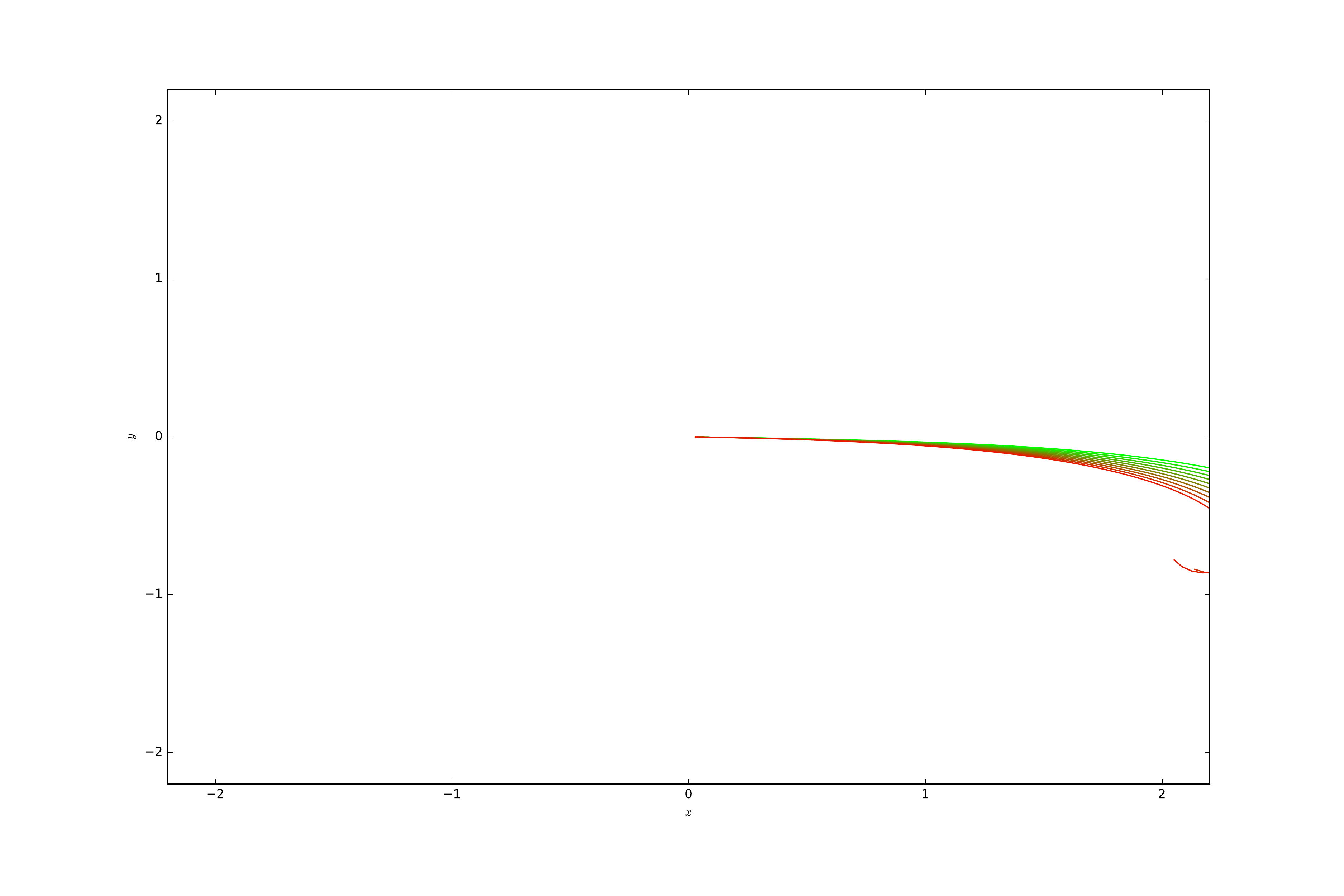}
      \includegraphics[width=.3\columnwidth,trim=0 0 0 0,clip=true]{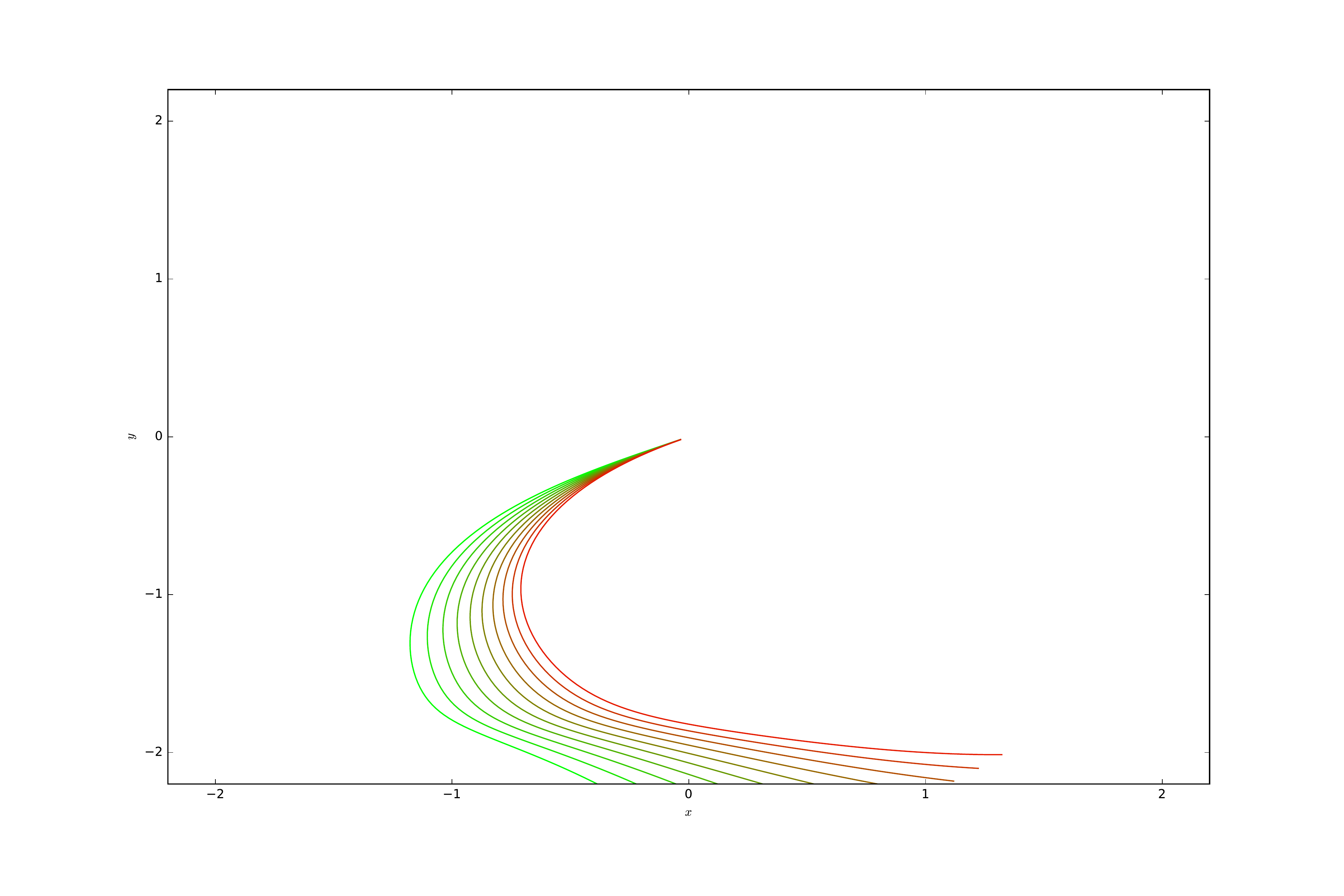}
      \includegraphics[width=.3\columnwidth,trim=0 0 0 0,clip=true]{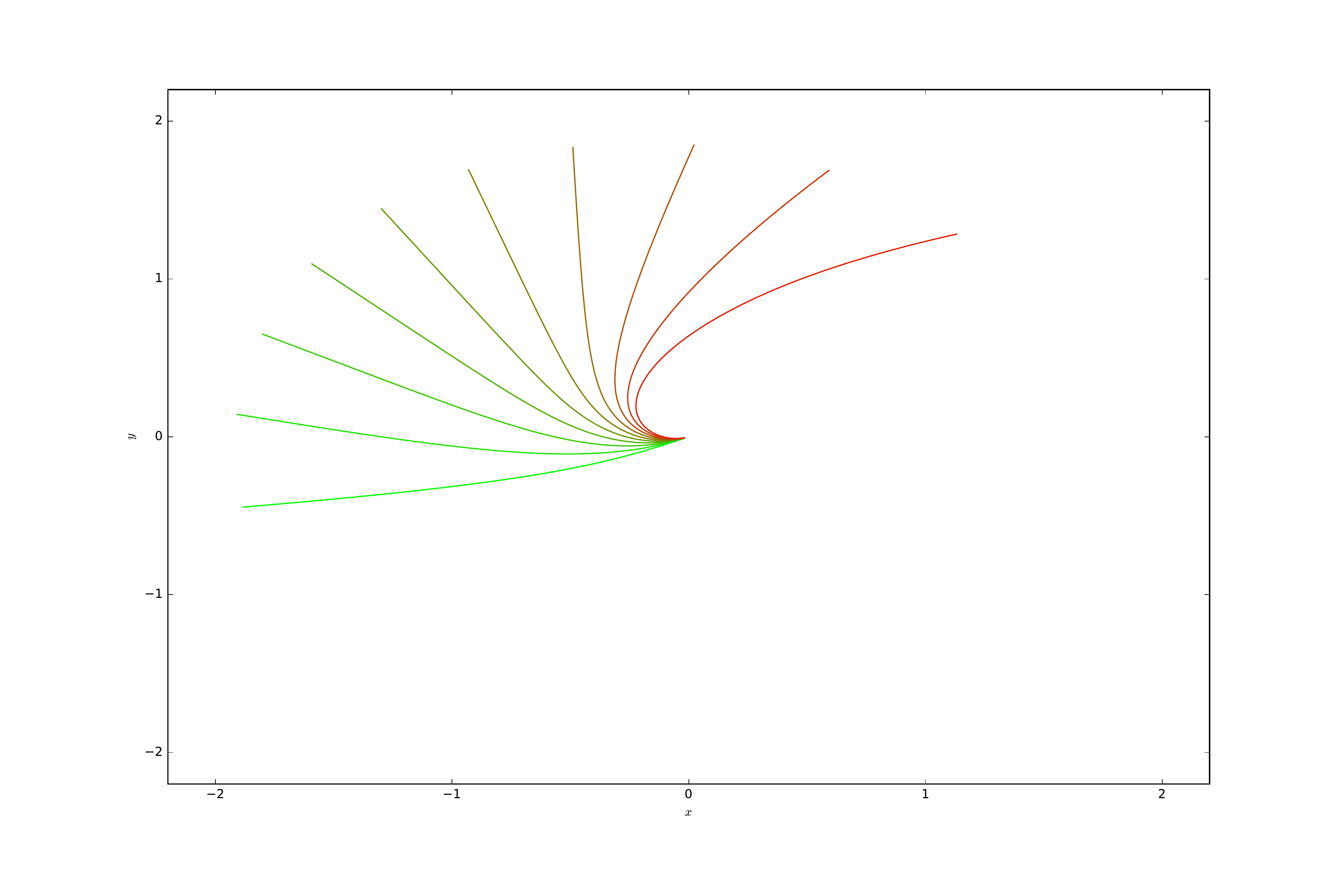}
      \\
  \vspace{-.1cm}
      \includegraphics[width=.3\columnwidth,trim=0 0 0 0,clip=true]{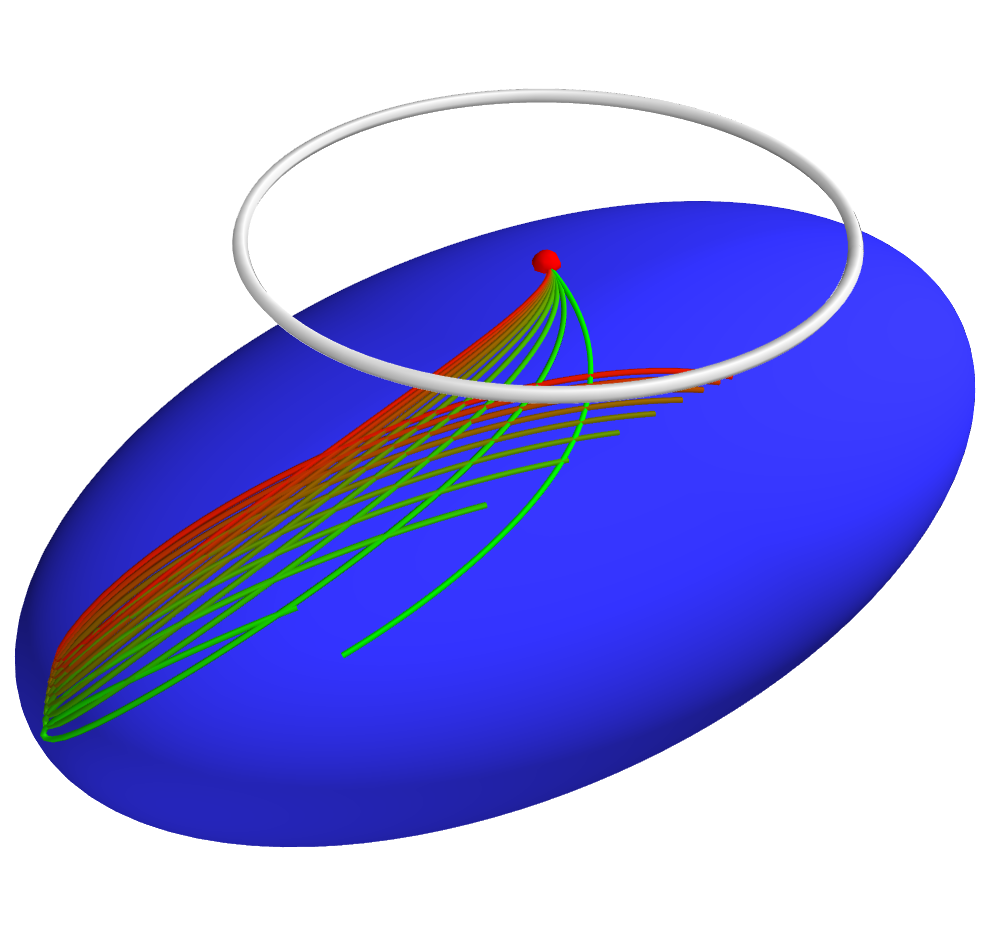}
      \includegraphics[width=.3\columnwidth,trim=0 0 0 0,clip=true]{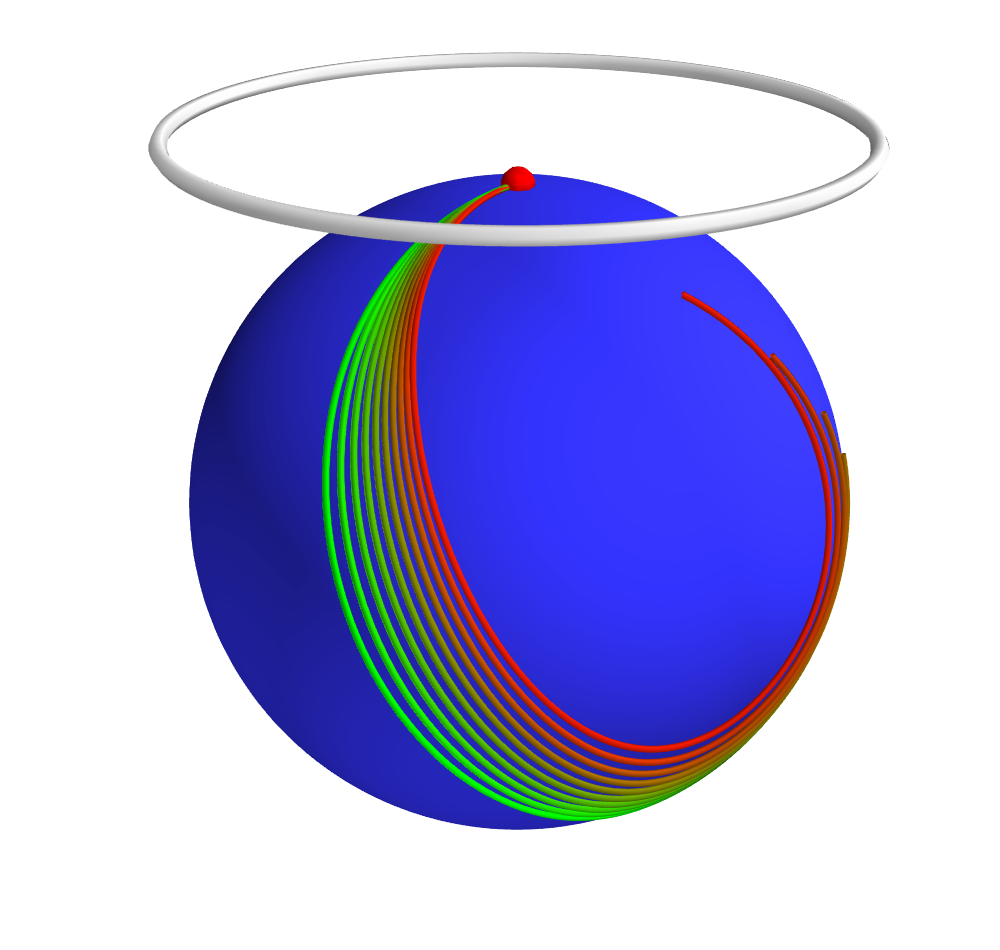}
      \includegraphics[width=.3\columnwidth,trim=100 100 100 100,clip=true]{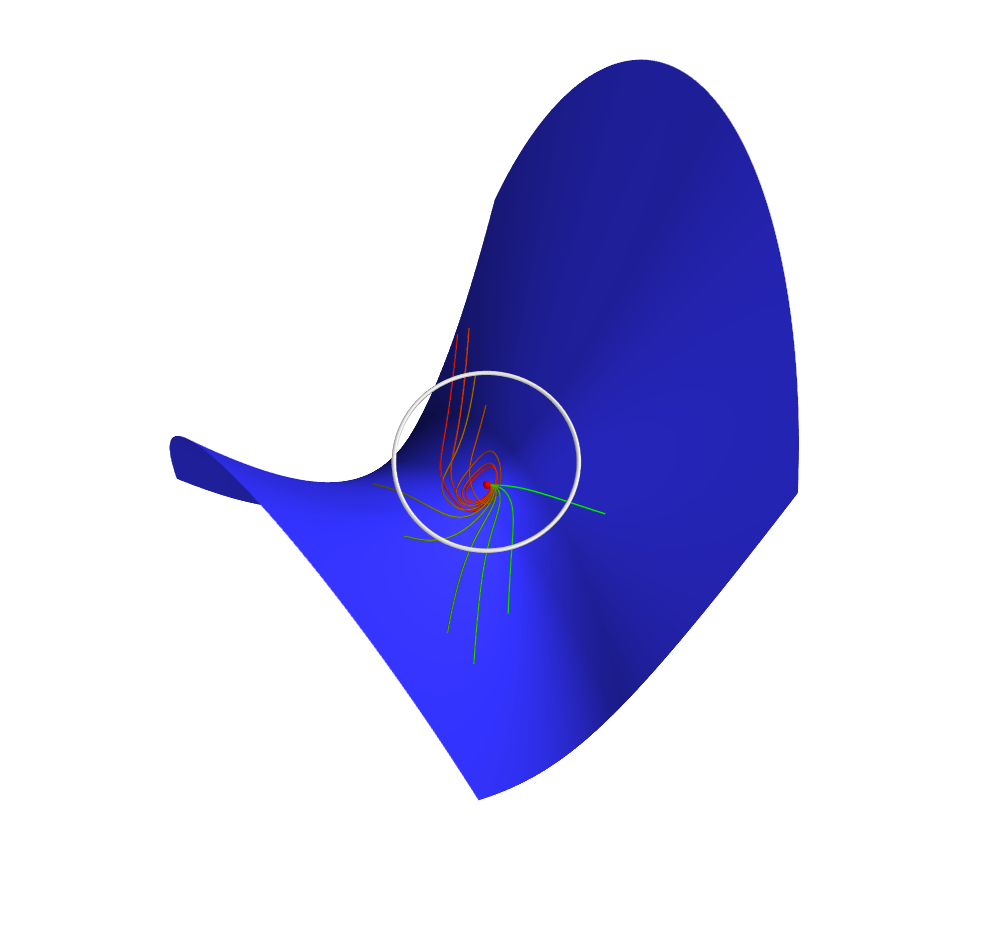}
      \\
  \vspace{-.1cm}
      \includegraphics[width=.3\columnwidth,trim=0 0 0 0,clip=true]{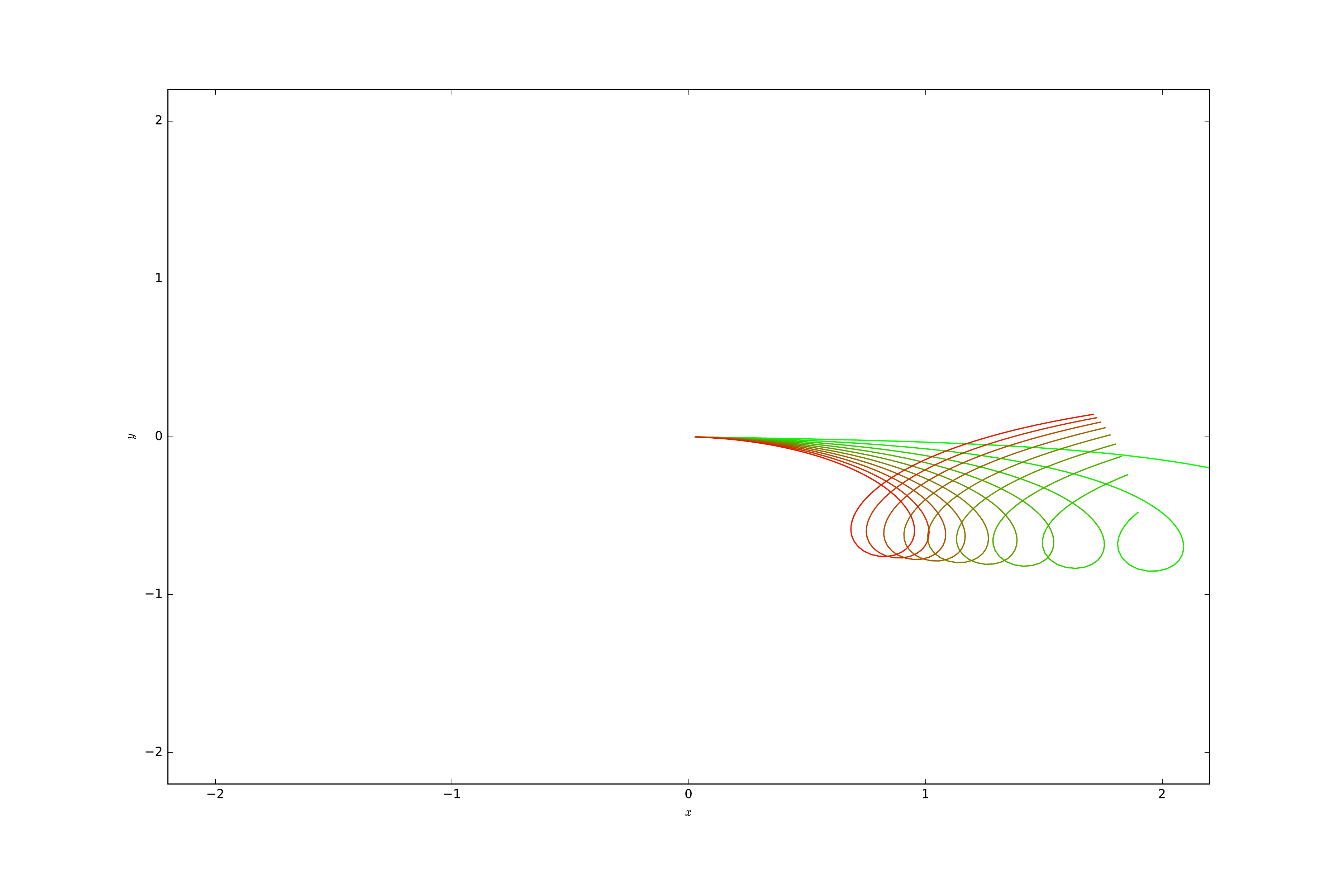}
      \includegraphics[width=.3\columnwidth,trim=0 0 0 0,clip=true]{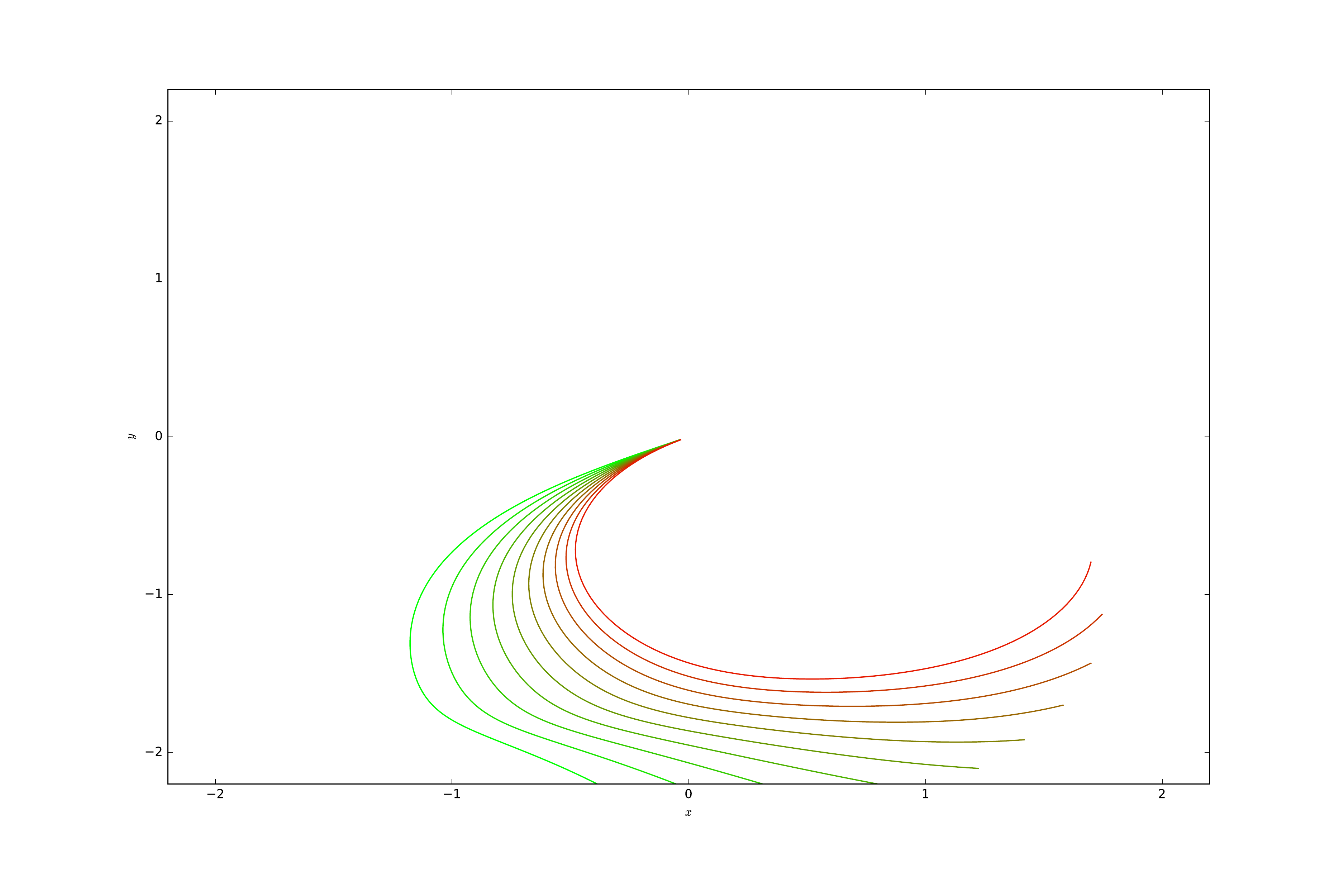}
      \includegraphics[width=.3\columnwidth,trim=0 0 0 0,clip=true]{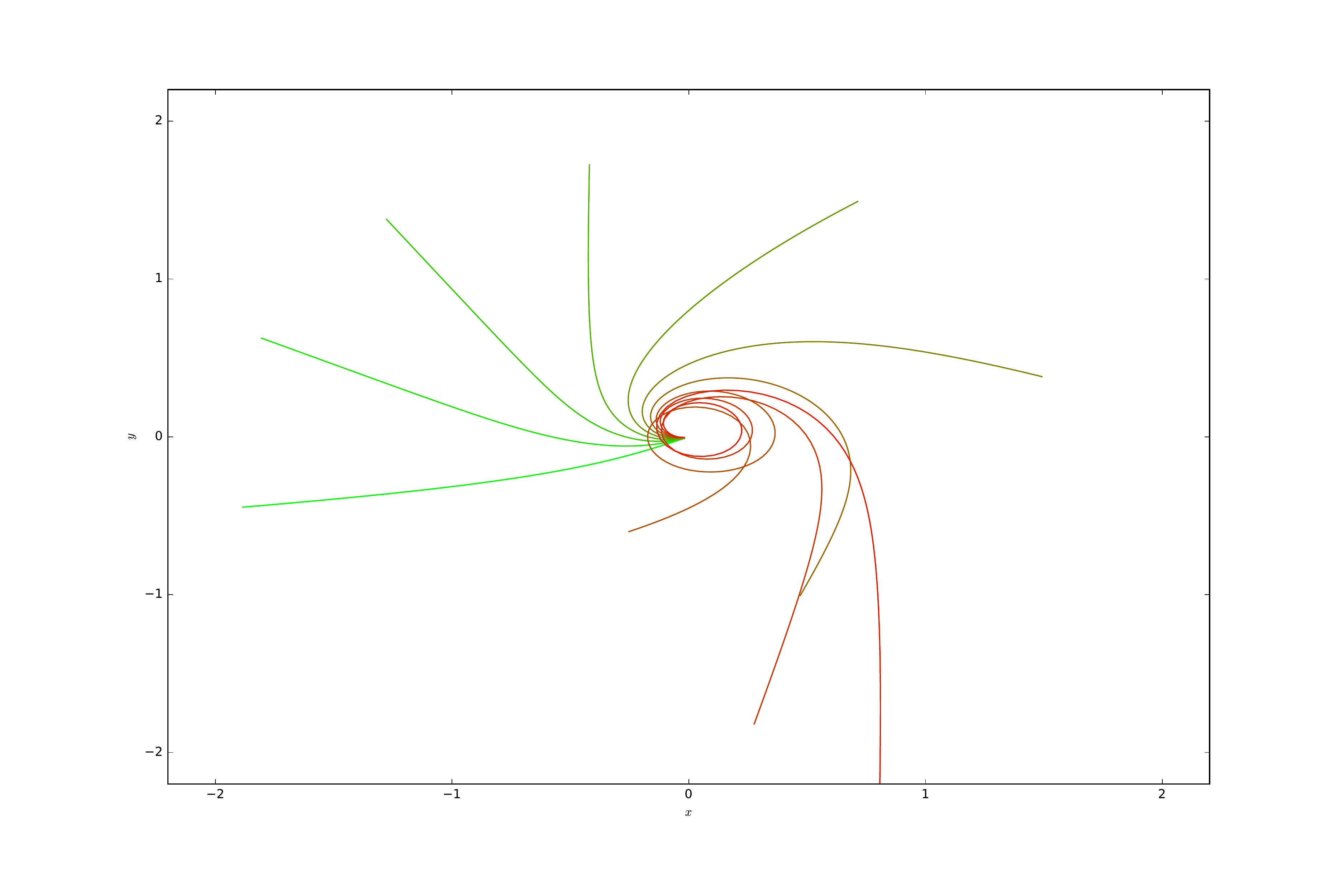}
  \end{center}
  \caption{
    With the setup of Figure~\ref{fig:varmom},\ref{fig:varlog}, generated
    families of curves by varying the vertical
    $V^*FM$ part of the initial momentum $\xi_0\in T^*FM$ but keeping the base
    point and frame $u_0$ fixed. The vertical part allows varying degree of
    ``twisting'' of the curve.
    }
  \label{fig:varxia}
\end{figure}

Figures~\ref{fig:varmom}-\ref{fig:varxia} shows families of curves satisfying the
MPP equations in three cases: 
(1) With fixed 
starting point $x_0\in M$ and initial velocity $\dot x_0\in TM$ but varying
anisotropy represented by changing frame $u$ in the fiber above $x_0$;
(2) minimizing normal MPPs with fixed starting point and endpoint $x_0,x_1\in M$ but changing 
frame $u$ above $x_0$; (3) fixed starting point $x_0\in M$ and frame $u$ but
varying $V^*FM$ vertical part of the initial momentum $\xi_0\in T^*FM$. The first 
and second cases 
thus shows the effect of varying anisotropy while the third case illustrates the
effect of the ``twist'' the $d^2$ degrees in the vertical momentum allows.
Note the displayed anti-developed curves in $\RR^2$ that for classical $\Cc$
geodesics would always be straight lines.

\subsection{LDDMM Landmark Equations}
We here give a example of the MPP equations using the
finite dimensional landmark manifolds that arise from right invariant metrics on
subsets of the diffeomorphism group in the Large Deformation Diffeomorphic
Metric Mapping (LDDMM) framework \cite{younes_shapes_2010}. 
The LDDMM metric can be conveniently expressed as a cometric, and, using a
rank-deficient inner product $g_{F^kM}$ as discussed in
Section~\ref{sec:rankdef}, we can obtain a reduction of the system
of equations to $2(2N+2Nk)$ compared to $2(2N+(2N)^2)$ with $N$ landmarks 
in $\RR^2$.
\begin{figure}[t!]
  \begin{center}
      \includegraphics[width=.19\columnwidth,trim=80 60 80 60,clip=true]{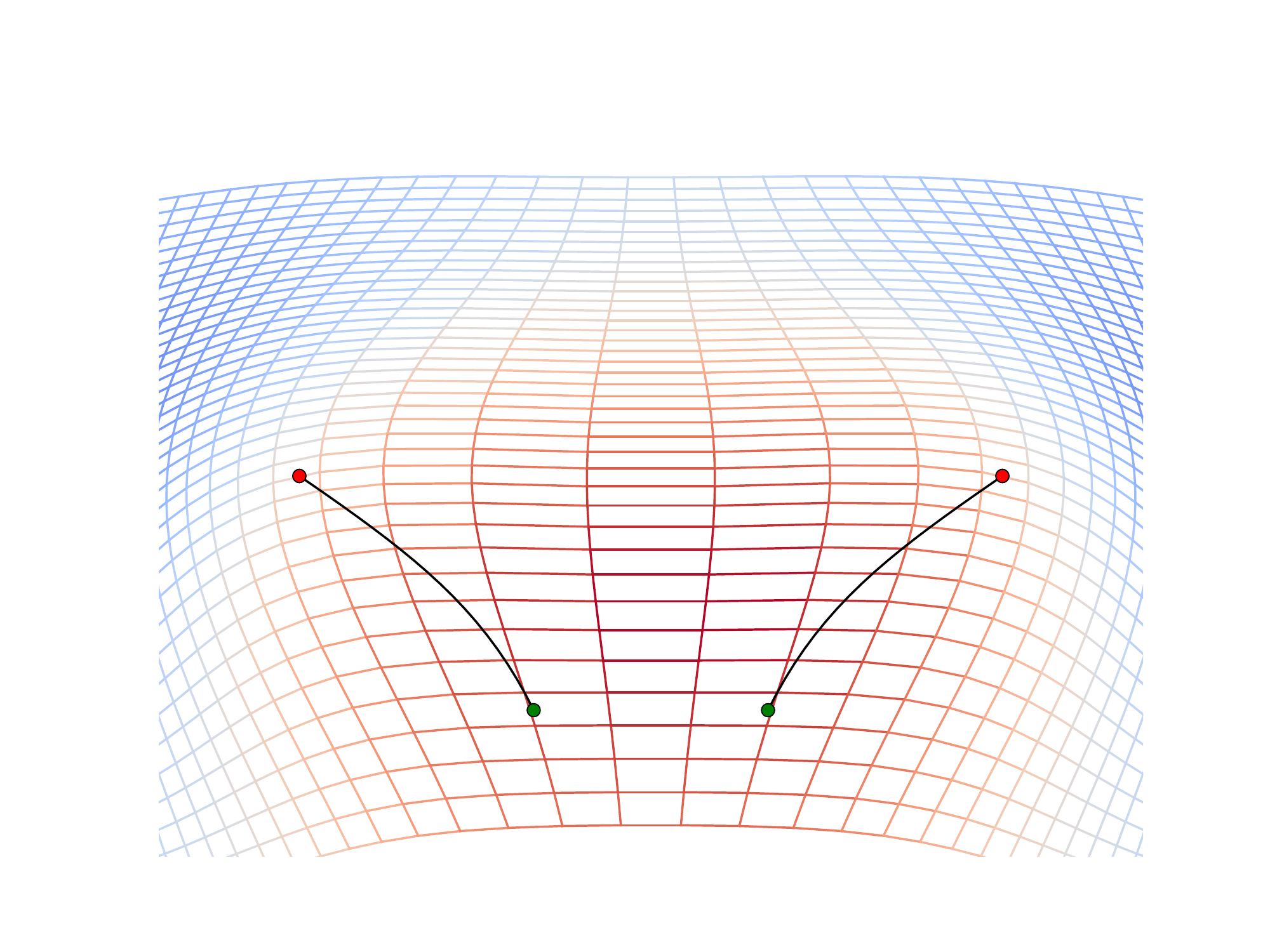}
      \includegraphics[width=.19\columnwidth,trim=80 60 80 60,clip=true]{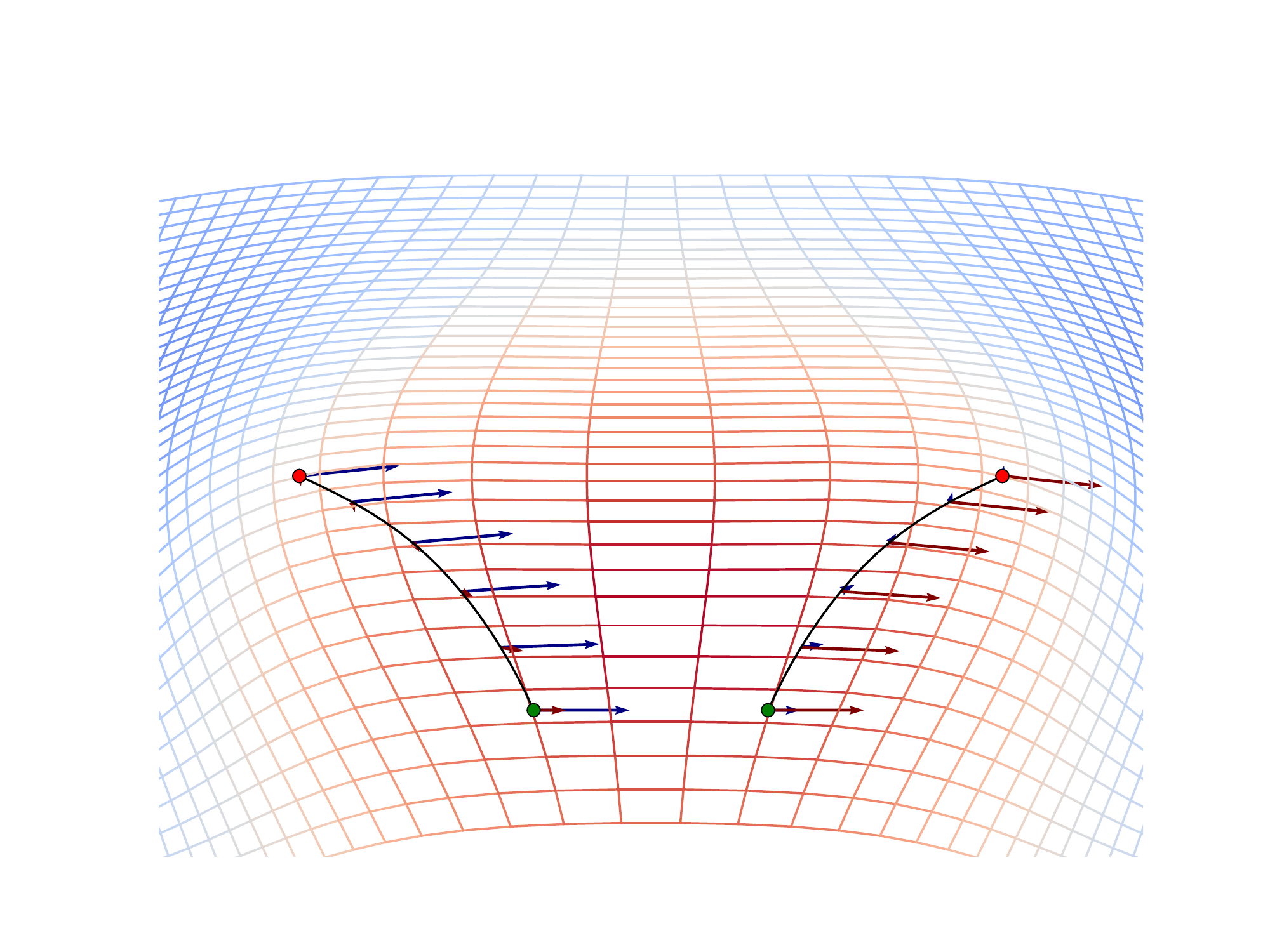}
      \includegraphics[width=.19\columnwidth,trim=80 60 80 60,clip=true]{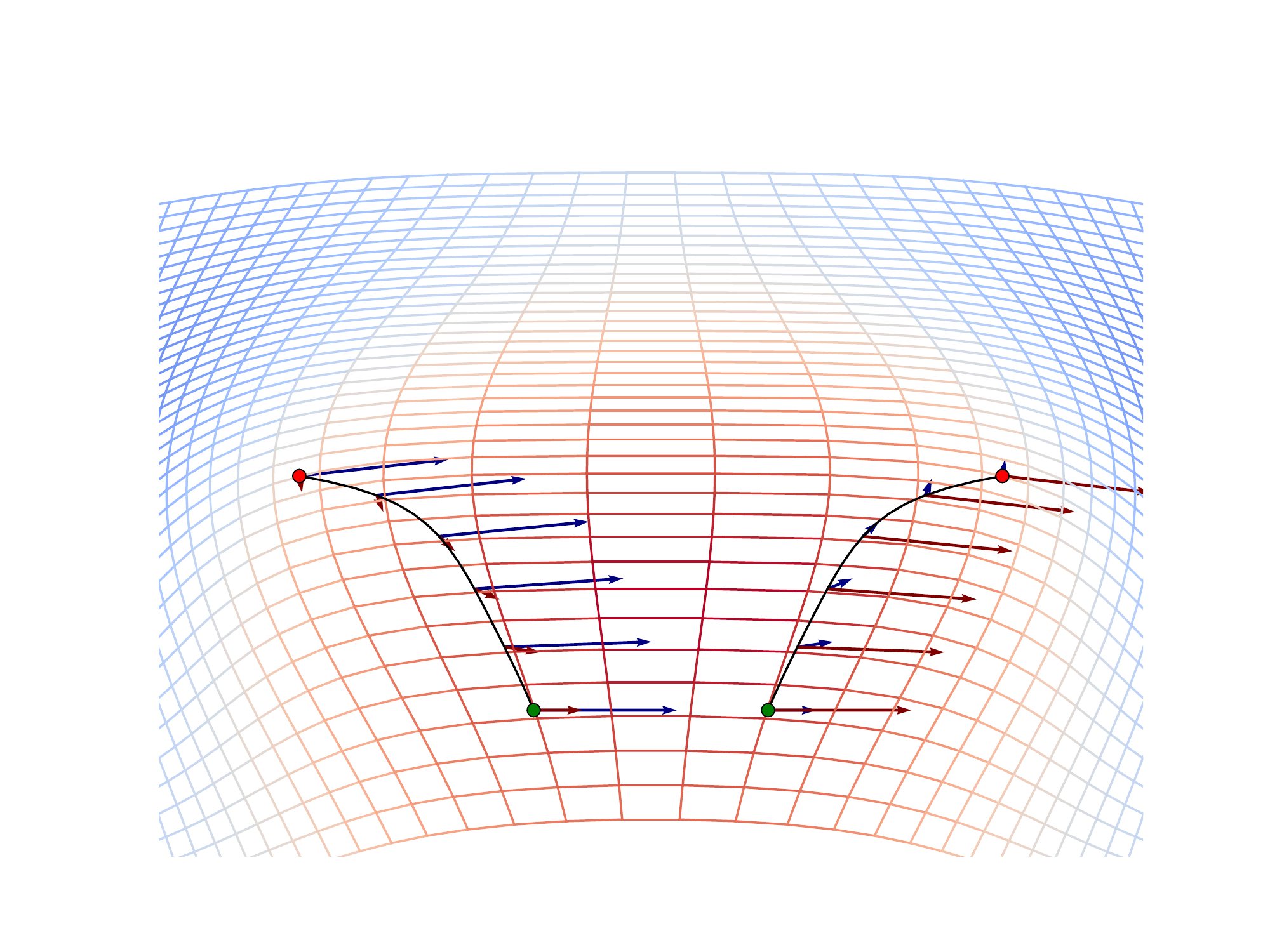}
      \includegraphics[width=.19\columnwidth,trim=80 60 80 60,clip=true]{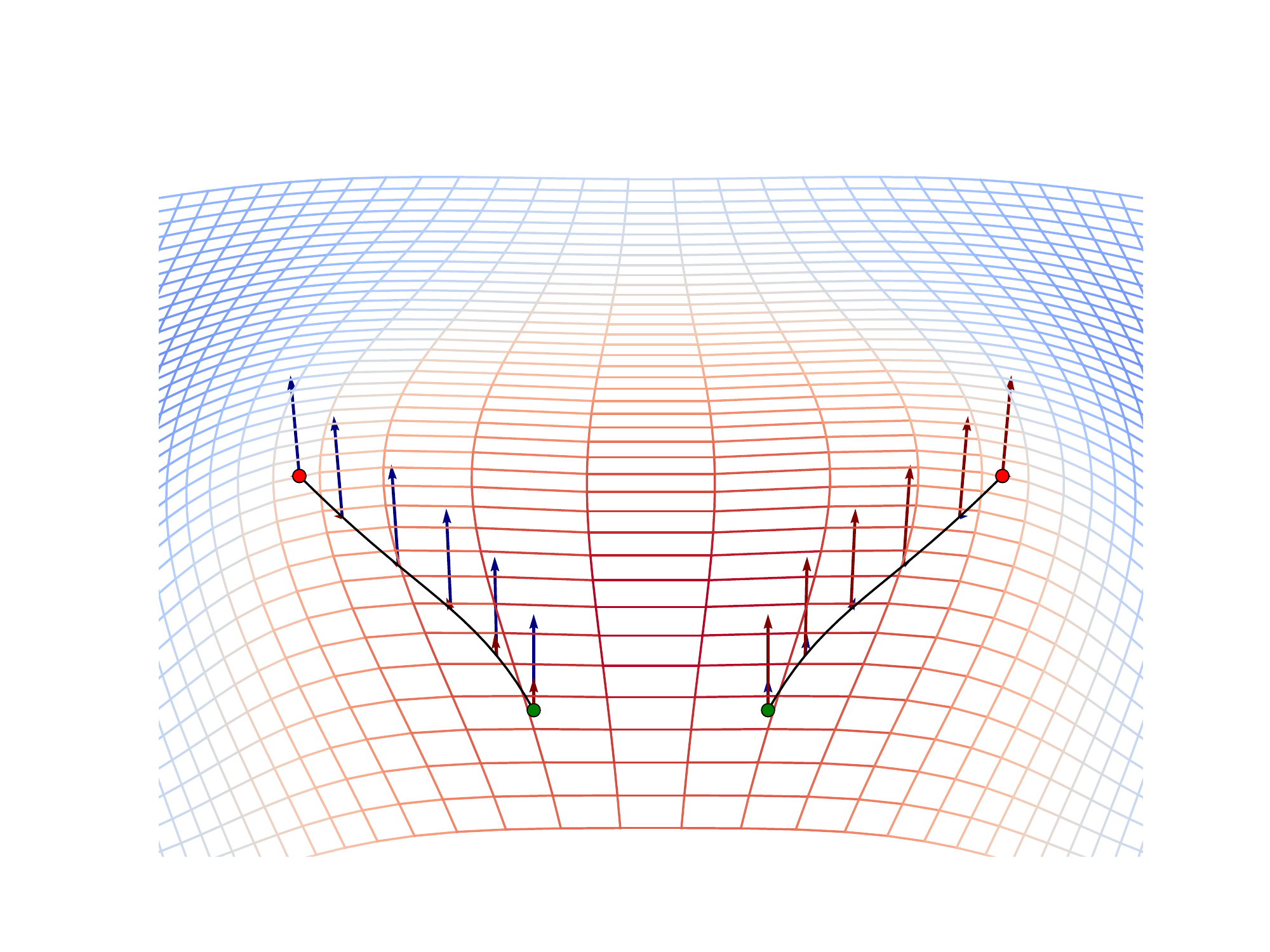}
      \includegraphics[width=.19\columnwidth,trim=80 60 80 60,clip=true]{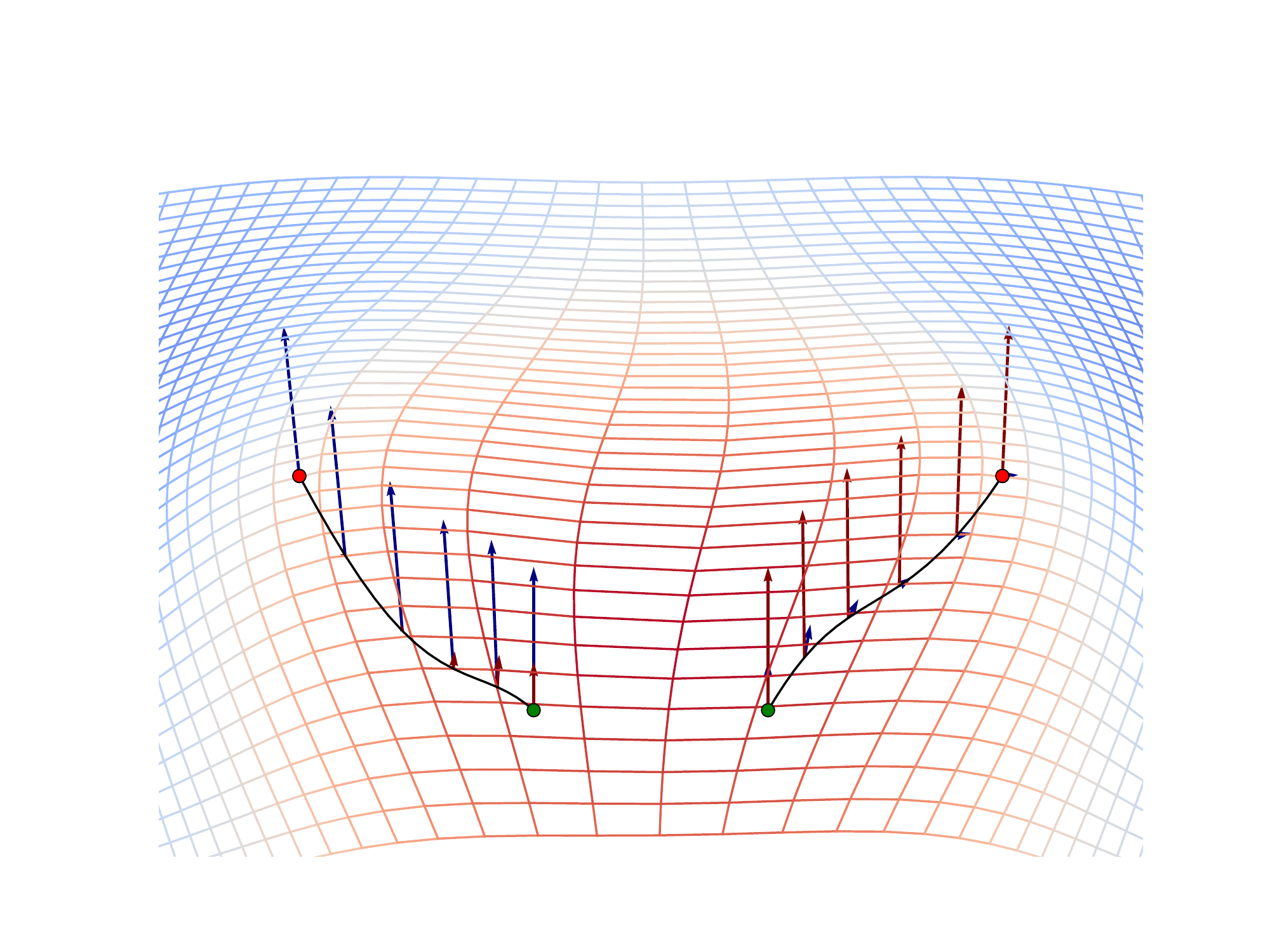}
      \\
      \hspace{.18cm}(a)
      \hspace{.163\columnwidth}(b)
      \hspace{.163\columnwidth}(c)
      \hspace{.163\columnwidth}(d)
      \hspace{.163\columnwidth}(e)
      \\
      \includegraphics[width=.19\columnwidth,trim=80 60 80 60,clip=true]{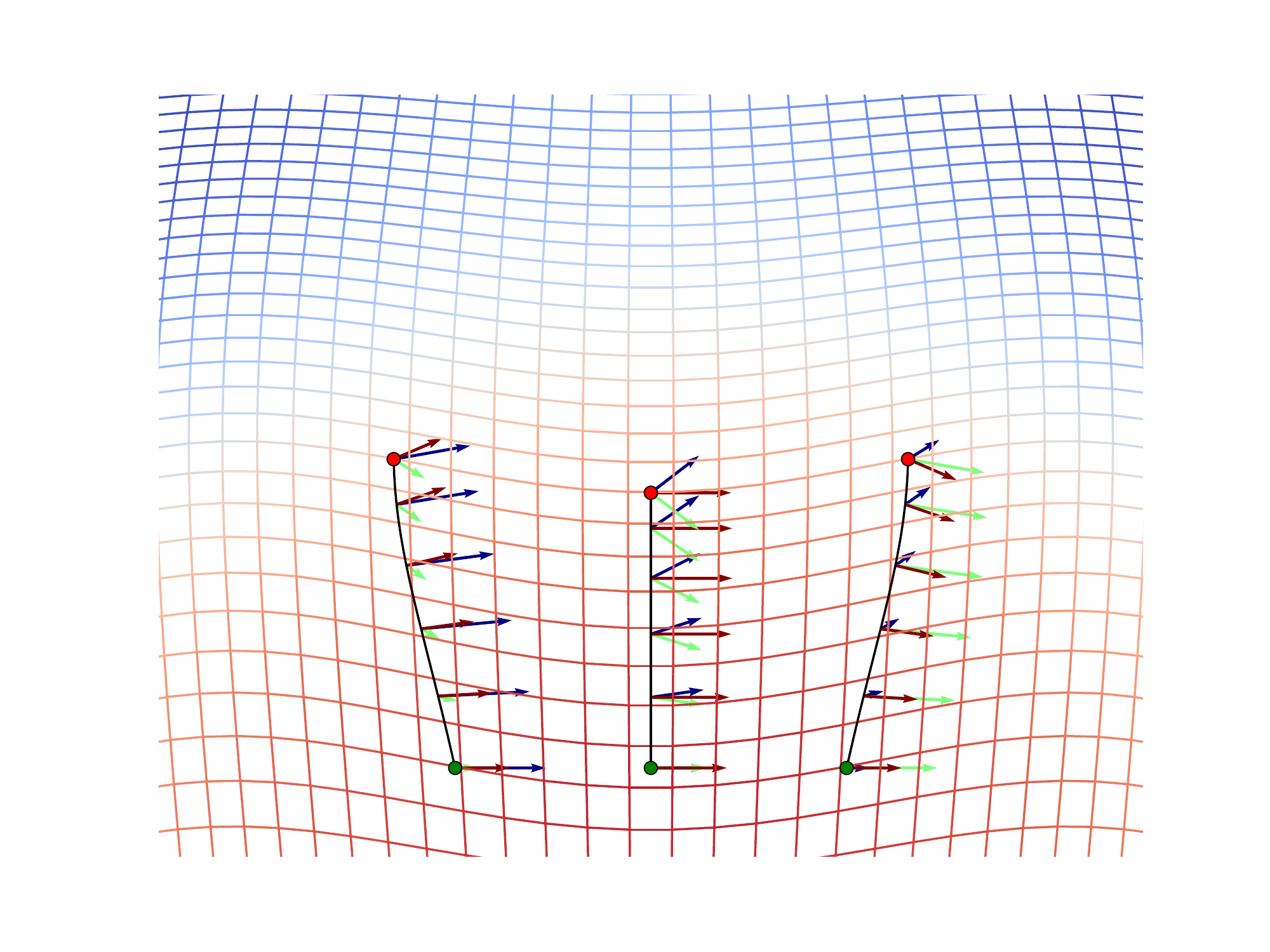}
      \includegraphics[width=.19\columnwidth,trim=80 60 80 60,clip=true]{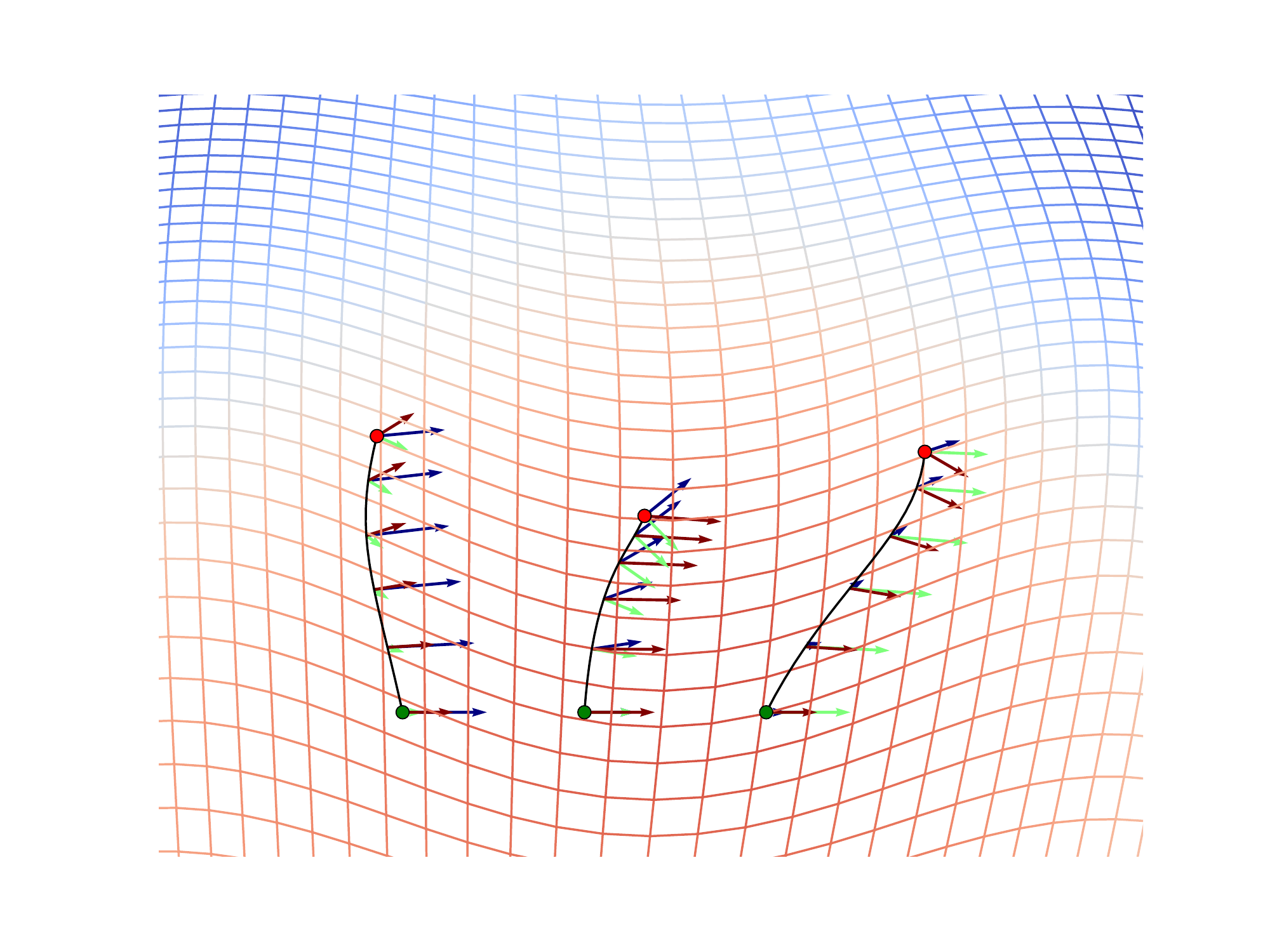}
      \includegraphics[width=.19\columnwidth,trim=80 60 80 60,clip=true]{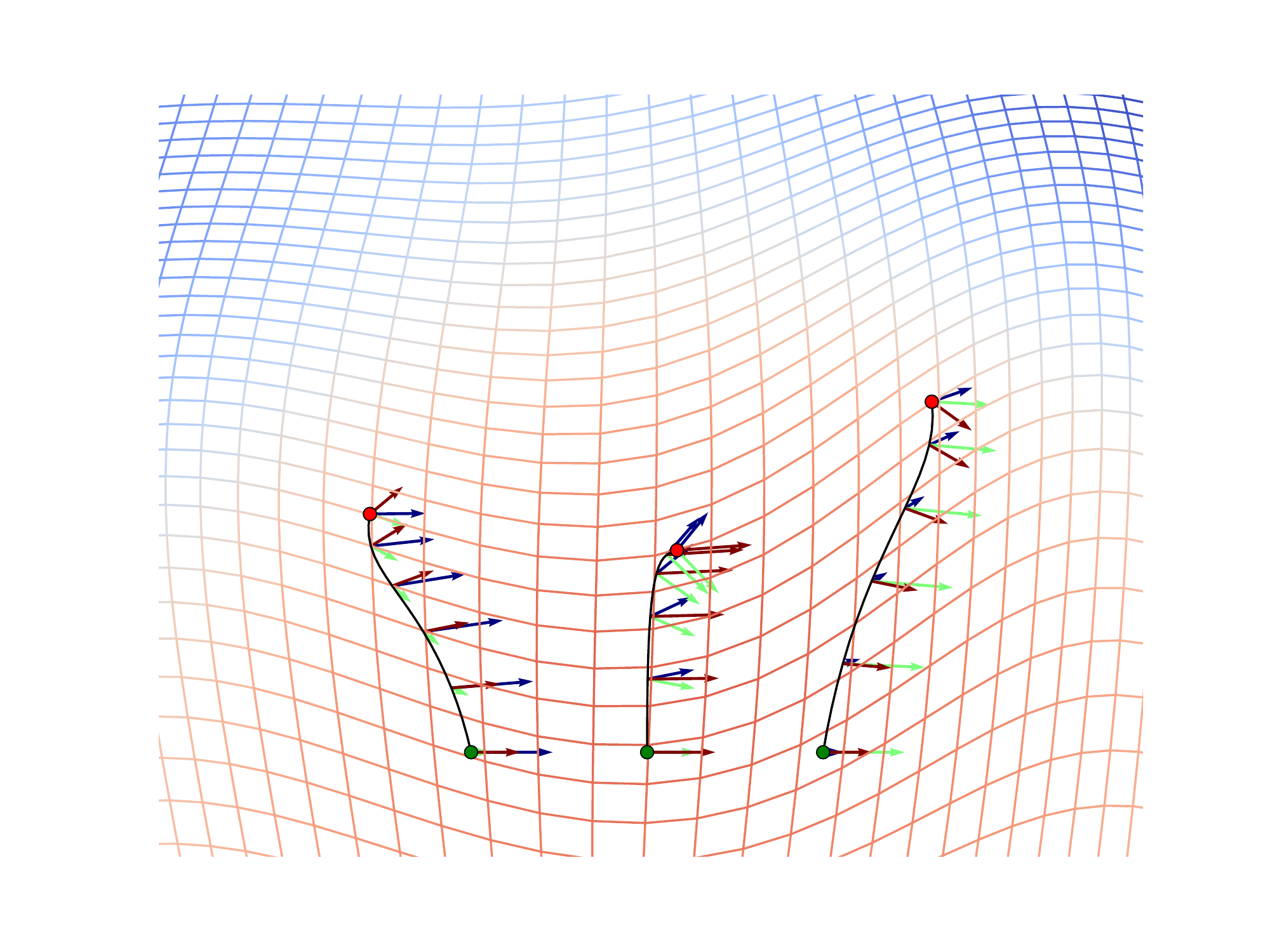}
      \includegraphics[width=.19\columnwidth,trim=80 60 80 60,clip=true]{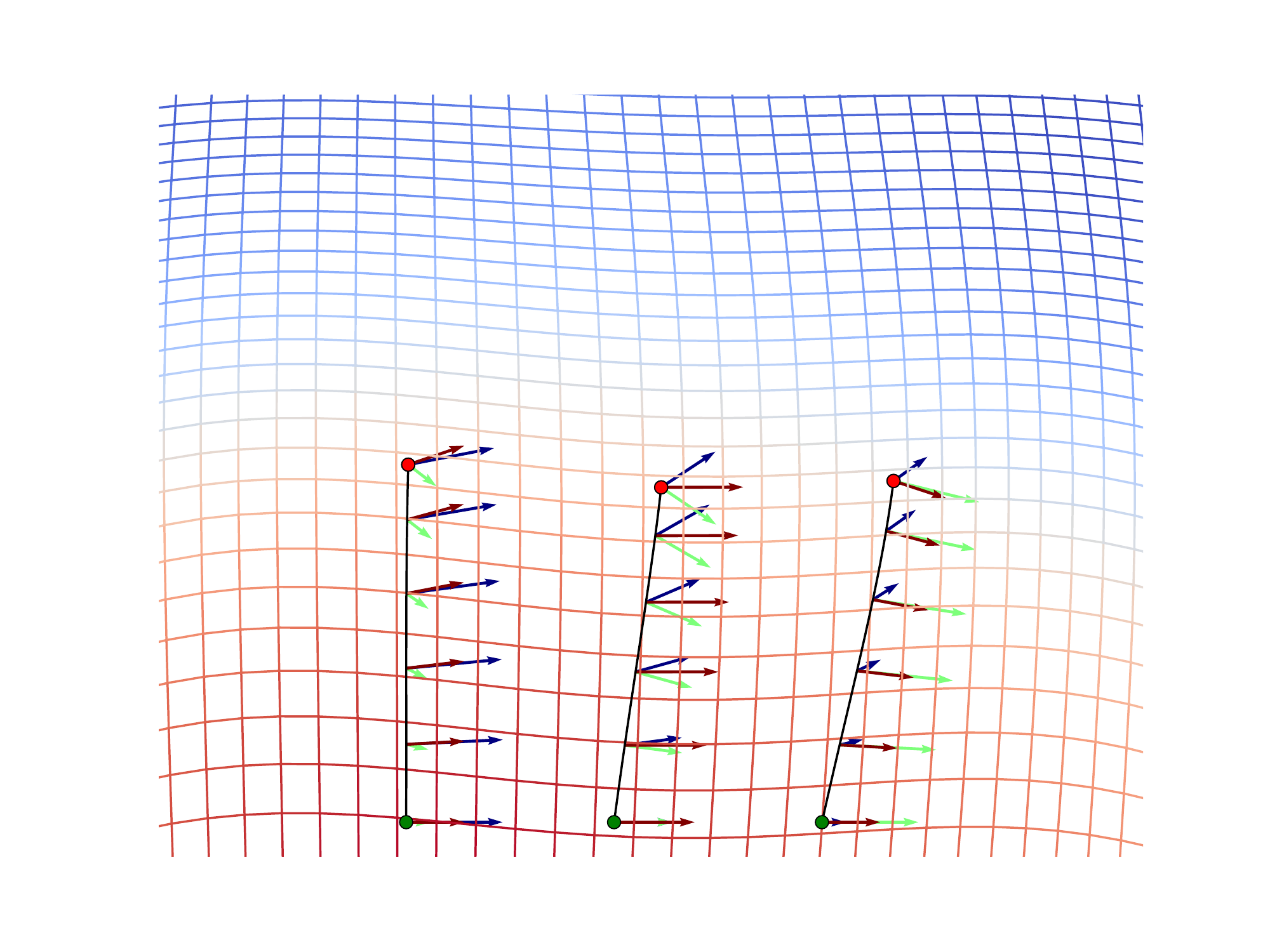}
      \includegraphics[width=.19\columnwidth,trim=80 60 80 60,clip=true]{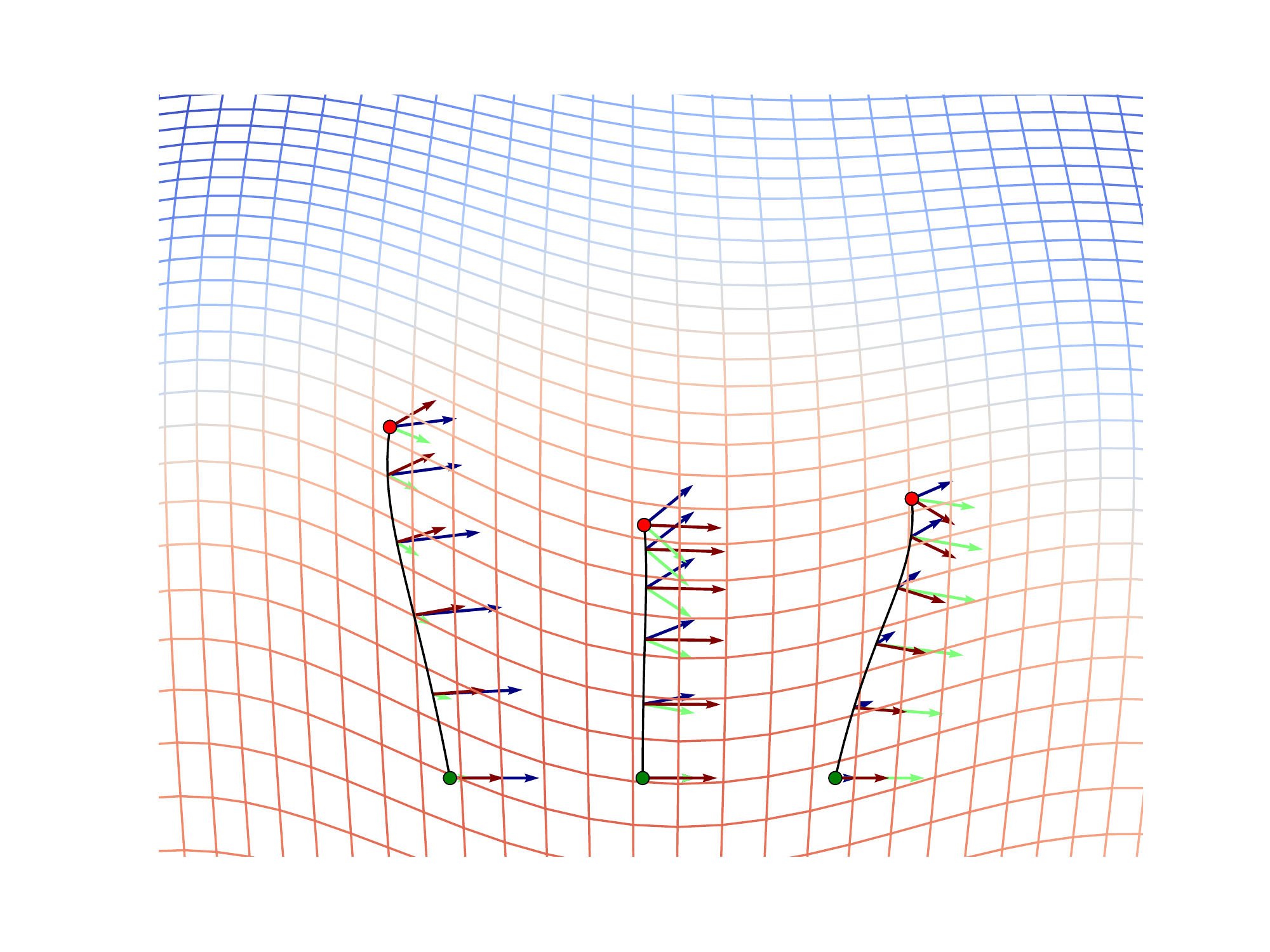}
  \end{center}
  \caption{
    (top row) Matching of two landmarks (green) to two landmarks (red) by
    (a) computing a minimizing Riemannian geodesic on the landmark manifold,
    and (b-e) minimizing MPPs with added covariance (arrows) in $\RR^2$ 
    horizontal direction
    (b-c) and vertical (d-e). The action of the
    corresponding diffeomorphisms on a regular grid is visualized by the
    deformed grid which is colored by the warp strain. The added covariance
    allows the paths to have more movement in the horizontal and vertical direction
    respectively
    because the anisotropically weighted metric penalizes high-covariance
    directions less. (bottom row) Five landmark trajectories with fixed initial
    velocity and anisotropic covariance but varying $V^*FM$ vertical initial
    momentum $\xi_0$. Changing the vertical momentum ``twists'' the paths.
    }
  \label{fig:landmarks}
\end{figure}

Let $\{p_1,\ldots,p_N\}$ be landmarks in a subset $\Omega\subset\RR^d$. The 
diffeomorphism group
$\mathrm{Diff}(\Omega)$ acts on the left on landmarks with the action
$\phi.\{p_1,\ldots,p_N\}=\{\phi(p_1),\ldots,\phi(p_N)\}$. In LDDMM, a Hilbert
space structure is imposed on a linear subspace $V$ of $L^2(\Omega,\RR^d)$ using
a self-adjoint operator $L:V\rightarrow V^*\subset L^2(\Omega,\RR^d)$ and
defining the inner product $\ip{\cdot,\cdot}_V$ by
\begin{equation*}
  \ip{v,w}_V=\ip{Lv,w}_{L^2}
  \ .
\end{equation*}
Under sufficient conditions on $L$, $V$ is reproducing and admits a kernel $K$ inverse
to $L$. $K$ is a Green's kernel when $L$ is a differential operator, or
$K$ can be a Gaussian kernel.
The Hilbert structure on $V$ gives a Riemannian metric on
a subset $G_V\subset\mathrm{Diff}(\Omega)$ by setting
$\|v\|_\phi^2=\|v\circ\phi^{-1}\|_V^2$, i.e. regarding $\ip{\cdot,\cdot}_V$ an inner
product on $T_{\mathrm{Id}}G_V$ and extending the metric to $G_V$ by
right-invariance. This Riemannian metric descends to a Riemannian metric on the
landmark space.

Let $M$ be the manifold
$M=\{(p_1^1,\ldots,p_1^d,\ldots,p_N^1,\ldots,p_N^d)|(p_i^1,\ldots,p_i^d)\in\RR^d\}$.
The LDDMM metric on the landmark manifold $M$ is directly related to the kernel
$K$ when written as a cometric 
$g_p(\xi,\eta)=\sum_{i,j=1}^N \xi^iK(p_i,p_j)\eta^j$. Letting
$i^k$ denote the index of the $k$th component of the $i$th landmark, the
cometric is in coordinates 
$g_p^{i^kj^l}=K(p_i,p_j)_k^l$.
The Christoffel symbols can be written in terms of derivatives of the cometric
$g^{ij}$ \cite{micheli_differential_2008} (recall that $\delta^i_j=g^{ik}g_{kj}=g_{jk}g^{ki}$)
\begin{equation}
  \Gamma\indices{^k_{ij}}
  =
  \frac{1}{2}
  g_{ir}
  \left(
  g^{kl}g\indices{^{rs}_{,l}}
  -
  g^{s l}g\indices{^{rk}_{,l}}
  -
  g^{r l}g\indices{^{ks}_{,l}}
  \right)
  g_{sj}
  \ .
  \label{eq:mario}
\end{equation}
This relation comes from the fact that
$g_{jm,k}
  =
  -g_{jr}g\indices{^{rs}_{,k}}g_{sm}
$ gives the derivative of the metric.
The derivatives of the cometric is simply
$g\indices{^{i^kj^l}_{,r^q}}=(\delta_{r}^{i}+\delta_{r}^{j})\partial_{p_r^q} K(p_i,p_j)_k^l$.
Using \eqref{eq:mario}, derivatives of the Christoffel symbols can be computed
\begin{align*}
  \Gamma\indices{^k_{ij}_{,\xi}}
  &=
  \frac{1}{2}
  g_{ir,\xi}
  \left(
  g^{kl}g\indices{^{rs}_{,l}}
  -
  g^{sl}g\indices{^{rk}_{,l}}
  -
  g^{rl}g\indices{^{ks}_{,l}}
  \right)
  g_{sj}
  +\frac{1}{2}
  g_{ir}
  \left(
  g^{kl}g\indices{^{rs}_{,l}}
  -
  g^{s l}g\indices{^{rk}_{,l}}
  -
  g^{r l}g\indices{^{ks}_{,l}}
  \right)
  g_{sj,\xi}
  \\
  &\qquad
  +\frac{1}{2}
  g_{ir}
  \left(
  g\indices{^{kl}_{,\xi}}g\indices{^{rs}_{,l}}
  +
  g^{kl}g\indices{^{rs}_{,l\xi}}
  -
  g\indices{^{sl}_{,\xi}}g\indices{^{rk}_{,l}}
  -
  g^{sl}g\indices{^{rk}_{,l\xi}}
  -
  g\indices{^{rl}_{,\xi}}g\indices{^{ks}_{,l}}
  -
  g^{rl}g\indices{^{ks}_{,l\xi}}
  \right)
  g_{sj}
   \ .
\end{align*}
This provides the full data for numerical integration of the evolution equations
on $F^kM$. 

In Figure~\ref{fig:landmarks} top row, we plot minimizing normal MPPs on the landmark
manifold with two landmarks and varying covariance in the $\RR^2$ horizontal
and vertical direction. The plot shows the landmark equivalent of experiment in
Figure~\ref{fig:varlog}. 
Note how adding covariance in the horizontal and
vertical direction respectively allows the minimizing normal MPP to varying more in
these directions because the anisotropically weighted metric penalizes high-covariance
directions less.

Figure~\ref{fig:landmarks} bottom row shows five curves satisfying the MPP
equations with varying vertical $V^*FM$ initial momentum similarly to the plots in
Figure~\ref{fig:varxia}. Again, we see how the extra degrees of freedom allows
the paths to twist generating a higher-dimensional family than classical
geodesics with respect to $\Cc$.

\section{Discussion and Concluding Remarks}
Incorporating anisotropy in models for data 
in non-linear spaces via the frame bundle as pursued in this paper leads to 
a sub-Riemannian structure and metric.
A direct implication is that most probable paths to observed data in the sense
of sequences of stochastic steps of a driving semi-martingale are not 
related to geodesics in the classical sense. Instead, a best estimate of the sequence of steps
$w_t\in\RR^d$ that leads to an observation $x=\phi_{u}(w_t)|_{t=1}$ is an MPP in the sense of
Definition~\ref{def:mpp}. As shown in the paper, these paths are generally not
geodesics or polynomials with respect to the connection on the manifold. In particular, if $M$
has a Riemannian structure, the MPPs are generally neither Riemannian geodesics
or Riemannian polynomials.
Below, we discuss statistical implications of this result.

\subsection{Statistical Estimators}
Metric distances and Riemannian geodesics have been the traditional
vehicle for representing observed data in non-linear spaces. Most fundamentally,
the sample Frech\'et mean
\begin{equation}
  \hat{x}=\textrm{argmin}_{x \in M} \sum_{i=1}^Nd_{g_R}\left(x, x_i\right)^2
  \label{eq:frechet}
\end{equation}
of observed data $x_1,\ldots,x_N\in M$ relies crucially on the Riemannian
distance $d_{g_R}$ connected to the metric $g_R$. Many PCA constructs, e.g. Principal Geodesics Analysis
\cite{fletcher_principal_2004-1}, uses the Riemannian $\mathrm{Exp}$ and $\Log$ maps to
map between linear tangent spaces and the manifold. These maps are
defined from the Riemannian metric and Riemannian geodesics. Distributions
modelled as in the random orbit model \cite{miller_statistical_1997} or Bayesian models
\cite{zhang_probabilistic_2013,zhang_bayesian_2013} again
rely on geodesics with random initial conditions.

Using the frame bundle sub-Riemannian metric $g_{FM}$, we can define an estimator 
analogous to the Riemannian Frech\'et mean estimator. Assuming the covariance
is a priori known, the estimator
\begin{equation}\label{minimize}
  \hat{x}=\textrm{argmin}_{u \in s(M)} \sum_{i=1}^Nd_{FM}\left(u, \pi^{-1}(x_i)\right)^2
\end{equation}
acts correspondingly to the Frech\'et mean estimator \eqref{eq:frechet}.
Here $s\in\Gamma(FM)$ is a (local) section of $FM$ that to $x\in M$ connects 
the known covariance represented by $s(x)\in FM$. The distances 
$d_{FM}\left(u, \pi^{-1}(x_i)\right)$, $u=s(x)$ are realized by MPPs from the
mean candidate $x$ to the fibers $\pi^{-1}(x_i)$. The Frech\'et mean problem
is thus lifted to the frame bundle with the anisotropic weighting incorporated in
the metric $g_{FM}$. This metric is not related to $g_R$ except for its 
dependence on the connection $\mathcal C$ that can be defined as the Levi-Civita
connection of $g_R$. The fundamental role of the distance $d_{g_R}$ and
$g_R$ geodesics in \eqref{eq:frechet} is thus removed.

Because covariance is an integral part of the model, sample covariance can also be 
estimated directly along with the sample mean. In \cite{sommer_modelling_2016}, 
the estimator
\begin{equation}
  \hat{u}=\textrm{argmin}_{u \in FM} \sum_{i=1}^Nd_{FM}\left(u, \pi^{-1}(x_i)\right)^2
  -N\log(\mathrm{det}_{g_R}u)
  \label{eq:frefm}
\end{equation}
is suggested. The normalizing term $-N\log(\mathrm{det}_{g_R}u)$ is derived such that
the estimator exactly corresponds to the maximum likelihood estimator of mean
and covariance for Euclidean Gaussian distributions. The determinant is defined
via $g_R$, and the term acts to prevent the covariance
from approaching infinity. Maximum likelihood estimators of mean and covariance for 
normally distributed Euclidean data have unique solutions in the sample mean and
sample covariance matrix, respectively. Uniqueness of the
Frech\'et mean \eqref{eq:frechet} is only ensured for sufficiently concentrated
data. For the estimator \eqref{eq:frefm}, existence and uniqueness properties are not
immediate, and more work is needed in order to find necessary and sufficient
conditions.

\subsection{Priors and Low-Rank Estimation}
The low-rank cometric formulation pursued in Section~\ref{sec:rankdef} gives a
natural restriction of \eqref{eq:frefm} to $u\in F^kM$, $1\le k\le d$. As for
Euclidean PCA, most variance is often captured in the span of the first
$k$ eigenvectors with $k\ll d$. Estimates of the remaining eigenvectors are
generally ignored as the variance of the eigenvector estimates increases as the noise 
captured in the span of the last eigenvectors becomes increasingly uniform. The 
low-rank cometric restricts the estimation to only the first $k$ eigenvectors
and thus builds the construction directly
into the model. In addition, it makes numerical implementation feasible because
a numerical representation need only store and evolve $d\times k$ matrices.
As a different approach for regularizing the estimator \eqref{eq:frefm}, the
normalizing term $-N\log(\mathrm{det}_{g_R}u)$ can be extended with other priors,
e.g. an $L^1$-type penalizing term. Such priors can potentially partly remove existence and
uniqueness issues and result in additional sparsity properties that can
benefit numerical implementations. The effects of such priors has yet to be
investigated.

In the $k=d$ case, the number of degrees of freedoms for the MPPs grows
quadratically in the dimension $d$. This naturally increases the variance of any
MPP estimate given only one sample from its trajectory. The low-rank
cometric formulation reduces the growth to linear in $d$. The number of degrees
of freedom is however still $k$ times larger than for Riemannian geodesics.
With longitudinal data, more samples per trajectory can be obtained reducing the
variance and allowing a better estimate of the MPP. However, for the estimators
\eqref{minimize} and \eqref{eq:frefm}
above, estimates of the actual optimal MPPs are not needed, only their squared length. It can be
hypothesized that the variance of the length estimates is lower than the variance of
the estimates of the corresponding MPPs. Further investigation regarding this will be the 
subject of future work.

\subsection{Conclusion}
The underlying model of anisotropy used in this paper originates from the anisotropic normal
distributions formulated in \cite{sommer_anisotropic_2015} and the
diffusion PCA framework \cite{sommer_diffusion_2014}. Because many
statistical models are defined using normal distributions, this approach for
incorporating anisotropy extends to
models such as linear regression. We expect that finding most probable 
paths in other statistical models such as regressions models can be carried out with a
program similar to the program presented in this paper.

The difference between MPPs and geodesics shows that the geometric and metric
properties of geodesics, zero acceleration and local distance minimization, are
not directly related to statistical properties such as maximizing path probability.
Whereas the concrete application and model determines if metric or statistical
properties are fundamental, most statistical models are formulated without
referring to metric properties of the underlying space. It can therefore be
argued that the direct incorporation of anisotropy and the resulting MPPs are 
natural in the context of many models of data variation in non-liner spaces.

\section*{Acknowledgements}
The author wishes to thank Peter W. Michor and Sarang Joshi for suggestions for the geometric
interpretation of the sub-Riemannian metric on $FM$ and discussions on diffusion
processes on manifolds. 
  The work was supported by the Danish Council for Independent Research,
the CSGB Centre for Stochastic Geometry and Advanced Bioimaging funded by a grant from the Villum foundation,
and the Erwin Schr\"odinger Institute in Vienna.

\bibliographystyle{amsalpha}
\bibliography{ss}   



\end{document}